\documentclass[reqno]{amsart}
\usepackage{amsmath}
\usepackage{amsthm}
\usepackage{amssymb}
\usepackage{stmaryrd}
\usepackage{afterpage}
\usepackage{rotating}    %sidewaystable
\usepackage{adjustbox}   
\usepackage{booktabs}
\usepackage[utf8]{inputenc}
\usepackage[T1]{fontenc}
\usepackage{lmodern}
\usepackage{eucal}
\pdfoutput=1 %Needed for microtype and ArXiV.
\usepackage{microtype}
\usepackage{url}
\usepackage{enumitem}
\usepackage{tikz}
\usetikzlibrary{arrows.meta,calc}
\usepackage{tikz-cd}
\usepackage{placeins}
\usepackage{bbm}
\usepackage{longtable,booktabs}
\usepackage{todonotes}

\usepackage[letterpaper, left = 1.5in, right = 1.5in, top = 1.5in, bottom = 1.5in]{geometry}

\usepackage{pdflscape} 

\usepackage{hyperref}

\usepackage{xcolor}
\hypersetup{
    colorlinks,
    linkcolor={red!50!black},
    citecolor={blue!50!black},
    urlcolor={blue!80!black}
}
\definecolor{chartgray}{gray}{0.4}
\definecolor{darkcyan}{rgb}{0, 0.7, 0.7}
\definecolor{darkgreen}{rgb}{0, 0.65, 0}
\definecolor{truemagenta}{rgb}{1, 0, 1}

\usepackage[capitalise,nameinlink]{cleveref}
\crefname{subsection}{Subsection}{Subsections}
\crefname{subsubsection}{Subsubsection}{Subsubsections}

\usepackage{multirow}

\usepackage{graphicx}% Required for inserting images

\usepackage{lipsum}     
\usepackage{booktabs} 
\usepackage{rotating}
\usepackage{tabularx}

\theoremstyle{definition}
\newtheorem{theorem}{Theorem}[section]

\newtheorem{defn}[theorem]{Definition}

\newtheorem{cor}[theorem]{Corollary}
\newtheorem{lemma}[theorem]{Lemma}
\newtheorem{prop}[theorem]{Proposition}
\newtheorem{rmk}[theorem]{Remark}

\newtheorem{conjecture}[theorem]{Conjecture}

\newtheorem*{rmk*}{Remark}
\newtheorem*{ex*}{Example}
\newtheorem*{theorem*}{Theorem}
\newtheorem*{defn*}{Definition}
\newtheorem{thmx}{Theorem}
\renewcommand{\thethmx}{\Alph{thmx}} 

\newcommand{\bbZ}{\mathbb{Z}}

\newcommand{\bbF}{\mathbb{F}}
\newcommand{\bbN}{\mathbb{N}}
\newcommand{\bbP}{\mathbb{P}}

\newcommand{\bbR}{\mathbb{R}}

\newcommand{\Fun}{\mathrm{Fun}}
\newcommand{\s}{\mathrm{s}}

\newcommand{\Hom}{\operatorname{Hom}}

\newcommand{\Ext}{\operatorname{Ext}}

\newcommand{\colim}{\operatorname*{colim}}

\renewcommand{\ker}{\operatorname{Ker}}

\newcommand{\Sq}{\mathrm{Sq}}

\newcommand{\id}{\operatorname{id}}

\newcommand{\tmf}{\operatorname{tmf}}
\newcommand{\CP}{\mathbb{C}\mathbb{P}}
\newcommand{\KO}{\operatorname{KO}}
\newcommand{\ko}{\operatorname{ko}}

\newcommand{\h}{{\mathrm{h}}}

\makeatletter
\DeclareRobustCommand{\tvdots}{%
  \vbox{\baselineskip4\p@\lineskiplimit\bbZ@\kern0\p@\hbox{.}\hbox{.}\hbox{.}}}
\makeatother

\makeatletter
\@namedef{subjclassname@2020}{\textup{2020} Mathematics Subject Classification}
\makeatother

\makeatletter
\newcommand{\raisemath}[1]{\mathpalette{\raisem@th{#1}}}
\newcommand{\raisem@th}[3]{\raisebox{#1}{$#2#3$}}
\makeatother

\makeatletter
\def\l@subsection{\@tocline{2}{0pt}{2pc}{5pc}{}}
\makeatother

\newcommand{\CPfourgrid}[4]{%
  \foreach \x in {#1,...,#2} {
    \draw[black!20] (\x-0.5,#3-0.5) -- (\x-0.5,#4-0.5);
  }
  \foreach \y in {#3,...,#4} {
    \draw[black!20] (#1-0.5,\y-0.5) -- (#2-0.5,\y-0.5);
  }
}

\newcommand{\CPfourarrow}[1][]{%
  \draw[-{Stealth[length=1.5mm]},shorten >=0.08cm,#1]}

\tikzset{CPfour red/.style={red!85!black}}

\begin{document}

\title[Exotic smooth structures on complex projective spaces]{Exotic smooth structures\\ on complex projective spaces} \author{Yueshi Hou} \address{Department of Mathematics, University of California San Diego, La Jolla, CA 92093, USA} \email{yuh091@ucsd.edu}

\author{Hana Jia Kong} \address{School of Mathematical Sciences, Zhejiang University, Hangzhou, Zhejiang 310058, China} \email{hana.jia.kong@gmail.com}

\author{Zhouli Xu} \address{Department of Mathematics, University of California Los Angeles, Los Angeles, CA 90095-1555, USA} \email{xuzhouli@ucla.edu}

\begin{abstract}
We study smooth structures on the topological manifold underlying complex projective space. We conjecture that every complex projective space $\CP^n$ with $n>3$ admits an exotic smooth structure.

We prove two results providing strong evidence for this conjecture. First, the set of complex dimensions $n$ for which $\CP^n$ admits an exotic smooth structure has natural density one. Second, we explicitly identify $266$ congruence classes modulo $288$ such that $\CP^n$ admits an exotic smooth structure whenever $n>3$ lies in one of these classes. These classes cover more than $92\%$ of complex dimensions, improving the previously known lower bound of approximately $36\%$.

We prove these results by first relating the concordance group of $\CP^n$ to suitable Atiyah--Hirzebruch filtration subgroups of its stable cohomotopy group, and then applying computations in chromatic homotopy theory at height $2$.
\end{abstract}

\maketitle

\tableofcontents

\section{Introduction}

We study smooth structures on the topological manifold underlying complex projective space. An \emph{exotic smooth structure on $\CP^n$} is a smooth manifold homeomorphic, but not diffeomorphic, to standard complex projective space.

\begingroup
\renewcommand{\thetheorem}{A}
\begin{conjecture}\label{conj:exotic-cpn}
For every integer $n\geq 4$, $\CP^n$ admits an exotic smooth structure.
\end{conjecture}
\endgroup

Historically, topological manifolds that have the same homotopy type as $\CP^n$ have been extensively studied. For each $n \geq 4$, Hsiang \cite{Hsiang1966homotopyCP} constructed infinitely many smooth manifolds that are homotopy equivalent to the standard $\CP^n$, with distinct rational Pontrjagin classes. Thus these smooth manifolds are not homeomorphic to the standard $\CP^n$. On the other hand, for each $n \geq 4$, Sullivan \cite{Sullivan1996triangulating} completely classified PL-homeomorphism classes of homotopy $\CP^n$.
% , while
% Brumfiel studied homotopy smoothings of almost smooth manifolds and their
% relation to free circle actions on homotopy spheres
% \cite[Appendix I]{Brumfiel1971almostsmooth}.
% These results concern a broader homotopy type; detecting a new smoothing
% of the fixed topological manifold $\CP^n$ requires the corresponding
% smoothing-theoretic information.

% Final review (2026-11-04), original line 458.

% When it comes to finding exotic $\CP^n$, connected sums $\CP^n\mathbin{\#}\Sigma^{2n}$ on standard $\CP^n$ give natural candidates , where $\Sigma^{2n}$ is a homotopy sphere.
% Even when $\Sigma^{2n}$ is exotic, the resulting manifold need not be
% exotic. The inertia group $I(\CP^n)$ records the oriented homotopy
% spheres for which this connected sum is orientation-preservingly
% diffeomorphic to $\CP^n$.
% Kawakubo proved that $I(\CP^n)=0$ for $n\leq 8$
% \cite{Kawakubo1968inertia}.
% In contrast, Basu--Kasilingam proved that the inertia groups of $\CP^9$
% and $\CP^{13}$ are nontrivial
% \cite{BasuKasilingam2018inertia}.
% Consequently, finding a homotopy sphere and detecting its contribution
% to the smooth structure of projective space are distinct tasks.

Regarding the topological manifold $\CP^n$, Farrell--Jones \cite{FarrellJones1994hyperbolic} proved that $\CP^{4n+1}$ admits exotic smooth structures for every $n\geq 1$ by showing there is a homotopy sphere $\Sigma^{8n+2}$ such that $\CP^{4n+1}\mathbin{\#}\Sigma^{8n+2}$ is not diffeomorphic to the standard smooth structure. Kasilingam \cite{Kasilingam2015classification} showed that $\CP^3$ has a unique smooth structure, and that for $4\leq n\leq8$, $\CP^n$ admits exotic smooth structures \cite{Kasilingam2026classification}. More recently, Bauer--Quigley \cite{BauerQuigley2026freeactions} constructed infinite families of homotopy spheres that do not bound parallelizable manifolds and admit smooth free $S^1$-actions. An immediate consequence of their computations is that there exist exotic $\CP^{6+36n}$ \cite{BauerQuigley2026freeactions}. Using $EO$-theory computations, Basu--Biswas \cite{BasuBiswas2026exotic} further extended the range of dimensions in which complex projective spaces admit exotic smooth structures, with examples including $\CP^{10+36n}$, $\CP^{11+36n}$, and $\CP^{20+36n}$.

% Final review (2026-11-04), original line 483.

% Final review (2026-11-04), original line 483.
%\todo[inline]{\textbf{Review -- Reference.} The key \texttt{BasuBiswas2026exotic} is missing from \texttt{refs.bib}. Add Basu--Biswas, \emph{Exotic structures on hyperbolic manifolds via the EO-theory of Projective Spaces}, arXiv:2610.10070. Its inertia-group result directly covers $10+36n$ and $20+36n$. For $11+36n$, make explicit the deduction from its free-circle-action construction through a nonzero class in $[\CP^{11+36n},Top/O]$; this is not one of the dimensions in its inertia-group statement.}
% Review source check: Basu--Biswas, 2610.10070.pdf, Theorems B and C, PDF p. 3 / printed p. 3, and Section 5.15, PDF pp. 44--45 / printed pp. 44--45.
% Primary URL: https://arxiv.org/pdf/2610.10070

Our first main result provides asymptotic evidence for Conjecture~\ref{conj:exotic-cpn}.

\setcounter{thmx}{1}
\begin{thmx}\label{thm:unus}
The set of positive integers $n$ for which $\CP^n$ admits an exotic smooth structure has natural density one. Equivalently,
\[
\lim_{N\to\infty}
\frac{\#\{1\leq n\leq N:
            \CP^n\text{ admits an exotic smooth structure}\}}{N}
=1.
\]
\end{thmx}

Aside from the density-one result, using chromatic homotopy theory at height two and at the primes $2$ and $3$, we show that $\CP^n$ admits an exotic smooth structure for every $n\geq4$ in one of $266$ explicit congruence classes modulo $288$, covering over $92\%$ of complex dimensions.

Let $[a,b]$ denote the set of integers from $a$ to $b$, including $a$ and $b$, and define
\[
\begin{aligned}
\mathcal R_2={}&[0,13]\cup\{15,17\}\cup[19,25]\cup[27,29]\cup[31,37]\cup[39,41]\\
&{}\cup[43,47]\cup\{49\}\cup[51,61]\cup[65,67]\cup[69,73]\\
&{}\cup[75,77]\cup[81,83]\cup\{85\}\cup[87,95],
\end{aligned}
\]
\[
\mathcal R_3=\{0\}\cup[5,6]\cup[13,14]\cup[19,24]\cup[37,38]\cup[41,50]\cup[53,71].
\]
Define also
\[
\mathcal R
=\{0\leq r<288:r\bmod96\in\mathcal R_2
                 \text{ or }r\bmod72\in\mathcal R_3\}.
\]
For the size of these sets, we have
\[
|\mathcal R_2|=77,\qquad |\mathcal R_3|=42,\qquad |\mathcal R|=266.
\]
Thus $\mathcal R$ contains a proportion $266/288=133/144>92\%$ of the residue classes modulo $288$.

\begin{thmx}\label{main theorem prime 2 and 3}

Let $n\geq4$. Then $\CP^n$ admits an exotic smooth structure if either
\begin{enumerate}[label=(\roman*),ref=\thethmx(\roman*),leftmargin=*,nosep]
\item\label[thmx]{thm:duo} $n\bmod96\in\mathcal R_2$, or
\item\label[thmx]{thm:tres} $n\bmod72\in\mathcal R_3$.
\end{enumerate}
Thus the conclusion holds whenever $n\bmod288\in\mathcal R$.
\end{thmx}

\begin{rmk}
Define
\[
\begin{aligned}
\mathcal R_2'
&=\{r\in[0,95]:r,r+96,r+192\in\mathcal R\}\\
&=\mathcal R_2\cup\{14,38,48,62,68,86\}.
\end{aligned}
\]
Then $|\mathcal R_2'|=83$, and Theorem~\ref{main theorem prime 2 and 3} yields exotic smooth structures on $\CP^n$ for every $n\geq4$ and $n\bmod96\in\mathcal R_2'$. The complement consists of the following $13$ residues:
\[
\begin{aligned}
\overline{\mathcal R_2'}&=[0,95]\setminus\mathcal R_2'\\
&=\{16,18,26,30,42,50,63,64,74,78,79,80,84\}.
\end{aligned}
\]

\par\medskip
\noindent\begin{minipage}{\linewidth}
\centering
\begingroup
\definecolor{coveragegrey}{RGB}{180,180,180}
\tikzset{
  coverage labels/.style={ font=\fontsize{6}{7}\selectfont, inner sep=0pt
  }
}

\begin{minipage}[t]{.34\linewidth}
\centering
\vspace{0pt} {\scriptsize $\mathcal R_2'\subseteq\mathbb Z/96$\par $83/96\approx86.46\%$\par
}
\smallskip
\begin{tikzpicture}[x=3mm,y=3mm]
  \fill[coveragegrey] (0,0) rectangle (12,-8);

  \foreach \r in {16,18,26,30,42,50,63,64,74,78,79,80,84} {
    \pgfmathtruncatemacro{\col}{mod(\r,12)}
    \pgfmathtruncatemacro{\row}{floor(\r/12)}
    \fill[white] (\col,-\row) rectangle ++(1,-1);
  }

  \draw[black!45,line width=.3pt,step=1]
    (0,0) grid (12,-8);

  \foreach \col in {0,...,11}
    \node[coverage labels,anchor=south]
      at (\col+.5,.2) {\col};

  \foreach \row in {0,...,7} {
    \pgfmathtruncatemacro{\base}{12*\row}
    \node[coverage labels,anchor=east]
      at (-.3,-\row-.5) {\base};
  }
\end{tikzpicture}
\end{minipage}\hfill
\begin{minipage}[t]{.63\linewidth}
\centering
\vspace{0pt} {\scriptsize $\mathcal R\subseteq\mathbb Z/288$\par $266/288\approx92.36\%$\par
}
\smallskip
\begin{tikzpicture}[x=3mm,y=3mm]
  \fill[coveragegrey] (0,0) rectangle (24,-12);

  \foreach \r in {
    16,18,26,30,74,79,80,84,112,146,159,160,
    170,174,175,176,180,218,234,242,255,256
  } {
    \pgfmathtruncatemacro{\col}{mod(\r,24)}
    \pgfmathtruncatemacro{\row}{floor(\r/24)}
    \fill[white] (\col,-\row) rectangle ++(1,-1);
  }

  \draw[black!45,line width=.3pt,step=1]
    (0,0) grid (24,-12);

  \foreach \col in {0,...,23}
    \node[coverage labels,anchor=south]
      at (\col+.5,.2) {\col};

  \foreach \row in {0,...,11} {
    \pgfmathtruncatemacro{\base}{24*\row}
    \node[coverage labels,anchor=east]
      at (-.3,-\row-.5) {\base};
  }
\end{tikzpicture}
\end{minipage}
\par\endgroup

\makeatletter
\def\@captype{figure}
\makeatother
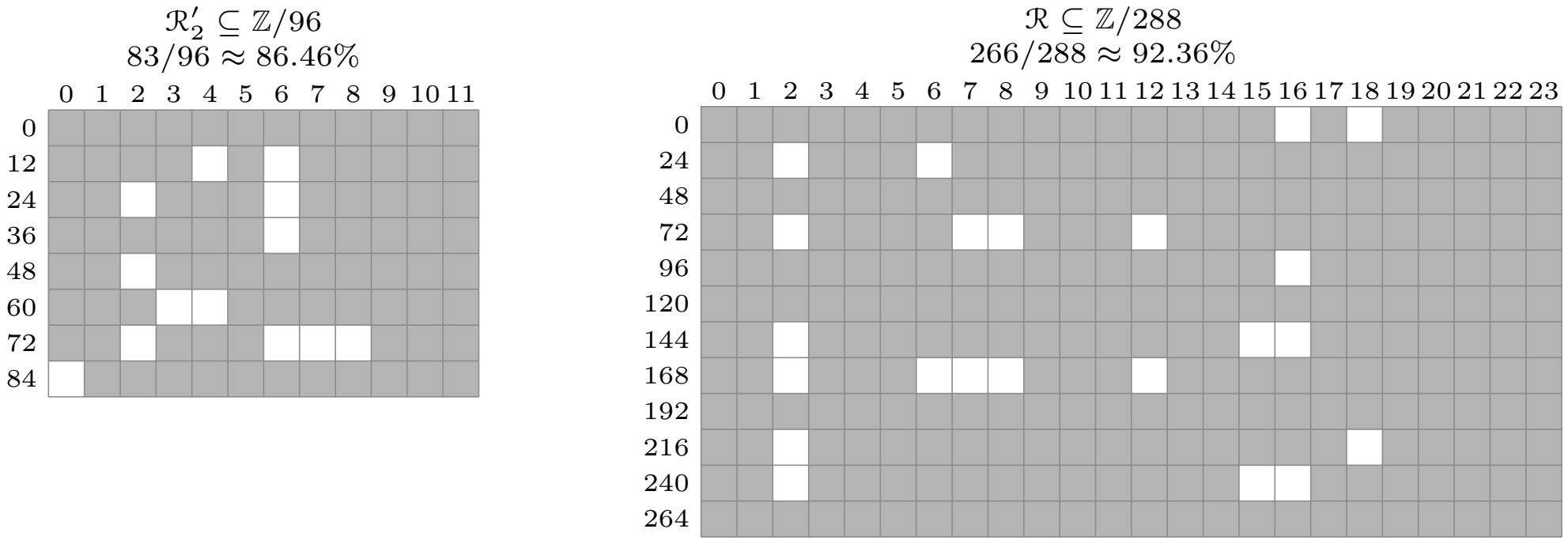
\caption{Residue classes in $\mathcal R_2'$ modulo $96$ and
$\mathcal R$ modulo $288$, shown in grey. The residue class that each cell represents is the sum of its row and column labels.}
\label{fig:residue-coverage}
\end{minipage}
\par\medskip

In future work, we plan to prove Conjecture~\ref{conj:exotic-cpn} in full by resolving the remaining cases among dimensions $n\geq4$ with $n\bmod96\in\overline{\mathcal R_2'}$.

\end{rmk}

\subsection{From manifold topology to stable homotopy theory}\label{subsec:intro-methods}
To prove Theorem \ref{main theorem prime 2 and 3}, we first use smoothing theory to translate the geometric existence problem into a computation in stable homotopy theory. Classically, the set of concordance classes of smooth structures of $\CP^n$ can be identified with a subset of its stable cohomotopy group. Nonetheless, we need a way to recognize which nonzero stable cohomotopy classes represent exotic smooth structures. We establish a detection criterion (Theorem~\ref{Detection}) that gives explicit bounds on the skeletal filtration under which a nonzero stable cohomotopy class determines an exotic smooth structure.

% Final review (2026-11-04), original line 568.

For a closed smooth manifold $M$ of real dimension at least five, the classification theorem of Kirby--Siebenmann \cite{KirbySiebenmannfoundational} identifies concordance classes of marked smoothings of its underlying topological manifold with $[M,Top/O]$. $Top/O$ is an infinite loop space whose definition we recall in Section~\ref{smoothingcomparison}. Such homotopy classes of maps thus form a group, which is the \emph{concordance group of $M$}, denoted by $C(M)$.

% Final review (2026-11-04), original line 574.

Precomposition with complex conjugation $\tau:\CP^n\to\CP^n$ defines an automorphism
\[\tau^\ast \colon C(\CP^n) \to C(\CP^n), \qquad (\tau^\ast)^2 = \id.\]
By Sullivan's correspondence \cite{Sullivan1996triangulating}, there is a bijection
\[
\{\text{The }C_2\text{-orbits under the action of } \tau^\ast\} 
\;\cong\;
\left\{
\begin{array}{c}
\text{unoriented diffeomorphism classes of smooth}\\
\text{manifolds homeomorphic to $\CP^n$}
\end{array}
\right\}.
\]
Consequently, $\CP^n$ admits an exotic smooth structure if and only if the concordance group $C(\CP^n)$ is nontrivial.

Brumfiel \cite{Brumfiel1971almostsmooth} shows that there is an injection of underlying sets from $C(\CP^n)$ to the stable cohomotopy group $[\Sigma^\infty\CP^n, S^0]$. We abbreviate this group by $[\CP^m,S^0]$ when there is no risk of ambiguity. When localizing at odd primes, Sullivan \cite{Sullivan1996triangulating} shows that the two groups coincide.

However, at the prime $2$, two distinctions require attention:
\begin{enumerate}
    \item In terms of their underlying sets, the $2$-component of $C(\CP^n)_{(2)}$ could be a proper subset of $[\CP^n, S^0]_{(2)}$.
    \item The group structures are also different, which we recall in Theorem~\ref{thm: distinct group structures}.
\end{enumerate}

Hence we need a criterion by which we can determine which elements of the stable cohomotopy group represent actual concordance classes. A particularly useful criterion we developed is based on the Atiyah-Hirzebruch spectral sequence (whose detailed setup we postpone until Subsection~\ref{spectralsequences}). For now, we recall only what is needed to state the criterion:

The skeletal filtration on $[\CP^n,S^0]$ is 
\begin{equation*}
F^{2s}[\CP^n,S^0]_{(2)}
=\ker\bigl([\CP^n,S^0]_{(2)}\longrightarrow [\CP^{s-1},S^0]_{(2)}\bigr),
\qquad 1\leq s \leq n.
\end{equation*}
It gives rise to the Atiyah-Hirzebruch spectral sequence (AHSS)
\[
\begin{aligned}
E_1^{2s,t}(\CP^n;S^0_{(2)})
&= \pi_{-t+2s}(S^{2s})_{(2)}
\Longrightarrow [\CP^n,S^{2s+t}]_{(2)}
 ,\\
d_{2r}:E_{2r}^{2s,t}&\longrightarrow E_{2r}^{2(s+r),t-2r+1},\qquad r\geq1.
\end{aligned}
\]
where $2s$ is the skeletal filtration and $2s+t$ is the total cohomological degree. We write $\alpha[2s] \in E_1^{2s,-t}$ for an element in the AHSS $E_1$-page with skeletal filtration $2s$ and total degree $2s-t$, where $\alpha \in \pi_{t}S^0$. In particular, any Atiyah--Hirzebruch cycle converging to $[\CP^n,S^0]$ is of the form $\alpha[2\ell], \alpha \in \pi_{2\ell}S^0$.

Under this convention, we prove in Subsection~\ref{filtrationsmoothing} that any element with sufficiently high filtration actually detects an exotic smooth structure, namely:

\begin{thmx}\label{Detection}
Let $n \geq 3$ and suppose $2^{r-1}\leq\ell\leq n\leq2^r-2.$ Every nonzero surviving cycle $\alpha[2\ell] \in E_\infty^{2\ell,-2\ell}$ determines an exotic smooth structure on $\CP^n$.
\end{thmx}

For example, in Theorem~\ref{thm:tmf-cp-10-45-58-82} we show that $2\bar \kappa[20]$ survives in the AHSS for $\CP^{10}$ and $\CP^{11}$. It follows that there are smooth manifolds $M^{20}$ and $M^{22}$ homeomorphic but not diffeomorphic to $\CP^{10}$ and $\CP^{11}$, respectively. Moreover, $M^{20}$ and $M^{22}$ can be constructed compatibly, with their smooth structures geometrically realizing a cohomotopy class on $\CP^{10}$ detected by $2\bar \kappa[20]$, and its extension to $\CP^{11}$, respectively. We explain this geometric interpretation of Theorem~\ref{Detection} and the resulting iterative construction in Subsection~\ref{subsec: methods introduction}.

\hypertarget{intro:three-strategies}{}%
% The above theorem now supplies a solution to the geometric problem via an algebraic topological computation. We next briefly explain the strategy of proving Theorem~\ref{thm:explicit-coverage}:

\subsection{Methods for finding nonzero cohomotopy classes}\label{subsec: methods introduction}

In principle, our methodology consists of the following two steps:
\begin{enumerate}
    \item\label{topcellclass} First find some nonzero $ \alpha[2n]$ supplied by the top cell in the $p$-local AHSS of $\CP^n$.
    \item\label{extension} Then we seek to further extend its survival in the $p$-local AHSS of $\CP^N$ for $N>n$. 
\end{enumerate}

At the prime $2$, after we establish the nontriviality of $\alpha[2n]$ in the AHSS of $\CP^N$, so long as the numerical condition of Theorem~\ref{Detection} is satisfied, then exotic $\CP^\ell$'s exist for all $n \leq \ell \leq N$. 

To prove Theorem~\ref{thm:duo}, we develop three methods at chromatic height $2$ for detecting nonzero stable cohomotopy classes. The height-$2$ input comes from $\tmf$, the connective spectrum of topological modular forms, and the $\beta$-families in the stable homotopy groups of spheres.

The first two methods use the unit map $u\colon S^0\to\tmf$. For $\alpha\in\pi_{2n}S^0$, if $u_*(\alpha)[2n]$ survives nontrivially in the $\tmf$-AHSS of $\CP^n$, it detects the nonzero composite
\[
\CP^n\xrightarrow{q}S^{2n}
\xrightarrow{\alpha}S^0
\longrightarrow\tmf.
\]
Consequently, $\alpha[2n]$ survives nontrivially in the AHSS for $[\CP^n,S^0]$.

The third method uses the \emph{$BP$-Mahowald square}. More precisely, it combines the $BP$-algebraic Atiyah--Hirzebruch spectral sequence (algAHSS) with the Adams--Novikov spectral sequence (ANSS) to establish the nontrivial survival of $\beta[2n]$ in the AHSS, where $\beta\in\pi_{2n}S^0$ is a $\beta$-family element.

We now outline the idea behind each method and the main difficulty it addresses.
\begin{enumerate}
    \item \textbf{The $\tmf$-AHSS differential method}, developed in Subsection~\ref{sec:tmf-detection} and Section~\ref{sec: tmf-differential method}, rests on the following observation: if $u_*(\alpha)[2n]$ supports a nonzero differential in the $\tmf$-AHSS of $\CP^\infty_n$, then $\alpha[2n]$ survives nontrivially in the AHSS for $[\CP^n,S^0]$. The main difficulties lie in computing long differentials and choosing suitable Hurewicz preimages whose associated classes survive nontrivially in the AHSS of larger projective spaces.
    
    The results obtained by this method cover $13$ residue classes in $\mathcal R_2$.
    
    \item \textbf{The $\tmf$-AHSS extension method}, developed in Subsection~\ref{sec: delta method} and Section~\ref{sec: tmf-extension method}, uses hidden $2$-extensions to establish the survival of candidate classes. We first construct these extensions in finite stunted projective spaces and then pull them back to the full projective spaces. The main difficulty is to rule out a jump of filtration under pullback, ensuring that the nonzero classes retain their prescribed AHSS representatives. 
    
    The results obtained by this method cover $18$ residue classes in $\mathcal R_2$.
    
     \item \textbf{The $\beta$-family method via the algAHSS}, developed in Subsection~\ref{section: ANSS methodology} and Section~\ref{sec:anss-detection}, uses the $BP$-Mahowald square to detect $\beta$-family elements. Let $\beta \in \Ext^{2,2n+2}_{BP_\ast BP}(BP_\ast, BP_\ast)$ be an permanent cycle in the ANSS that detects $b \in \pi_{2n}S^0$. Write $\beta[2n]$ for its image under the induced map
    \[\Ext^{2,2n+2}_{BP_\ast BP}(BP_\ast, BP_\ast) \xrightarrow{q^\ast} \Ext^{2,2}_{BP_\ast BP}(BP_\ast\CP^n, BP_\ast)\]
    where $q$ is the pinch map onto the top cell. We prove that $\beta[2n]$ is a nonzero permanent cycle in the ANSS of $\CP^n$, and that the homotopy class it detects has AHSS name $b[2n]$.

    The main difficulty is to show that $\beta[2n]$ survives nontrivially. In the algebraic Atiyah--Hirzebruch spectral sequence (algAHSS), it could potentially be hit by an incoming differential from the ANSS $1$-line. We first compute the $d_2$-differentials and verify that none hits $\beta[2n]$. We then construct nonzero classes in the $\Ext$ group whose algAHSS names are precisely the remaining $d_2$-cycles. The strategy is first to use the logarithm of the universal $2$-typical formal group law to construct classes in $\Ext^{0,2}_{BP_\ast BP}(BP_\ast\CP^\infty,(BP_\ast)_{\mathbb Q}/BP_\ast)$, and then to show that their images under the rational Bockstein in $\Ext^{1,2}_{BP_\ast BP}(BP_\ast\CP^\infty,BP_\ast)$ have the desired algAHSS names. Naturality under restriction from $\CP^\infty$ to $\CP^n$ then guarantees that none of these $d_2$-cycles can support a higher differential hitting $\beta[2n]$ in the algAHSS converging to $\Ext^{\ast,\ast}_{BP_\ast BP}(BP_\ast\CP^n,BP_\ast)$. The results obtained by this method cover $22$ residue classes in $\mathcal R_2$.

\end{enumerate}
All three methods supply nonzero classes $x[2n] \in E_\infty^{2n,-2n}$ in the AHSS of $\CP^n$. Should they also survive in the AHSS of $\CP^N$ for $N>n$, applying Theorem~\ref{Detection} will supply exotic complex projective spaces in all dimensions from $n$ to $N$. Together with $24$ residue classes $1,5,\ldots,93$ modulo $96$ of Farrell-Jones \cite{FarrellJones1994hyperbolic}, these results yield all $77$ classes in $\mathcal R_2$. 

Theorem~\ref{thm:tres} is obtained by applying the above $\beta$-families in algAHSS method at the prime $3$. For Theorem~\ref{thm:unus}, we use divided beta families at all $p \geq 5$. The set of complex dimensions covered by these families has natural density one.

\begin{rmk}
The survival of $\alpha[2n]$ and its extension to $\CP^N$ admit the following geometric interpretation, within the range of Theorem~\ref{Detection}.

For the classes considered here, the image of $\alpha\in\pi_{2n}S^0$ in $\operatorname{coker}J_{2n}$ is nonzero and has vanishing Kervaire invariant. By Kervaire--Milnor \cite{kervairemilnor1963groups}, this image is represented by an exotic sphere $[\Sigma^{2n}]\in\Theta_{2n}$. The class $\alpha[2n]$ survives nontrivially in the AHSS for $\CP^n$ if and only if
\[
M=\CP^n\#\Sigma^{2n}
\]
is homeomorphic but not diffeomorphic to $\CP^n$. We assume this nontriviality in what follows.

Let $f\colon M\to\CP^n$ be the connected-sum homeomorphism, and choose a smooth map $f'$ homotopic to $f$. Let $L$ be the dual of the tautological line bundle over $\CP^n$, and set
\[
L'=(f')^*L.
\]
Then $L'$ is a smooth complex line bundle over $M$. Write $D(L')$ and $S(L')$ for its unit disk bundle and unit circle bundle, respectively. The bundle $L'$ is topologically isomorphic to $f^*L$, and hence the total space $S(L')$ is homeomorphic to $S^{2n+1}$. Within the range of Theorem~\ref{Detection}, we have
\[
\alpha[2n]\text{ survives in the AHSS for }\CP^{n+1}
\quad\Longleftrightarrow\quad
S(L')\cong_{\mathrm{diff}}S^{2n+1}.
\]
When this holds, choose a boundary diffeomorphism $\varphi\colon S^{2n+1}\to S(L')$ and form
\[
M_{n+1}=D(L')\cup_{\varphi}D^{2n+2}.
\]
This manifold is homeomorphic to $\CP^{n+1}$. The gluing can be chosen so that its smooth structure realizes the prescribed cohomotopy extension; by Theorem~\ref{Detection}, it is then an exotic smooth structure on $\CP^{n+1}$.

By Theorem~\ref{finiteinfinitecriterion}, the class $\alpha[2n]$ eventually supports a nonzero outgoing differential in the spherical AHSS for $\CP^\infty/\CP^{n-1}$. Geometrically, any compatible iteration must reach a finite stage at which the total space of the corresponding circle bundle is an exotic sphere and hence cannot be capped off by a standard disk. Consequently, a compatible construction realizing this nonzero cohomotopy class cannot continue indefinitely.

\end{rmk}

\subsection{Organization}
The remainder of this paper consists of the following parts:
\begin{itemize}
    \item Section~\ref{smoothingcomparison} gives a brief review of smoothing theory. Using this comparison, Section~\ref{sec:filtration-criteria} establishes Theorem~\ref{Detection} which detects exotic $\CP^n$ by elements in its stable cohomotopy group.
    \item Section~\ref{sec:ahss-two} elaborates on the mechanism of the aforementioned three techniques at the prime $2$. Sections~\ref{sec: tmf-differential method}, \ref{sec: tmf-extension method}, and~\ref{sec:anss-detection} then carry out the detailed computations which prove Theorem~\ref{thm:duo} in combination.
    \item Section~\ref{sec:three-primary-detection} proves Theorem~\ref{thm:tres} using $3$-primary computations. Section~\ref{sec:density-one} proves our density one Theorem~\ref{thm:unus}, by combining beta-families at all primes $p\geq5$.

\end{itemize}

\subsection*{Acknowledgement} 
The second author was supported by “the Fundamental Research Funds for the Central Universities” No. 226-2025-00092 (China). The third author is partially supported by the AMS Centennial Research Fellowship and NSF Grant DMS-250624. The authors thank Agnes Beaudry, Mark Behrens, Guchuan Li, Sihao Ma, XiaoLin Danny Shi, Guozhen Wang, Yuchen Wu, and Shangjie Zhang for helpful discussions. The authors thank J.D. Quigley for his talk on smooth $S^1$-actions on exotic spheres at the 2025 IWoAT conference at Zhejiang University, where this project began.

\section{Concordance classes and stable cohomotopy}\label{smoothingcomparison}

In this section we recall the smoothing-theoretic comparison between concordance classes of smooth structures and reduced stable cohomotopy of $\CP^n$.

A \emph{marked smoothing of $\CP^n$} is a pair $(N,h)$, where $N$ is a closed smooth manifold and $h:N\to\CP^n$ is a homeomorphism. Two marked smoothings are concordant if the smooth structures they induce on $\CP^n$ extend to a smooth structure on $\CP^n\times[0,1]$ with the prescribed structures at its two ends. Let $C(\CP^n)$ denote their concordance classes.

Next we recall the definition and properties of some infinite loop spaces that play an essential role in computing concordance classes.

\begin{defn}[{\cite{BrumfielMadsenMilgram1973PLclass,BoardmanVogt1973loopspace}}]\label{classifyingspaces}
Let $G=\colim_r G(r)$ be the stable monoid of based homotopy self-equivalences of $S^r$, stabilized by suspension, and let $SG$ be its degree-one component, based at the identity. Let $PL=\colim_r PL(r)$, where $PL(r)$ consists of the piecewise linear homeomorphisms of $\bbR^r$ fixing the origin. Let $Top=\colim_r Top(r)$, where $Top(r)$ consists of the homeomorphisms of $\bbR^r$ fixing the origin. For the forgetful maps between the stable classifying spaces, write
\[
Y/X=\operatorname{hofib}(BX\longrightarrow BY).
\]
The forgetful maps
\[
BO\longrightarrow BPL\longrightarrow BTop\longrightarrow BG
\]
are maps of infinite loop spaces, so the homotopy fibers defined above are infinite loop spaces. In particular, the composite $BO\to BG$ is the stable $J$-map.
\end{defn}

The following theorem of Kirby-Siebenmann identifies concordance classes of $M$ with maps into the infinite loop space $Top/O$:

\begin{theorem}[Kirby--Siebenmann {\cite{KirbySiebenmannfoundational}}]\label{KSclassification}
Let $M$ be a closed topological manifold of real dimension $d\geq5$, and choose a reference marked smoothing $(M_0,h_0)$. There is a bijection
\[
C(M)\longrightarrow[M,Top/O]
\]
from concordance classes of marked smoothings to homotopy classes of maps, sending $(M_0,h_0)$ to the constant map. We give $C(M)$ the corresponding group structure via this bijection.
\end{theorem}

For $\CP^n$ and $n \geq 3$, we take the standard smooth structure with its identity marking as reference. We base $\CP^n$ at $\CP^0$ and use based homotopy classes from now on. The next theorem of Sullivan's shows that finding an exotic $\CP^n$ is equivalent to finding a nonzero concordance class.

\begin{theorem}\cite{Sullivan1996triangulating}\label{concordanceCP}
For $n\geq3$, every nonzero class in $C(\CP^n)$ is represented by a smooth manifold homeomorphic but not diffeomorphic to the standard $\CP^n$. In particular, an exotic $\CP^n$ exists if and only if $C(\CP^n)\neq0$.
\end{theorem}

\begin{prop}\label{LES}

Let $j:SG\to G/O$ be the map in the oriented fiber sequence
\[
SO\xrightarrow{J}SG\xrightarrow{j}G/O\longrightarrow BSO.
\]
Let $q$ be the induced map in the fiber sequence
\[
Top/O \to G/O \xrightarrow{q} G/Top.
\]
Denote the composite $q \circ j$ by $i$. There is a natural isomorphism of abelian groups
\[\ker\bigl(i_*:[\CP^n,SG]\longrightarrow
[\CP^n,G/Top]\bigr)\xrightarrow{\ \cong\ }C(\CP^n).\]
\end{prop}

\begin{proof}
By \cite{Brumfiel1971almostsmooth}, for every $n\geq3$, there are exact sequences of abelian groups
\[
0\longrightarrow C(\CP^n)\xrightarrow{u_*}[\CP^n,G/O]
\xrightarrow{q_*}[\CP^n,G/Top]
\]
\[
0\longrightarrow[\CP^n,SG]\xrightarrow{j_*}[\CP^n,G/O]
\longrightarrow\bbZ^{\lfloor n/2\rfloor}\longrightarrow0.
\]
The two exact sequences, together with the fact that $C(\CP^n)$ and $[\CP^n,SG]$ are torsion groups, prove the proposition.

\end{proof}

Furthermore, Brumfiel's calculation identifies the underlying sets of $[\CP^n,SG]$ and the stable cohomotopy group $[\CP^n, S^0]$. The stable monoid $G$ is a model for $GL_1(S)$, the degree $\pm1$ components of $QS^0$. Translation by the degree-one point gives a based homotopy equivalence
\[
t \colon Q_0S^0\xrightarrow{\ \simeq\ }SG,
\qquad x\longmapsto1+x.
\]
Postcomposition therefore gives a natural pointed-set bijection $[\CP^n, Q_0S^0]\to[\CP^n,SG]$. The group structures require distinction: the group structure of $[\CP^n,SG]$ is induced by composition of stable sphere self-equivalences, whereas the group structure on the stable cohomotopy group is induced by addition on $QS^0$.

\begin{theorem}[Brumfiel {\cite[\S 4]{Brumfiel1971almostsmooth}}]\label{thm: distinct group structures}
For $x, y \in [\CP^n,S^0]$ define
\[
x\star y=x+y+xy,
\]
where the product is induced by multiplication on $S^0$. Consequently, $([\CP^n,S^0],\star)$ identifies with $[\CP^n,SG]$ as groups.
\end{theorem}

\begin{rmk}\label{translation preserves filtration}
Attaching further cells is compatible with translation, in the sense that the cofiber sequence $\CP^{n-1}\to\CP^n\to S^{2n}$ gives the commutative diagram
\[
\begin{tikzcd}[column sep=large]
\pi_{2n}(SG) \arrow[r] \arrow[d,"\cong"'] &
{[\CP^n,SG]} \arrow[r] \arrow[d] &
{[\CP^{n-1},SG]} \arrow[d] \\
\pi_{2n}S^0 \arrow[r] &
{[\CP^n,S^0]} \arrow[r] &
{[\CP^{n-1},S^0]}.
\end{tikzcd}
\]    
\end{rmk}

\section{Filtration and smoothing criteria}\label{sec:filtration-criteria}

By Theorem~\ref{concordanceCP} and Proposition~\ref{LES}, every nonzero concordance class of $\CP^n$ determines an exotic $\CP^n$, whereas the set of concordance classes is a subset of the stable cohomotopy group. In this section we establish criteria by which we can determine which stable cohomotopy classes of $\CP^n$ belong to the concordance group.

Our criteria divide into an odd primary part and a $2$-primary part. The odd primary part is an easy corollary of \cite{Sullivan1996triangulating}. In the $2$-local case, we show that so long as an element has sufficiently high Atiyah-Hirzebruch filtration, then it belongs to the concordance group. 

\begin{theorem}\label{thm:surviving-smoothings}
For $n \geq 3$, 

\begin{enumerate}[label=(\roman*),ref=\thetheorem(\roman*),leftmargin=*]

\item\label[theorem]{oddAHSSsmoothing} Let $p$ be an odd prime. Every nonzero element of $[\CP^n,S^0]_{(p)}$ determines a nonzero element of $C(\CP^n)_{(p)}$ and hence an exotic smooth structure of $\CP^n$.

\item\label[theorem]{AHSSsmoothing} In the $2$-local AHSS for $[\CP^n,S^0]$, suppose that $2^{r-1}\leq\ell\leq n\leq2^r-2$. Every nonzero class in the degree-zero abutment $[\CP^n,S^0]_{(2)}$ detected by
\[
\alpha[2\ell]\in E_\infty^{2\ell,-2\ell}
(\CP^n;S^0_{(2)})
\]
determines a nonzero element of $C(\CP^n)_{(2)}$ and hence an exotic smooth structure of $\CP^n$.

\end{enumerate}
\end{theorem}

 In the remainder of this section, we first set up the Atiyah-Hirzebruch spectral sequences (AHSS) in Subsection~\ref{spectralsequences}, then briefly recall Sullivan's theorem at odd primes in Subsection~\ref{oddprimarycomparison}. Finally, we prove our $2$-primary criterion in Subsection~\ref{filtrationsmoothing}.

\subsection{The Atiyah-Hirzebruch spectral sequence (AHSS)}\label{spectralsequences}

For a finite based CW complex $Y$ and a coefficient spectrum $E$, write
\[
E^q(Y)=[\Sigma^\infty Y,\Sigma^q E].
\]
Let $Y^{s}$ be the skeleton in real dimension $s$, with $Y^{-1}=*$. Applying the function-spectrum functor to the skeletal inclusions gives the tower
\[
\cdots\longrightarrow\Fun(\Sigma^\infty Y^{s},E)
\longrightarrow\Fun(\Sigma^\infty Y^{s-1},E)
\longrightarrow\cdots.
\]
The associated cohomological Atiyah--Hirzebruch spectral sequence thus begins with
\begin{equation*}
\begin{aligned}
E_1^{s,t}(Y;E)
&= E^{s+t}(Y^{s}/Y^{s-1}) = \bigoplus_{i=1}^{r_s} \pi_{-t}E
 \\
d_r:E_r^{s,t}&\longrightarrow E_r^{s+r,t-r+1},\qquad r\geq1.
\end{aligned}
\end{equation*}
Here $s$ is the skeletal filtration, $s+t$ is the total (cohomological) degree, and $r_s$ is the number of $s$-dimensional cells. Starting from the $E_2$-page, 
\[
E_2^{s,t}(Y;E)=H^s(Y;\pi_{-t}E)
\Longrightarrow E^{s+t}(Y).
\]
Since $Y$ is finite, the abutment has the finite skeletal filtration
\begin{equation*}
\begin{aligned}
F^s E^n(Y)
&=\ker\bigl( E^n(Y)\longrightarrow
                  E^n(Y^{s-1})\bigr),\\
E_\infty^{s,t}(Y;E)
&\cong F^s E^{s+t}(Y)/F^{s+1} E^{s+t}(Y).
\end{aligned}
\end{equation*}

For a cell of dimension $\ell$, we write $\alpha[\ell]\in E_1^{\ell,-k}(Y;E)$ for the class with coefficient $\alpha\in\pi_kE$. Its position on $E_1$ is $(s,t)=(\ell,-k)$, and its total degree is $\ell-k$. In particular, in the AHSS of $\CP^n$, all the elements of total degree zero have the form $\alpha[2\ell]$ where $\alpha\in\pi_{2\ell}E$.

In particular, under the Adams convention of spectral sequence charts where we take $-(s+t)$ as the $x$-coordinate and $s$ as the $y$-coordinate, a $d_{r}$-Atiyah-Hirzebruch differential is an arrow going left by one unit and going up by $r$ units.

\begin{rmk}
All explicit AHSS computations in this paper involve finite CW complexes, so convergence is automatic. In Theorem~\ref{finiteinfinitecriterion}, we also consider the AHSS of the infinite stunted projective space $\CP^\infty_n$, where $n\geq1$. In this case the inverse systems consist of finite groups and therefore satisfy the Mittag--Leffler condition, which ensures strong convergence.
\end{rmk}

For explicit formulas that compute the Atiyah-Hirzebruch differentials and hidden extensions in certain finite cell complexes, we refer the reader to \cite{wangxu2017triviality}, \cite{wangxu201851stem} and \cite{HouZhang2026C3}. These results apply to complexes with at most four cells. For spaces or spectra with more complicated cell structures, we construct suitable zigzags of maps relating them to simpler complexes. We then compute the relevant differentials or hidden extensions in these simpler complexes and use naturality to pull them back or push them forward to the original objects. A detailed analysis will be carried out in Sections~\ref{sec: tmf-differential method} and~\ref{sec: tmf-extension method}.

\subsection{The comparison at odd primes}\label{oddprimarycomparison}

With $2$ inverted, the concordance group of $\CP^n$ is isomorphic to its stable cohomotopy group. The identification is an easy corollary of a theorem of Sullivan.

\begin{theorem}\cite{Sullivan1996triangulating}\label{Sullivanodd}
For every odd prime $p$ and $n>2$, there is an isomorphism of finite abelian $p$-groups
\[
C(\CP^n)_{(p)}\cong [\Sigma^\infty\CP^n, S^0]_{(p)}.
\]
\end{theorem}
\begin{proof}
Sullivan's proposition originally identifies
\[
[\CP^n,PL/O]\otimes\bbZ[1/2]
\cong [\Sigma^\infty\CP^n, S^0]\otimes\bbZ[1/2].
\]
The result then follows immediately from the facts (see \cite{BrumfielMadsenMilgram1973PLclass}) that
\[
PL/O\longrightarrow Top/O\longrightarrow Top/PL \text{ is a fiber sequence and }
Top/PL\simeq K(\bbZ/2,3).
\]

\end{proof}

\begin{proof}[Proof of Theorem~\ref{oddAHSSsmoothing}]
This follows immediately by combining Theorems~\ref{concordanceCP} and~\ref{Sullivanodd}.
\end{proof}

\subsection{A two-primary filtration criterion}\label{filtrationsmoothing}

By Proposition~\ref{LES}, the concordance group is the kernel of $i_\ast: [\CP^n,SG] \to [\CP^n,G/Top]$, where $i$ is the composition
\[SG \xrightarrow{j} G/O \xrightarrow{q} G/Top.\]
The following theorem gives a complete description of the $2$-local $G/Top$.

\begin{theorem}[Sullivan {\cite{BrumfielMadsenMilgram1973PLclass}}]\label{G/Top}
There are fundamental classes
\[
k_{4a-2}\in H^{4a-2}((G/Top)_{(2)};\bbF_2),
\qquad L_{4a}\in H^{4a}((G/Top)_{(2)};\bbZ_{(2)}),
\quad a\geq1,
\]
whose induced map gives a homotopy equivalence of spaces
\[
\chi:(G/Top)_{(2)}\xrightarrow{\ \simeq\ }
\prod_{a\geq1}K(\bbF_2,4a-2)
\times\prod_{a\geq1}K(\bbZ_{(2)},4a).
\]
\end{theorem}

\begin{theorem}[{\cite{BrumfielMadsenMilgram1973PLclass}}]\label{map of cohomology}
Write $k_{4a}$ for the mod-two reduction of $L_{4a}$ in Theorem~\ref{G/Top}. In $H^\ast(SG;\mathbb{F}_2)$ there are relations
\begin{equation}\label{degree restraint}
i^* k_{2m}=0 \quad\text{unless }m=2^r\text{ or }m=2^r-1 \text{ for some }r\geq0,    
\end{equation}

and
\[
i^*k_{2^{r+1}}=(i^*k_2)^{2^r},\qquad r\geq0.
\]
\end{theorem}

\begin{proof}[Proof of Theorem~\ref{AHSSsmoothing}]
Let $f$ be a nonzero homotopy class detected by $\alpha[2\ell]$, and set $g=i\circ  {t}_\ast f$, where $ {t}_\ast$ is the translation of Remark~\ref{translation preserves filtration}. By Theorem~\ref{G/Top}, the based map $g:\CP^n\to(G/Top)_{(2)}$ is nullhomotopic if and only if
\begin{equation*}
\vphantom{\textnormal{(5)}}
g^*k_{4a-2}=0\quad\text{and}\quad g^*L_{4a}=0
\qquad\text{for every }a\geq1.
\end{equation*}

Every integral coordinate $g^*L_{4a}$ vanishes: since $g$ factors through the rationally contractible space $SG$, its image in $H^{4a}(\CP^n;\mathbb Q)$ is zero. The coefficient map
\[
H^{4a}(\CP^n;\mathbb Z_{(2)})
\longrightarrow H^{4a}(\CP^n;\mathbb Q)
\]
is injective, so $g^\ast L_{4a}=0$. Now consider $g^\ast k_{4a-2}$. By Remark~\ref{translation preserves filtration} $t_\ast f $ maps to zero under restriction to $[\CP^{\ell-1},SG]$, so $g^*k_{4a-2} =0 $ for all $4a-2 \leq 2\ell-2$. The degrees $4a-2 >2\ell$ fall in the range of Equation~\ref{degree restraint} where $i^\ast k_{4a-2} = 0$. It follows that $\tau_\ast f$ is in the kernel of 
\[i_\ast: [\CP^n, SG] \to [\CP^n, G/Top],\]
hence it determines an exotic $\CP^n$ by Theorem~\ref{concordanceCP} and Proposition~\ref{LES}.

\end{proof}

\begin{rmk}
One can slightly enhance the statement of Theorem~\ref{thm:surviving-smoothings}: $F^{2\ell}[\CP^n,S^0]_{(2)}$ is in fact a subgroup of $C(\CP^n)_{(2)}$. The product of any two elements of $F^{2\ell}[\CP^n,S^0]_{(2)}$ vanishes since $4\ell > 2n$; hence the addition $x\star y$ of Theorem~\ref{thm: distinct group structures} actually agrees with $x+y$ on this subgroup.
\end{rmk}

\section{Detection methods at $p=2$}\label{sec:ahss-two}

At the prime $2$, by Theorem~\ref{thm:surviving-smoothings}, any nonzero $\alpha[2\ell] \in E_\infty^{2\ell,-2\ell}$ of the AHSS of $\CP^n$ determines an exotic smooth structure of $\mathbb{CP}^n$, where $\alpha \in \pi_{2\ell}S^0$ with $\ell$ sufficiently large. In practice, direct computations in the AHSS of $\mathbb{CP}^n$ are often difficult. The complexity of AHSS computations arises from the following two aspects:

\begin{enumerate}[label=(\roman*),ref=\roman*,leftmargin=*]
\item The coefficients lie in the stable stems, which are known to be intricate;
\item Determining differentials requires detailed knowledge of the stable attaching maps in $\CP^n$.
\end{enumerate}
A natural idea is to replace the sphere spectrum by another ring spectrum $E$, which should have the following two features:
\begin{itemize}
    \item The Hurewicz image $u_\ast: \pi_\ast S^0 \to \pi_\ast E$ should be nonzero in sufficiently many degrees;
    \item The homotopy ring $\pi_\ast E$ should be well-understood.
\end{itemize}
For a spectrum $E$ that satisfies these two requirements, consider $\alpha \in \pi_{2n}S^0$ with nonzero Hurewicz image $u_\ast(\alpha) \in \pi_{2n}E$. If $u_\ast(\alpha)[2n]$ is a surviving cycle in the $E$-AHSS of $\CP^n$, then it detects the following nonzero homotopy class 
\[\CP^n \xrightarrow{q} S^{2n} \xrightarrow{\alpha} S^0 \xrightarrow{u} E.\]
Consequently, 
\[\CP^n \xrightarrow{q} S^{2n} \xrightarrow{\alpha} S^0\]
will be nonzero, i.e., $\alpha[2n]$ survives in the AHSS for $[\CP^n,S^0]$. If such $\alpha[2n]$ further survives in the AHSS for $[\CP^N,S^0]$ with $n < N$ satisfying the numerical conditions in Theorem~\ref{AHSSsmoothing}, then $\alpha[2n]$ will further determine an exotic smooth $\CP^N$. In summary, using the auxiliary spectrum $E$, we turn the problem of finding nonzero cohomotopy classes of $\CP^n$ into two steps:
\begin{enumerate}
    \item Finding surviving cycles $u_\ast(\alpha)[2n]$ in the $E$-AHSS of $\CP^n$;
    \item Further establishing the survival of $\alpha[2n]$ in the AHSS of $\CP^N$ for $N > n$.
\end{enumerate}

A natural candidate of $E$ is the spectrum of \emph{connective topological modular forms} ($\tmf$). Both its homotopy groups \cite{BrunerRognes2021tmf} and its Hurewicz image \cite{BMQhurewicz2021} have been completely computed, while its Hurewicz image has a 192-periodic pattern. This is the reason why our 2-primary computations produce families of exotic $\CP^N$'s, in dimensions differing by multiples of $96$. 

Our $\tmf$-AHSS computations develop two branches: 
\begin{enumerate}
    \item the $\tmf$-differential method of Subsection~\ref{sec:tmf-detection};
    \item the $\tmf$-extension method of Subsection~\ref{sec: delta method}.
\end{enumerate} 
The first method is based on the following observation:
\begin{itemize}
    \item If a filtration $\ell$ class $\alpha[\ell]$ supports a differential in the AHSS of $X$, then its restriction to the $\ell$-th skeleton $X^\ell$ cannot be hit by a differential.
\end{itemize}
Hence nonzero differentials
\[d_r(\alpha[2n]) = \gamma[2N]\]
in the $\tmf$-AHSS of $\CP^N$ imply the survival of $\alpha[2n]$ in the $\tmf$-AHSS of $\CP^n$.

\begin{figure}[htbp]
\centering
\begin{tikzpicture}[x=0.50cm,y=0.46cm,font=\small]
  \CPfourgrid{0}{3}{0}{10}
  \CPfourgrid{6}{9}{0}{10}
  \CPfourgrid{12}{15}{0}{10}

  \coordinate (xa) at (1,5);
  \coordinate (ya) at (0,9);
  \coordinate (xb) at (7,5);
  \coordinate (yb) at (6,9);
  \coordinate (ub) at (8,3);
  \coordinate (vb) at (8,1);
  \coordinate (xc) at (13,5);
  \coordinate (uc) at (14,3);
  \coordinate (vc) at (14,1);

  % Solid differentials.
  \CPfourarrow (xa) -- (ya);
  \CPfourarrow (xb) -- (yb);

  % Horizontal comparisons.
  \foreach \a/\b in {ya/yb,xa/xb,xb/xc,ub/uc,vb/vc} {
    \CPfourarrow[CPfour red,shorten <=0.08cm]
      (\a) -- (\b);
  }

  % Hypothetical incoming differentials.
  \CPfourarrow[CPfour red,dashed] (uc) -- (xc);
  \CPfourarrow[CPfour red,dashed] (vc) -- (xc);

  % Classes and labels.
  \foreach \p in {xa,ya,xb,yb,ub,vb,xc,uc,vc} {
    \fill (\p) circle (1.5pt);
  }
  \foreach \p in {xa,xb} {
    \node[below right] at (\p) {$x$};
  }
  \foreach \p in {ya,yb} {
    \node[above] at (\p) {$y$};
  }
  \node[above] at (xc) {$x$};
  \node[right=2pt] at (uc) {$x'$};
  \node[right=2pt] at (vc) {$x''$};

  \node[anchor=east] at (-0.7,9) {$n+r$};
  \node[anchor=east] at (-0.7,5) {$n$};
  \node[anchor=east] at (-0.7,1) {$n-r'$};

  \node at (1,-1.3) {$E_1(X_{n})$};
  \node at (7,-1.3) {$E_1(X)$};
  \node at (13,-1.3) {$E_1(X^n)$};

  \node[anchor=west,align=left] at (16,4.7) {
    $\nexists\,x'\text{ with }\operatorname{fil}(x')<n$\\[4pt]
    $\text{such that }d_r(x')=x.$
  };
\end{tikzpicture}
\caption{A pictorial demonstration of $\tmf$-differential methods. The dashed arrows depict trivial differentials.}
\label{fig:section4-E1-comparison}
\end{figure}
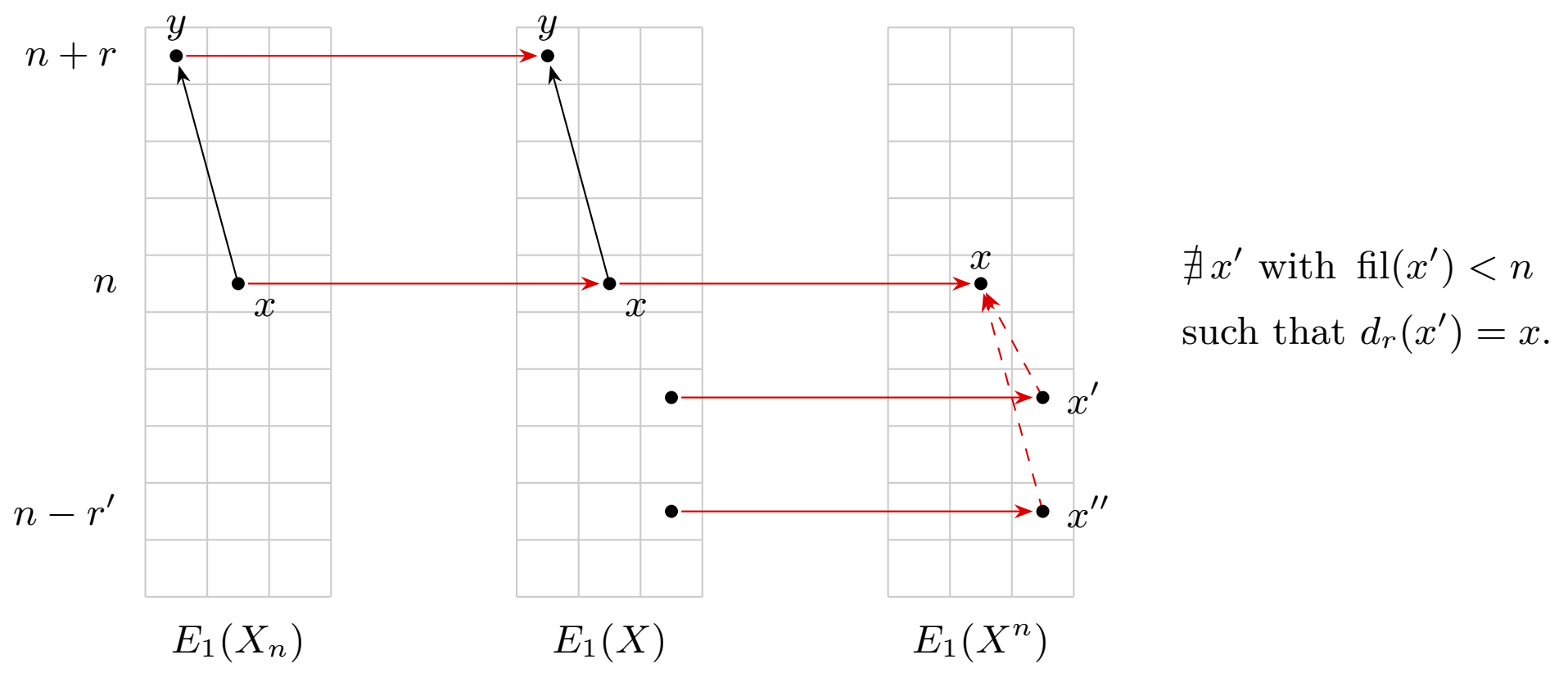

In order to show $\alpha[2n]$ survives in the $\tmf$-AHSS of $\CP^n$, the second method establishes hidden $2$-extensions
\[2^s \cdot \alpha[2n] = \delta[2N]\]
in the $\tmf$-AHSS of $\CP^N_n$. The classes $\delta \in \pi_{2N}\tmf$ are integral, detected by the powers of $\Delta$ (see Definition~\ref{defn:delta power classes}). The difference $N-n$ is small, so that the stunted $\CP^N_n$ is a finite complex where such an extension can be established.  

% Furthermore, as $\tmf$-modules, for $n \equiv -1 (\text{mod } 8)$ $\tmf \wedge \CP^n$ has particular nice properties such as top-cell splits off (Lemma~\ref{lem:tmf-projective-top-cell-splitting}). 

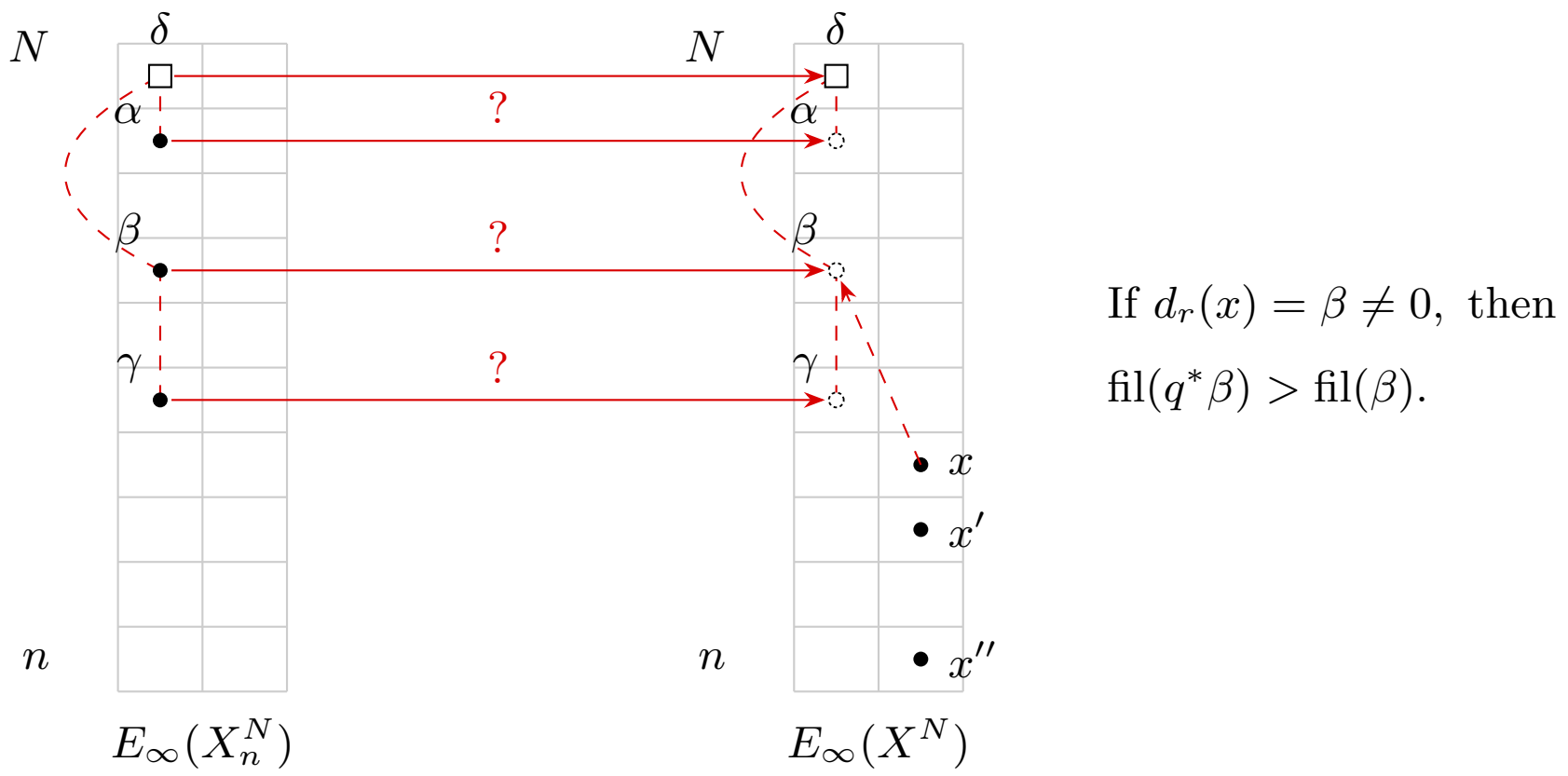
\begin{figure}[htbp]
\centering
\begin{tikzpicture}[x=0.60cm,y=0.46cm,font=\small]
  \CPfourgrid{0}{2}{0}{10}
  \CPfourgrid{8}{10}{0}{10}

  % Dashed extension connections and classes in each chart.
  % Their arrowheads are unspecified in the sketch.
  \foreach \s in {0,8} {
    \draw[CPfour red,dashed] (\s,8) -- (\s,9);
    \draw[CPfour red,dashed] (\s,4) -- (\s,6);
    \draw[CPfour red,dashed]
      (\s,6)
      .. controls (\s-1.5,7) and (\s-1.5,8) ..
      (\s,9);

        \foreach \h/\lab in {8/\alpha,6/\beta,4/\gamma} {
      \ifnum\s=8\relax
        \draw[draw=black,fill=white,line width=0.35pt,
              dash pattern=on 0.7pt off 0.6pt]
          (\s,\h) circle (1.5pt);
      \else
        \fill (\s,\h) circle (1.5pt);
      \fi
      \node[above left] at (\s,\h) {$\lab$};
    }

    \node[draw,fill=white,rectangle,inner sep=0pt,
          minimum size=4.5pt] at (\s,9) {};
    \node[above=3pt] at (\s,9) {$\delta$};

    \node[anchor=east,yshift=6pt] at (\s-1.1,9) {$N$};
    \node[anchor=east] at (\s-1.1,0) {$n$};
  }

  % Horizontal comparisons.
  \CPfourarrow[CPfour red,shorten <=0.10cm]
    (0,9) -- (8,9);

  \foreach \h in {8,6,4} {
    \CPfourarrow[CPfour red,shorten <=0.08cm]
      (0,\h) -- (8,\h) node[midway,above] {$?$};
  }

  % Lower-filtration classes in the right-hand chart.
  \foreach \h/\lab in {3/x,2/{x'},0/{x''}} {
    \fill (9,\h) circle (1.5pt);
    \node[right=2pt] at (9,\h) {$\lab$};
  }
  \CPfourarrow[CPfour red, dashed]
    (9,3) -- (8,6);

  \node at (0.5,-1.3) {$E_\infty(X_n^N)$};
  \node at (8.5,-1.3) {$E_\infty(X^N)$};

  \node[anchor=west,align=left] at (11,4.8) {
    $\text{If }d_r(x)=\beta\ne0,\text{ then}$\\[6pt]
    $\operatorname{fil}(q^*\beta)>\operatorname{fil}(\beta).$
  };
\end{tikzpicture}
\caption{Comparison of the $E_\infty$-pages. The question marks and dashed dots imply that, before applying the extension method, it is not justified that the classes $\alpha,\beta,\gamma$ survive to the $E_\infty$-page of $\tmf$-AHSS of $X^N$.}
\label{fig:section4-Einfty-comparison}
\end{figure}

We then consider the induced map
\[q^\ast \colon \tmf\text{-AHSS of } \CP^N_n \to \tmf\text{-AHSS of } \CP^N.\]
$\delta$ being integral shows that $\delta[2N]$ must survive to the $E_\infty$-page of both spectral sequences. Consequently, 
\[\text{The homotopy class detected by }\alpha[2n] \text{ must have a nontrivial image under } q^\ast.\]
In $\tmf^0(\CP^N)$, this nontrivial homotopy class could still be detected by some element whose filtration is higher than $2n$. After ruling out such filtration jumps, naturality would force the survival of $\alpha[2n]$ in the $\tmf$-AHSS of $\CP^N$ and $\CP^n$, respectively.

\begin{rmk}
    We emphasize the following distinction between the $S^0$- and $\tmf$-AHSSs:
    
    let $u \colon S^0 \to \tmf$  and $\alpha \in \pi_{2n}S^0$. A nonzero class $u_\ast\alpha[2n]$ supplied by the $\tmf$-extension method could either support a differential in the $\tmf$-AHSS of $\CP^\infty$, or it could survive in the $\tmf$-AHSS of $\CP^\infty$. However, the class $\alpha[2n]$ should nevertheless support a differential in the AHSS of $\CP^\infty$. There is no contradiction in such a phenomenon as the target of the AHSS differential could have trivial Hurewicz image.
\end{rmk}

Both of our $\tmf$-motivated strategies establish the survival of top-cell elements without ruling out every potential incoming differential page by page. Indeed when the dimension $n$ grows larger, analyzing the stable attaching maps of $\CP^n$ grows increasingly difficult. To circumvent this difficulty we developed a third strategy by a \emph{$BP$-Mahowald square} argument. Namely, by $BP$-Mahowald square we refer to the following diamond-shaped diagram of spectral sequences:
\[
\begin{tikzpicture}[
  baseline=(current bounding box.center),
  every node/.style={inner sep=2pt},
  ssarrow/.style={
    double, double distance=.55pt, line width=.28pt,
    -{Implies[length=4.5pt,width=4.5pt]}
  },
  sslabel/.style={
    fill=white, inner xsep=2pt, inner ysep=1.5pt,
    font=\scriptsize\bfseries
  }
]
\node (top) at (0,0)
  {$\textstyle\bigoplus_{i=1}^{n}
    \operatorname{Ext}_{BP_\ast BP}^{s,t+2i}(BP_*,BP_*)$};

\node (left) at (-4.3,-1.95)
  {$\textstyle\bigoplus_{i=1}^{n}\pi_{t-s+2i}S^0_{(2)}$};

\node (right) at (4.3,-1.95)
  {$\operatorname{Ext}_{BP_\ast BP}^{s,t}(BP_*\CP^n,BP_*)$};

\node (bottom) at (0,-3.90)
  {$[\CP^n,\Sigma^{s-t}S^0]_{(2)}.$};

\draw[ssarrow] (top) --
  node[sslabel,midway]{Cellwise ANSS} (left);

\draw[ssarrow] (top) --
  node[sslabel,midway]{Algebraic AHSS} (right);

\draw[ssarrow] (left) --
  node[sslabel,midway]{Cohomological AHSS} (bottom);

\draw[ssarrow] (right) --
  node[sslabel,midway]{Adams--Novikov SS} (bottom);
\end{tikzpicture}
\]
Each double arrow depicts a spectral sequence labeled by its name. The upper right double arrow depicts the cohomological \emph{algebraic Atiyah--Hirzebruch spectral sequence (algAHSS)} 
\[
\begin{aligned}
E_1^{2i,s,t}(\Ext(BP_\ast\CP^n))
&=
\Ext_{BP_\ast BP}^{s,t+2i}(BP_*,BP_*)
\Longrightarrow
\Ext_{BP_\ast BP}^{s,t}(BP_*\CP^n,BP_*),
\\
d_{2r}:E_{2r}^{2i,s,t}
&\longrightarrow
E_{2r}^{2(i+r),s+1,t},
\qquad r\geq1.
\end{aligned}
\]
Here $2i$ is the (algebraic) skeletal filtration, $s$ is the Adams-Novikov filtration, and $t$ is the internal degree. For $\alpha\in\Ext_{BP_\ast BP}^{s,t}(BP_*,BP_*)$, we write $\alpha[2i]\in E_1^{2i,s,t-2i}$ for the corresponding element in the algAHSS $E_1$-page with skeletal filtration $2i$, Adams-Novikov filtration $s$ and stem $t-s-2i$. We give a more detailed account of its construction in subsection~\ref{subsec:anss-cobar}.

The goal of our third strategy is to show 
\[\CP^n \xrightarrow{q} S^{2n} \xrightarrow{\beta} S^0\]
is nonzero where $\beta$ is a surviving cycle in the divided beta-families of the sphere, detected by
\[b \in \Ext^{2,2+2n}_{BP_\ast BP}(BP_\ast,BP_\ast)\]
in the ANSS of the sphere. Combining the two spectral sequences on the right side of the $BP$-Mahowald square, we show that
\begin{enumerate}
    \item There is a nonzero class in $\Ext^{2,2}(BP_\ast \CP^n,BP_\ast)$ which is detected by $b[2n]$ in the algAHSS.
    \item This class further survives the ANSS, converging to a nonzero class in $[\CP^n,S^0]$. Furthermore, this class is detected by $\beta[2n]$ in the AHSS.
\end{enumerate}
The first step is entirely an algebraic computation building on the $BP$-formal group law, whereas the second step is a naturality argument. Neither of these steps requires knowing attaching maps longer than $\eta$, which is detected by a nonzero $Sq^2$ in mod-$2$ cohomology. We then seek to extend the $\beta[2n]$'s to $\CP^N$ with $N > n$.

We will elaborate on the mechanism of each of these methods in the remainder of this section. Also, we conclude this section by recollecting the computations of Farrell-Jones \cite{FarrellJones1994hyperbolic} that have shown $\mathbb{CP}^{4k+1}$ admits an exotic smooth structure for $k \geq 1$.

\subsection{\texorpdfstring{$\tmf$}{tmf}-differential method}
\label{sec:tmf-detection}

For a CW complex $X$, the inclusion of a finite skeleton induces a cofiber sequence
\[X^n \to X^N \to X^N_{n+1}\]
which gives rise to maps of AHSSs. A simple diagram chase shows that, if a class $\alpha[n]$ supports a differential in the AHSS of $X_n^N$, then it must not be hit by a differential in the AHSS of $X^n$. Consequently such $\alpha[n]$ survives in the AHSS of $X^n$.

\begin{theorem}\label{thm:tmf-detection}
Let $X$ be a based CW complex, let $E$ be a spectrum, and fix a map $u:S^0\to E$. Choose an $n$-cell of $X$, with $n\geq1$, and let
\[
q:X^{n}\longrightarrow S^n
\]
be the pinch map that collapses $X^{n-1}$ and all the other $n$-cells. For $\alpha\in\pi_nS^0$, put $x=u_*\alpha\in\pi_nE$, and write $x[n]\in E_1^{n,-n}(X;E)$ for the class with coefficient $x$ on the chosen cell and zero on the other $n$-cells.

If $x[n]$ supports a nonzero outgoing $d_r$ on some page then $x[n]$ survives in the $E$-AHSS of $X^{n}$, whereas $\alpha[n]$ survives in the AHSS of $X^{n}$. The same statement also holds in the $p$-local stable homotopy category.
\end{theorem}

\begin{proof}

The survival of $x[n]$ in the $E$-AHSS of $X^{n}$ is equivalent to the composition
\[f \colon X^{n} \xrightarrow{q} S^n \xrightarrow{\alpha} S^0 \xrightarrow{u} E\]
being nonzero. Consider the following diagram where the rows are cofiber sequences.
\begin{center}
    \begin{tikzcd}
X^{n-1} \arrow[r] \arrow[d, equal ] & X^{n} \arrow[r] \arrow[d] & \vee S^n \arrow[d] \\
X^{n-1} \arrow[r]                                & X \arrow[r]                 & X/X^{n-1}.       
\end{tikzcd}
\end{center}
Suppose, for a contradiction, that the composition is zero. Then in the induced long exact sequences of $E$-cohomology, there exists some $g \in E^{-1}(X^{n-1})$ whose image in $E^0(X/X^{n-1})$ maps to $x \in \pi_n E\subset E^0(X^{n}/X^{n-1})$. By naturality, the image of $g$ in $E^0(X/X^{n-1})$ is nontrivial. It is detected by $x[n]$ in the AHSS of $X^{n}/X^{n-1}$, thus all outgoing differentials on $x[n]$ vanish in the AHSS of $X/X^{n-1}$. The induced map
\[E\text{-AHSS of } X/X^{n-1} \longrightarrow E\text{-AHSS of } X\]
and naturality therefore give a contradiction. Consequently, the composition
\[
X^{n} \xrightarrow{q} S^n \xrightarrow{\alpha} S^0
\]
is also nonzero, which is equivalent to the survival of $\alpha[n]$ in the AHSS.
     
\end{proof}

For complex projective spaces with sphere coefficients, the converse of Theorem~\ref{thm:tmf-detection} is also true.

\begin{theorem}\label{finiteinfinitecriterion}
Let $n\geq1$ and $\alpha\in\pi_{2n}S^0$. Then $\alpha[2n] \neq 0$ in $[\Sigma^\infty\CP^n,S^0]$ if and only if $\alpha[2n]$ supports a nonzero outgoing differential in the cohomological AHSS for $\CP_n^\infty$.
\end{theorem}

\begin{proof}
The "if" direction is Theorem~\ref{thm:tmf-detection}, whereas the "only if" direction follows from Nishida's \cite{Nishida1983S1Segal} solution to the $S^1$-Segal conjecture . Namely, Nishida identifies
\[
\varprojlim_N[\Sigma^\infty_+\CP^N,S^0]
\cong A(S^1)\cong\bbZ,
\]
with the integer summand supplied by the disjoint basepoint. Consequently
\[
\varprojlim_N[\Sigma^\infty\CP^N,S^0]
\cong 0.
\]
Now suppose that $\alpha[2n]$ supports no nonzero outgoing differential in $\CP^\infty_n$. Then naturality supplies a homotopy class $f_k \in [\Sigma^\infty \CP^k,S^0]$ for all $k \in \mathbb{N}$ such that 
\[f_k \text{ is detected by } \alpha[2n],~~ i_{k+1}^\ast f_{k+1} = f_k\]
where $i_k: \CP^k \to \CP^{k+1}$ is the inclusion. This provides a nonzero element in the based inverse limit, which is a contradiction.
\end{proof}

\begin{rmk}
Theorem~\ref{finiteinfinitecriterion} holds true integrally (without localizing at any primes). Heuristically, it suggests that in order to find exotic $\CP^N$ for increasing $N$, one has to look for nonzero $\alpha[2\ell]$ with $\ell$ increasing accordingly with $N$. 

Also, Theorem~\ref{finiteinfinitecriterion} guarantees the existence of such a differential in theory. However, in practice it is often not easy to explicitly compute such a differential. Rather, a more practical strategy is: after finding some nonzero $\alpha[2n]$ in the AHSS of $\CP^n$, we extend to $\CP^{n+r}$ where $r$ is usually no greater than $5$. $\CP^{n+r}_n$ involves only $r+1$ cells, hence for small $r$ we can ensure its survival.
\end{rmk}

\subsection{\texorpdfstring{$\tmf$}{tmf}-extension method}\label{sec: delta method}
In this subsection we first describe our strategy for establishing hidden extensions, postponing the auxiliary lemmas until the end. 

We first fix terminology for a family of homotopy classes represented by integral modular forms.

\begin{defn}\label{defn:delta power classes}
    For $n \geq 0$, we refer to the classes inside the 2-local $\pi_\ast\tmf$ detected by
    \[
\begin{gathered}
8\Delta\cdot\Delta^{8n},\quad
4\Delta^2\cdot\Delta^{8n},\quad
8\Delta^3\cdot\Delta^{8n},\quad
2\Delta^4\cdot\Delta^{8n},\\
8\Delta^5\cdot\Delta^{8n},\quad
4\Delta^6\cdot\Delta^{8n},\quad
8\Delta^7\cdot\Delta^{8n},\quad
\Delta^{8n+8}
\end{gathered}
\]
    as \emph{$\Delta$-power classes}.    
\end{defn}
Our second method shows that there exist certain classes $ 0 \neq \alpha[2n] \in \tmf^0(\CP^{N})$ for $N > n$, where $N$ is always divisible by 12. Thus $\alpha[2n]$ is nonzero in $ \tmf^0(\CP^n)$ and $[\CP^n,S^0]$. The main difficulty lies in ruling out the differentials hitting such $\alpha[2n]$. Our strategy divides into the following steps:
\begin{itemize}
    \item Showing hidden 2-extensions 
    \[2^\ell\cdot \alpha[2n] = \delta[2N]\]
    exist in finite quotients $\CP^N_{n'}$ where $n' \leq n \leq N$. In other words, we show that $\tmf^0(\CP^N_{n'})$ has an integral subgroup generated by a homotopy class $A$ detected by $\alpha[2n]$, where $2^\ell \cdot A$ is detected by $\delta[2N]$.
    
    \item We next consider the induced map
    \[q^\ast: \tmf-\text{AHSS of }\CP^N_{n'} \longrightarrow \tmf-\text{AHSS of }\CP^N\]
    of spectral sequences. On the level of spectral sequences, $\delta[2N]$ must map to the element of the same name, where both elements are surviving cycles owing to Lemma~\ref{delta no jump of filtration} .
    
    \item On the level of homotopy groups, $q^\ast A$ hence must be nonzero in $\tmf^0(\CP^N)$ with filtration possibly greater than $2n$. We exclude the possible jump of filtration, which forces the survival of $\alpha[2n]$ in $\tmf^0(\CP^N)$. 
    
    \end{itemize}

In closer detail, there are two types of hidden extensions. We refer to the first type as \emph{unbranched extensions}. The name suggests the phenomenon that all the classes we consider $2$-extend to one another successively. Under such circumstances, Lemma~\ref{delta no jump of filtration} first precludes jump of filtration of the class $\delta[2N]$. We then apply Lemma~\ref{lem:tmf-integral-no-jump-reformulation} to rule out filtration jumps for the classes $\alpha[2n]$ by descending induction on $n$.

The second type of extension, which we call \emph{branched extensions}, is slightly more involved. Namely, there could be two classes $\alpha[2n],\epsilon[2\ell]$ where both $\alpha$ and $\epsilon$ are torsion, whose filtrations satisfy $\ell > n$, yet both extend to $ \delta[2N]$ by $2$:
\[2 \cdot \alpha[2n] = \delta[2N] = 2 \cdot \epsilon[2\ell].\]

On the level of $\tmf$-cohomology groups, the homotopy classes detected by $\alpha[2n],\epsilon[2\ell]$ generate a subgroup consisting of both torsion and integral components. In particular, when displayed by a chart, the classes involved present a branching pattern from the top cell class $\delta[2N]$. 
    
These branching patterns do not satisfy the hypotheses of Lemma~\ref{lem:tmf-integral-no-jump-reformulation}. We use Lemma~\ref{lem:tmf-projective-top-cell-splitting} which provides a splitting of the top cell in the $\tmf$-module $\tmf \wedge \CP^\ell$ to preclude a jump of filtration. Figure~\ref{fig:tmf-hidden-extensions-cp12-cp24} below illustrates hidden extensions in the $\tmf$-AHSS $E_\infty$-page, where the squares represent integral classes, the dots represent $\mathbb{Z}/2 \in \tmf_\ast$, and the concentric circles represent $\mathbb{Z}/4 \in \tmf_\ast$. Dashed lines display hidden extensions. Note that we have omitted classes that do not admit extension to the $\Delta$-power classes from the figure.

\begin{figure}[htbp]
\centering
\begin{tikzpicture}[x=0.60cm,y=0.46cm,font=\small]

  % Left panel: CP^12_7.
  \begin{scope}
    \CPfourgrid{0}{3}{0}{11}

    \foreach \h/\lab in {10/24,8/22,6/20} {
      \node[anchor=east] at (-0.8,\h) {$\lab$};
    }

    % Hidden 2-extensions.
    \draw[CPfour red,dashed] (1,6) -- (1,8);
    \draw[CPfour red,dashed] (1,8) -- (1,10);

    % A dot at filtration 22.
    \fill (1,8) circle (1.5pt);
    \node[right=3pt] at (1,8) {$\bar{\kappa}\eta^2$};

    % A central dot surrounded by a circle at filtration 20.
    \draw[draw=black,fill=white,line width=0.4pt]
      (1,6) circle (3.2pt);
    \fill (1,6) circle (1.5pt);
    \node[right=5pt] at (1,6) {$2\bar{\kappa}$};

    \node[draw,fill=white,rectangle,inner sep=0pt,
          minimum size=4.5pt] at (1,10) {};
    \node[right=3pt] at (1,10) {$8\Delta$};

    \node[anchor=north,align=center] at (1,-1.3) {
      $E_\infty$-page of the\\[2pt]
      $\tmf$-AHSS of $\CP^{12}_7$
    };
  \end{scope}

  % Right panel: CP^24_19.
  \begin{scope}[shift={(8,0)}]
    \CPfourgrid{0}{3}{0}{11}

    \foreach \h/\lab in {10/48,8/46,6/44,4/42,2/40} {
      \node[anchor=east] at (-0.8,\h) {$\lab$};
    }

    % Hidden 2-extensions.
    \draw[CPfour red,dashed] (1,2) -- (1,4);
    \draw[CPfour red,dashed] (1,8) -- (1,10);
    \draw[CPfour red,dashed]
      (1,4)
      .. controls (-0.4,6) and (-0.4,8) ..
      (1,10);

    \foreach \h/\lab in {
      8/{\eta\eta_1\bar{\kappa}},
      4/{\kappa^3},
      2/{\bar{\kappa}^2}
    } {
      \fill (1,\h) circle (1.5pt);
      \node[right=3pt] at (1,\h) {$\lab$};
    }

    \node[draw,fill=white,rectangle,inner sep=0pt,
          minimum size=4.5pt] at (1,10) {};
    \node[right=3pt] at (1,10) {$4\Delta^2$};

    \node[anchor=north,align=center] at (1,-1.3) {
      $E_\infty$-page of the\\[2pt]
      $\tmf$-AHSS of $\CP^{24}_{19}$
    };
  \end{scope}

\end{tikzpicture}
\caption{Some hidden $2$-extensions in the $\tmf$-AHSS of stunted projectives.}
\label{fig:tmf-hidden-extensions-cp12-cp24}
\end{figure}

\begin{lemma}\label{delta no jump of filtration}
Under the induced map
    \[q^\ast: \tmf-\text{AHSS of }\CP^N_{n'} \longrightarrow \tmf-\text{AHSS of }\CP^N\]
of spectral sequences, $\delta[2N]$ must map to the element of the same name, where both elements are surviving cycles.    
\end{lemma}

\begin{proof}
$\delta[2N]$ is a permanent cycle in $\tmf$-AHSS of $\CP^N$. The sources of potential differentials hitting it are of the form $\alpha[2k]$, where $\alpha \in \tmf_{2k+1}$ is torsion. The Leibniz rule excludes such possibilities.
\end{proof}

    \begin{lemma}\label{lem:tmf-integral-no-jump-reformulation}
Suppose a surviving cycle $\beta[2n]$ admits a $2^s$-extension to $\delta[2N]$ in the $\tmf$-AHSS of $\CP^N_n$, with $\delta$ belonging to the $\Delta$-power classes of Definition~\ref{defn:delta power classes}. Suppose also that multiplication by $2^{s-1}$ in $\tmf^0(\CP^N)$ sends every homotopy class of AH filtration greater than $2n$ into filtration $2N$. Then $\beta[2n]$ survives in the $\tmf$-AHSS of $\CP^N$; in other words, the homotopy classes detected by $\beta[2n]$ map nontrivially under
\[q^\ast:\tmf^0(\CP^N_n) \to \tmf^0(\CP^N),\]
where they are detected by the elements with the same AH names. 
\end{lemma}

\begin{proof}

The existence of such a hidden extension amounts to an integral subgroup of $\tmf^0(\CP^N_n)$ whose generator has AH filtration $2n$, which we denote by 1. $2^s \cdot 1$ is detected by $\delta[2N]$, which we denote by $D$. According to Lemma~\ref{delta no jump of filtration} , $0 \neq q^\ast D = 2^s \cdot q^\ast 1$, with AH filtration possibly greater than $2n$. If filtration jump actually occurred, then $2^{s-1} q^\ast 1$ must have AH filtration $2N$ in $\tmf^0(\CP^N)$, which contradicts $2^s q^\ast 1$ being detected by $\delta[2N]$.

\end{proof}

\begin{lemma}\label{lem:tmf-projective-top-cell-splitting}
For every $n\geq1$, $\CP^{8n-1}$ is $\tmf$-orientable. It follows that there is a $\tmf$-fundamental class $[\CP^{8n-1}]_{\tmf} \in \tmf_{16n-2}\CP^{8n-1}$ such that
\begin{equation*}
\bigl\langle q^*(z),[\CP^{8n-1}]_{\tmf}\bigr\rangle=z
\qquad (z\in\pi_{16n-2}\tmf).
\end{equation*}
where $q*:\pi_{16n-2}\tmf\longrightarrow \tmf^0(\CP^{8n-1})$ is induced by the top-cell collapsing map. Furthermore, there is a splitting
\begin{equation*}
\tmf\wedge\CP^{8n-1}
\simeq (\tmf\wedge\CP^{8n-2})\vee\Sigma^{16n-2}\tmf,
\end{equation*}

\end{lemma}

\begin{proof}

Lemma~\ref{lem:tmf-periodic-thom} shows $8nL$ is $\tmf$-orientable. The Euler sequence gives
\[T\CP^{8n-1}\oplus\mathbb C\cong8nL,\]
so the tangent and stable normal bundles are also $\tmf$-orientable, which supplies the fundamental class in $\tmf_{16n-2}(\CP^{8n-1})$ satisfying the desired pairing.

Extending the fundamental class in the category of $\tmf$-modules gives
\[
s_n:\Sigma^{16n-2}\tmf\longrightarrow\tmf\wedge\CP^{8n-1},
\qquad (1\wedge q)s_n=\mathrm{id}.
\]
The splitting thus follows from the section of the quotient map.
\end{proof}

\subsection{The $\beta$-families in algAHSS method}\label{section: ANSS methodology}

We first briefly recall the ANSS and the $BP$-algebraic Atiyah-Hirzebruch spectral sequence (algAHSS). Fix a prime $p$. For bounded-below spectra $X$ and $Y$, the ANSS is
\begin{equation*}
\begin{aligned}
E_2^{s,t}(X)&=\Ext_{BP_\ast BP}^{s,t}(BP_*X,BP_*Y)
                 \Longrightarrow [X,\Sigma^{s-t}Y]_{(p)},\\
d_r:E_r^{s,t}&\longrightarrow E_r^{s+r,t+r-1},\qquad r\geq2.
\end{aligned}
\end{equation*}
The index $s$ is Adams--Novikov filtration, $t$ is internal degree, and $t-s$ is the stem.

For a detailed construction of the algebraic Atiyah--Hirzebruch spectral sequence, we refer the reader to \cite{behrens2006root, wangxu2016algebraic,wangxu2017triviality}. Here we recall only the setup needed for our computations.

We filter $\Ext_{BP_\ast BP}^{s,t}(BP_*\CP^n,BP_*)$ by 
\[F^{2i}\Ext_{BP_\ast BP}^{s,t}(BP_*\CP^n,BP_*)
=
\ker\!\left(
\Ext_{BP_\ast BP}^{s,t}(BP_*\CP^n,BP_*)
\longrightarrow
\Ext_{BP_\ast BP}^{s,t}(BP_*\CP^{i-1},BP_*)
\right).\]
The filtration gives rise to the algebraic Atiyah--Hirzebruch spectral sequence (algAHSS)
\[
\begin{aligned}
E_1^{2i,s,t}(\Ext(BP_\ast\CP^n))
&=
\Ext_{BP_\ast BP}^{s,t+2i}(BP_*,BP_*)
\Longrightarrow
\Ext_{BP_\ast BP}^{s,t}(BP_*\CP^n,BP_*),
\\
d_{2r}:E_{2r}^{2i,s,t}
&\longrightarrow
E_{2r}^{2(i+r),s+1,t},
\qquad r\geq1.
\end{aligned}
\]
where $2i$ is the (algebraic) skeletal filtration, $s$ is the Adams-Novikov filtration, and $t$ is the internal degree. For $\alpha\in\Ext_{BP_\ast BP}^{s,t}(BP_*,BP_*)$, we write $\alpha[2i]\in E_1^{2i,s,t-2i}$ for the corresponding class on the $2i$-cell. If $\alpha[2i]$ survives nontrivially to $E_\infty$, it detects a class of filtration $2i$ in $\Ext_{BP_\ast BP}^{s,t-2i}(BP_*\CP^n,BP_*)$. Just as in the AHSS, the differentials all have even length for degree reasons.

As mentioned in Subsection~\ref{subsec: methods introduction}, we need to first show that, given a nonzero algebraic $\beta$-family element $\beta \in \Ext^{2,2n+2}_{BP_\ast BP}(BP_\ast, BP_\ast)$ in the ANSS of the sphere, $\beta[2n] \in E_1^{2n,2,2}(\Ext(BP_\ast\CP^n))$ survives in the algAHSS. For degree reasons, the only possible sources of algAHSS differentials that could hit $\beta[2n]$ have tridegree $(2\ell,1,2)$ and hence take the form
\[
\alpha[2\ell], \qquad \alpha\in \Ext^{1,2\ell+2}_{BP_\ast BP}(BP_\ast,BP_\ast), \qquad 1\leq\ell<n.
\]
For every such candidate $\alpha[2\ell] \in E_1^{2\ell,1,2}(\Ext(BP_\ast\CP^n))$, we show that
\begin{enumerate}
    \item Either $\alpha[2\ell]$ supports a nonzero $d_2$-algAH differential
    \[
    d_2(\alpha[2\ell])=\ell\alpha_1\alpha[2\ell+2],
    \]
    whose target is not $\beta[2n]$;
    \item or $\alpha[2\ell]$ survives to $E_\infty$ in the algAHSS.
\end{enumerate}
The relevant $d_2$-kernels are generated by $\alpha_{\ell+1}[2\ell]$ for even $\ell\geq2$, $g_2[2]$, and $2g_{\ell+1}[2\ell]$ for odd $\ell\geq3$, within the range $1\leq\ell<n$.

By naturality, if we can show that all such $d_2$-cycles survive the algAHSS of $\Ext_{BP_\ast BP}(BP_\ast \CP^\infty, BP_\ast)$, then they will also survive the algAHSS of $\Ext_{BP_\ast BP}(BP_\ast \CP^n, BP_\ast)$. We do this in two steps:
\begin{itemize}
    \item We first use the logarithm of the universal $2$-typical formal group to construct classes in $\Ext^{0,2}_{BP_\ast BP}(BP_\ast \CP^\infty, (BP_\ast)_{\mathbb Q}/BP_\ast)$.
    \item We then show that their images under the rational Bockstein
    \[
    \delta_{\mathbb Q}:
    \Ext^{0,2}_{BP_\ast BP}
    \bigl(BP_\ast\CP^\infty,(BP_\ast)_{\mathbb Q}/BP_\ast\bigr)
    \longrightarrow
    \Ext^{1,2}_{BP_\ast BP}
    \bigl(BP_\ast\CP^\infty,BP_\ast\bigr)
    \]
    have the desired algAHSS names in $\Ext^{1,2}_{BP_\ast BP}(BP_\ast\CP^\infty,BP_\ast)$.
\end{itemize}
The Bockstein is the boundary map in the long exact sequence of $\Ext$-groups, induced by the short exact sequence of $BP_\ast BP$-comodules
\[
0 \to BP_\ast \to (BP_\ast)_{\mathbb{Q}} \to (BP_\ast)_{\mathbb Q}/BP_\ast \to 0.
\]

By construction, suppose a zero-cocycle in $\Ext^{0,2}_{BP_\ast BP} \bigl(BP_\ast\CP^\infty,(BP_\ast)_{\mathbb Q}/BP_\ast\bigr)$ is represented by $\overline H$ in the cobar complex. It admits a rational lift $H$ whose coboundary $d^0H$ is integral. This coboundary represents $\delta_{\mathbb Q}(\overline H)$.

The key point is that the Bockstein respects the algAHSS filtration. We use the power-series description
\[
\mathcal C_{BP_\ast BP}^{s,t}(BP_\ast\CP^\infty,R)
\cong
\left\{
\sum_{j\geq1}u_jz^j:
u_j\in
\bigl((BP_\ast BP)^{\otimes_{BP_\ast}^s}
      \otimes_{BP_\ast}R\bigr)_{t+2j}
\right\},
\]
where $R$ is one of the coefficient comodules $BP_\ast$, $(BP_\ast)_{\mathbb Q}$, or $(BP_\ast)_{\mathbb Q}/BP_\ast$. Here the zeroth tensor power is $BP_\ast$. The algAHSS filtration is given by divisibility by powers of $z$, with $z^\ell$ in filtration $2\ell$. The Bockstein respects this filtration and acts on the associated graded through the coefficients. Explicitly, for a rational zero-cochain $H(z)$
\[
H(z)=az^\ell+O(z^{\ell+1})
\]
of internal degree $2$, with $a\in((BP_\ast)_{\mathbb Q})_{2\ell+2}$,
\[
d^0H(z)=(\eta_R(a)-\eta_L(a))z^\ell+O(z^{\ell+1}).
\]
If $d^0H$ is integral and $\delta_{\mathbb Q}(\bar a)\ne0$, Proposition~\ref{prop: algAHSS by cobar filtration} identifies the algAHSS name of $[d^0H]$ as $\delta_{\mathbb Q}(\bar a)[2\ell]$, with filtration exactly $2\ell$.

Thus, to lift a $1$-line class, we seek a rational power series whose
\begin{itemize}
    \item coboundary is integral;
    \item leading coefficient maps to the desired class under the coefficient Bockstein.
\end{itemize}

The reduced algebraic $S^1$-transfer provides a model for this construction. Consider the short exact sequence of $BP_\ast BP$-comodules:
\[
0 \to \Sigma^{-2}BP_\ast \to BP_\ast(\CP^\infty_{-1}) \to BP_\ast(\CP^{\infty}_+) \to 0.
\]
Yoneda composition with its extension class induces the associated connecting homomorphisms. We denote the restriction of this class to the reduced summand by
\[
e_\tau\in
\Ext^{1,2}_{BP_\ast BP}(BP_\ast\CP^\infty,BP_\ast).
\]
Miller's description \cite{Miller1982Bernoulli} gives
\[
e_\tau=\delta_{\mathbb Q}(\overline{P_F}),
\qquad
P_F(z)=\frac{1}{\lambda_F(z)}-\frac{1}{z}+\lambda_1,
\]
where $\lambda_F$ is the logarithm of the universal $2$-typical formal group law. Under the power-series description of the cobar complex, $P_F(z)=\sum_{j\geq1}c_jz^j$ corresponds to the $BP_\ast$-linear map
\[
\Phi_F:BP_\ast\CP^\infty
\longrightarrow\Sigma^{-2}(BP_\ast)_{\mathbb Q},
\qquad b_j\longmapsto c_j,
\]
where $b_j$ is dual to the $j$th power of the chosen complex orientation $x \in BP^2(\CP^\infty)$. We extend $\Phi_F$ to $BP_\ast(\CP^\infty_+)$ by setting $\Phi_F(b_0)=0$. The reduced transfer has algAHSS name $g_2[2]$, as an explicit computation gives
\[
P_F(z)=\frac{v_1^2}{4}z+O(z^2).
\]

The series $P_F$, which we call the reciprocal logarithm, thus supplies a natural starting point. Write $\lambda_F(z)=z+\sum_{i\geq1}\lambda_i z^{2^i}$. Hazewinkel's recursion controls the denominators of the $c_j$. For even $\ell\geq2$, multiplying $P_F$ by $2^\ell$ makes all terms of degree less than $\ell$ integral. Subtracting these terms leaves a series whose leading coefficient is congruent to $-v_1^{\ell+1}/2$ modulo $BP_\ast$, yielding a cobar lift of the algAHSS $d_2$-cycle $\alpha_{\ell+1}[2\ell]$.

For the remaining lifts, we compose the reciprocal logarithm with the $a$-series $[a]_F$. Miller's transfer calculation identifies
\[
-\Phi_F(b_1^\ell)=\frac{B_{\ell+1}(F)}{\ell+1},
\]
where $b_1^\ell$ denotes the $\ell$th Pontryagin power of $b_1$ \cite{Miller1982Bernoulli}. For odd $\ell\geq3$,
\[
\delta_{\mathbb Q}
\left(\overline{\frac{B_{\ell+1}(F)}{\ell+1}}\right)
\]
equals $2g_{\ell+1}$ up to an odd unit. In Proposition~\ref{prop:divided}, we construct a rational series $G_\ell(z)$ as a linear combination of $P_F([a]_F(z))$, with integral coboundary and leading term $-\Phi_F(b_1^\ell)z^\ell$. For odd $\ell\geq3$, these series yield lifts such as $2g_{\ell+1}[2\ell]$. Thus every $d_2$-cycle in tridegree $(2\ell,1,2)$, with $1\leq\ell<n$, survives in the algAHSS converging to $\Ext^{\ast,\ast}_{BP_\ast BP}(BP_\ast\CP^n,BP_\ast)$. Consequently, none can support a higher differential hitting $\beta[2n]$.

Next we show that $\beta[2n]$ is a permanent cycle in the ANSS. The quotient $q \colon \CP^n \to S^{2n}$ induces maps between ANSSs
\[\Ext^{s,t+2n}_{BP_\ast BP}(BP_\ast, BP_\ast) \xrightarrow{q^\ast}  \Ext^{s,t}_{BP_\ast BP}(BP_\ast\CP^n, BP_\ast),\]
under which the image of $\beta \in \Ext_{BP_\ast BP}^{2,2n+2}$ is $\beta[2n] \in \Ext_{BP_\ast BP}^{2,2}(BP_\ast\CP^n,BP_\ast)$. The permanence of $\beta[2n]$ follows from naturality.
    
Finally, we show that the homotopy class in $[\CP^n,S^0]$ detected by $\beta[2n]$ has AHSS name $b[2n]$. In other words, let $b\in \pi_{2n}S^0$ be a beta-family element detected by $\beta$. Then this homotopy class is detected by $b[2n]$ in the AHSS of $\CP^n$. We seek to extend $b[2n]$ further in the AHSS of $\CP^{n+r}$. If this extension exists, Theorem~\ref{AHSSsmoothing} will supply exotic complex projective spaces in dimensions $n$ through $n+r$.

\subsection{\texorpdfstring{$\ko$}{ko} detection}
\label{sec:ko-detection}

A first application of coefficient detection comes from real $K$-theory.

\begin{theorem}[Farrell--Jones {\cite{FarrellJones1994hyperbolic}}]
\label{thm:ko-farrell-jones}
For every $k\geq1$, there is a homotopy sphere $\Sigma^{8k+2}\in\Theta_{8k+2}$ such that $\CP^{4k+1}\mathbin{\#}\Sigma^{8k+2}$ is not concordant to standard $\CP^{4k+1}$. Hence there are exotic complex projectives in dimension $4k+1$ for $k \geq 1$.
\end{theorem}
\begin{proof}
Let $h_{\KO}:\pi_*S^0\to\KO_*$ be the unit map on homotopy. The Adams $\mu$-element $\mu_{8k+2}\in\pi_{8k+2}S^0$ has order two. Its image under $h_{\KO}$ is the generator of
\[
\KO_{8k+2}=\widetilde\KO^0(S^{8k+2})\cong\bbZ/2.
\]
We denote this generator by $\beta^k\eta^2$. Here $\beta \in \KO_8$ is the Bott class and $\eta \in \KO_1 = \mathbb{Z}/2\mathbb{Z}$ is the generator. By Wood's theorem there is an equivalence of $\KO$-modules:
\[F(\CP^{4k+1},\KO) \cong \bigvee_{i=1}^{2k} \Sigma^{-4i}\mathrm{KU} \vee \Sigma^{-8i-2}\KO.\]
Hence in the $\KO$-AHSS $\beta^k\eta^2[8k+2]$ is a surviving cycle, which detects an exotic $\CP^{4k+1}$ by Theorem~\ref{AHSSsmoothing} and Theorem~\ref{thm:tmf-detection}. In particular, the quotient map $q: \CP^{4k+1} \to S^{8k+2}$ induces
\[ q^\ast: \theta_{8k+2} \cong \pi_{8k+2}(Top/O) \to C(\CP^{4k+1}),\]
under which the exotic sphere $\Sigma^{8k+2}$ corresponding to $\mu_{8k+2}$ maps to $\CP^{4k+1} \# \Sigma^{8k+2}$, where the first $\CP^{4k+1}$ is the standard smooth complex projective space.

\end{proof}

\begin{rmk}\label{rmk:ko-family-scope}
Alternatively, owing to the Steenrod operation $\operatorname{Sq}^2(t^{4k+1})=t^{4k+2}$ in $H^*(\mathbb{CP}^\infty;\mathbb{F}_2) \cong \mathbb{F}_2[t]$ one has
\[
\CP_{4k+1}^{4k+2}
 \simeq \Sigma^{8k+2} C(\eta).
\]
It follows that there is a $d_2$-differential
\[
d_2\bigl(\mu_{8k+2}[8k+2]\bigr)
 =(\eta\mu_{8k+2})[8k+4]\ne0
\]
in the cohomological AHSS of $\mathbb{CP}^\infty_{4k+1}$. Hence by Theorem~\ref{finiteinfinitecriterion}, $\mu_{8k+2}[8k+2]$ survives in the AHSS of $\mathbb{CP}^{4k+1}$ and supplies an exotic smooth structure. 
\end{rmk}

Each nonconcordant connected sum in Theorem~\ref{thm:ko-farrell-jones} represents a nonzero smoothing class. By Theorem~\ref{concordanceCP}, its underlying manifold is homeomorphic but not diffeomorphic to standard projective space. In particular, this gives exotic $\CP^{4k+1}$ for $k\geq1$.

\section{Computations of \texorpdfstring{$\tmf$}{tmf}-AHSS differentials}
\label{sec: tmf-differential method}

All spectra and groups in this section are localized at $2$.

\begin{theorem}
\label{thm:tmf-primarymain}
The spaces $\CP^{m+96k}$ admit exotic smooth structures for
$$
m\in\{3,4,6,7,27,28,51,52,54,55,56,75,76\},\qquad k\geq0,
$$
except for $(m,k)=(3,0)$.
\end{theorem}

We use the notations of \cite[Table 3]{isaksenkongliruanzhu2024ANSStmf}. Recall the classes $\nu\in\pi_3{\tmf_{(2)}}$, $\kappa\in\pi_{14}{\tmf_{(2)}}$, and
$$
\nu_j=\{\Delta^j\nu\}\in\pi_{24j+3}{\tmf_{(2)}}\quad j=2,4,6,
\qquad \kappa_4=\{\Delta^4\kappa\}\in\pi_{110}{\tmf_{(2)}}.
$$
The braces specify the classes represented by the indicated modular detectors. Set $\nu_0=\nu$. The required products are
$$
\begin{gathered}
\eta\kappa\neq0,\qquad \eta\kappa_4\neq0,\\
\eta(\nu\nu_j)=0,\qquad \nu^2\nu_j\neq0
\quad j=0,2,4,6,
\end{gathered}
$$
which hold after multiplication by $\Delta^{8k}$ on sources and targets.

\begin{theorem}
\label{thm:tmf-primary}
The spaces $\CP^{m+96k}$ admit exotic smooth structures for
$$
m\in\{3,7,27,51,55,75\},\qquad k\geq0,
$$
except for $(m,k)=(3,0)$.
\end{theorem}

\begin{proof}
For $m=7,55$, let $x$ be $\kappa$ or $\kappa_4$ respectively. Put $N=m+96k$. Since $N$ is odd, there exists
$$
d_2\bigl(\Delta^{8k}x[2N]\bigr)=\Delta^{8k}\eta x[2N+2]\neq0.
$$
The existence of exotic $\CP^N$ follows from Theorem~\ref{AHSSsmoothing} and Theorem~\ref{thm:tmf-detection}.

For $m=3,27,51,75$, use $x=\nu\nu_j$ with $j=0,2,4,6$, respectively. Put $N=m+96k$. $\Sigma^{-2N}\CP^{N+2}_N$ is equivalent to
\[\text{Cof}(S^3 \vee S^1 \xrightarrow{(\nu,\eta)} S^0)\]
Consequently
\[d_2(\Delta^{8k}x[2N]\bigr)=\eta \cdot \Delta^{8k}x[2N+2] = 0,\]
$$
d_4\bigl(\Delta^{8k}x[2N]\bigr)=\Delta^{8k}\nu x[2N+4]\neq0.
$$
where $\Delta^{8k}\nu x[2N+4]$ is not a target of a $d_2$-differential as $Sq^2(t^{N+1}) = 0$ in $H^*(\mathbb{CP}^\infty;\mathbb{F}_2) \cong \mathbb{F}_2[t]$. The existence of exotic $\CP^N$ follows from Theorem~\ref{AHSSsmoothing} and Theorem~\ref{thm:tmf-detection}.
\end{proof}

\begin{theorem}\label{thm:tmf-eta-extensions}
For every $k\geq0$, $\CP^{56+96k}$ admits an exotic smooth structure.
\end{theorem}
\begin{proof}
According to the splitting of Lemma~\ref{lem:tmf-projective-top-cell-splitting}, $2\Delta^{8k}\kappa_4[110+192k]$ detects the nonzero homotopy class
\[\CP^{55+96k} \xrightarrow{q} S^{110+192k} \xrightarrow{2\Delta^{8k}\kappa_4} \tmf.\] 
Choose a lift $\sigma_k$ of \cite{BQeta2026}, with $u_*\sigma_k=2\Delta^{8k}\kappa_4$ and $\eta\sigma_k=0$. Theorem~\ref{thm:tmf-detection} shows $\sigma_k[110+192k]$ is a surviving cycle in the AHSS of $[\CP^{55+96k},S^0]$.

In the AHSS of $\CP_{55+96k}^\infty$,
\[d_2(\sigma_k[110+192k]) = \eta \sigma_k[112+192k] = 0.\]
Hence $\sigma_k[110+192k]$ extends over $\CP^{56+96k}$. Theorem~\ref{AHSSsmoothing} now gives the required exotic smooth structure.
\end{proof}

\begin{theorem}\label{thm:tmf-doubled}
The spaces $\CP^{28+96k}$ and $\CP^{76+96k}$ admit exotic smooth structures for every $k\geq0$.
\end{theorem}
\begin{proof}
Let $j\in\{2,6\}$, $k\geq0$, and $m=12j+3+96k$. In the $\tmf_{(2)}$--AHSS of $\CP_m^{m+3}$, one has
\[
d_6\bigl(2\Delta^{8k}\nu\nu_j[2m]\bigr)
 =\Delta^{8k}\epsilon\nu_j[2m+6]\ne0.
\]

For $m\equiv3\pmod8$, a cell-diagram depiction of $\Sigma^{-2m}\CP_m^{m+3}$ is the following:
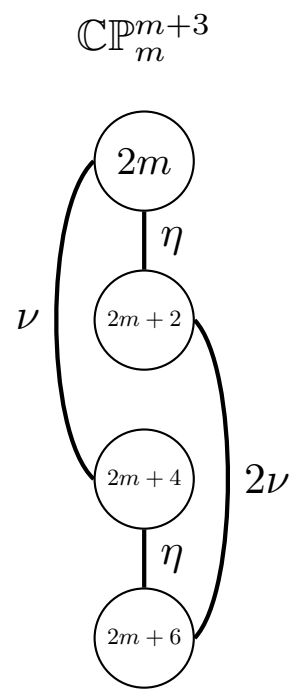
\[
\begin{tikzpicture}[x=.9cm,y=.95cm,scale=.82,transform shape,
 cell/.style={circle,draw,fill=white,minimum size=8mm,
              inner sep=.4pt,font=\small},
 eta/.style={black,line width=.8pt},
 nuedge/.style={black,line width=.8pt},
 twonuedge/.style={black,line width=.8pt}]
\node[cell] (c0) at (0,0) {$2m$};
\node[cell] (c2) at (0,-1.35) {\resizebox{6mm}{!}{$2m+2$}};
\node[cell] (c4) at (0,-2.7) {\resizebox{6mm}{!}{$2m+4$}};
\node[cell] (c6) at (0,-4.05) {\resizebox{6mm}{!}{$2m+6$}};
\draw[eta] (c0)--node[right]{$\eta$}(c2);
\draw[eta] (c4)--node[right]{$\eta$}(c6);
\draw[nuedge] (c0.west)
 .. controls (-.9,-.4) and (-.9,-2.3)
 .. node[midway,left]{$\nu$}(c4.west);
\draw[twonuedge] (c2.east)
 .. controls (.85,-1.75) and (.85,-3.65)
 .. node[midway,right]{$2\nu$}(c6.east);
\node at (0,1.05) {$\CP_m^{m+3}$};
\end{tikzpicture}
\]
The sources and targets of the proposed differentials survive to the $E_6$-page. In $\pi_*\tmf_{(2)}$, \cite[Appendix A]{BauerTmf2008} leads to the following matrix Toda bracket
\[
\left\langle2\nu,(\eta,\nu),
       \begin{pmatrix}2\nu\\\eta\end{pmatrix}\right\rangle
 =\{0+\epsilon\}=\{\epsilon\}.
\]
The Atiyah-Hirzebruch differential
\[d_6(2\nu[2m]) = \epsilon[2m+6]\]
thus follows from a slight generalization of \cite{wangxu2017triviality}. The proposed differentials follow by applying the Leibniz rule. 

By Theorem~\ref{thm:tmf-detection}, for any lift $\gamma_k$ of $\Delta^{8k}\nu\nu_2$, respectively $\Delta^{8k}\nu\nu_6$, $ 0 \neq 2\gamma_k[2m] \in [\CP^{m},S^0]$. $2\gamma_k$ is annihilated by $\eta$, so it extends to $[\CP^{m+1},S^0]$. Theorem~\ref{AHSSsmoothing} hence provides the desired exotic complex projectives.
\end{proof}

We next show that $\CP^{4+96k}$ and $\CP^{52+96k}$ admit exotic smooth structures. The proof makes use of a ``zigzags'' argument. Namely, we first prove there are differentials
\[d_8(\epsilon[8]) = d_8(B[8]) = \kappa\eta[16] \qquad \text{in} \qquad \mathbb{HP}^4_2,\]
where $B\in\pi_8\tmf$ is the Bott class as in Bruner--Greenlees--Rognes \cite{BGRlocal2021}. In the classical Adams spectral sequence $B$ is detected by $w_1$ in filtration $4$.

An $H\mathbb{F}_2$-quotient map $\widetilde p: \CP^9_4 \to \mathbb{HP}^4_2$ together with the "total differential"
\[\delta(\kappa[14]) = \kappa\eta[16]+\kappa\nu[18] \qquad \text{in} \qquad \CP^9_4\]
would supply the $d_{10}$-differential 
\[d_{10}(\epsilon[8]) = \kappa\nu[18]\]
in $\CP^{9}_4$, hence by Theorem~\ref{thm:tmf-detection}, $\CP^4$ admits an exotic smooth structure.

\begin{theorem}\label{thm: cp4}
For every $k\ge0$, in the $\tmf$--AHSS of $\CP_{4+96k}^{9+96k}$, there are $d_{10}$-differentials
\[
d_{10}(\Delta^{8k}\epsilon[8+192k])=d_{10}(\Delta^{8k}B[8+192k])=\Delta^{8k}\kappa\nu[18+192k]\ne0.
\]    
Hence $\CP^{4+96k}$ admits exotic smooth structures for every $k\geq0$.    
\end{theorem}
\begin{proof}
% Once we have shown in the $\tmf$--cohomology AHSS of $\CP_4^9$
It suffices to prove the following differential for $k=0$ 
\[
d_{10}(\epsilon[8])=d_{10}(B[8])=\kappa\nu[18]\ne0.
\]
The proposed differential immediately follows from the $\tmf$-module equivalence
\begin{equation*}
\tmf\wedge\CP_{4+96n}^{9+96n}
\simeq\Sigma^{192n}(\tmf\wedge\CP_4^9)
\end{equation*}
of Lemma~\ref{lem:tmf-periodic-thom}, combining Leibniz rule. The existence of exotic smooth structure would readily follow from Theorem~\ref{AHSSsmoothing} and Theorem~\ref{thm:tmf-detection}.
\end{proof}

\begingroup
\emergencystretch=3em

\def\HP{\mathbb{HP}}
\def\F{\mathbb F_2}
\def\Z{\mathbb Z_{(2)}}
\def\bk{\bar\kappa}
\def\br#1{\langle #1\rangle}

\colorlet{nublue}{black} \colorlet{twonu}{black}

\tikzset{
  cell/.style={ circle,draw,fill=white, minimum size=8mm,inner sep=.4pt,font=\small }, etaedge/.style={black,line width=.8pt}, nuedge/.style={nublue,line width=.8pt}, twoedge/.style={twonu,line width=.8pt}, map/.style={black,-{Stealth[length=2mm]},line width=.8pt}
}

\def\cells#1#2#3{%
  \foreach \d in {#3}{%
    \node[cell] (#1\d) at (#2,{-(\d-8)*.675}) {$\d$};%
  }%
}

\def\etaedge#1#2{%
  \draw[etaedge] (#1)--
    node[right]{$\eta$}(#2);%
}

\def\nuedge#1#2#3{%
  \ifdim#3pt<0pt \draw[nuedge] (#1.west) .. controls ($(#1)+(#3,-.4)$) and ($(#2)+(#3,.4)$) .. node[midway,left]{$\nu$}(#2.west); \else \draw[nuedge] (#1.east) .. controls ($(#1)+(#3,-.4)$) and ($(#2)+(#3,.4)$) .. node[midway,right]{$\nu$}(#2.east); \fi
}

\def\twoedge#1#2#3{%
  \ifdim#3pt<0pt \draw[twoedge] (#1.west) .. controls ($(#1)+(#3,-.4)$) and ($(#2)+(#3,.4)$) .. node[midway,left]{$2\nu$}(#2.west); \else \draw[twoedge] (#1.east) .. controls ($(#1)+(#3,-.4)$) and ($(#2)+(#3,.4)$) .. node[midway,right]{$2\nu$}(#2.east); \fi
}

Now we elaborate on the proof that
\[d_{10}(\epsilon[8])=d_{10}(B[8])=\kappa\nu[18]\ne0.\]
By Wang--Xu \cite[Definition~4.1]{wangxu2017triviality}, an $H\F$--subcomplex is a map injective on mod-$2$ homology, and an $H\F$--quotient is a map surjective on mod-$2$ homology. A cofiber of the former is the latter; the fiber of the latter is the former \cite{wangxu2017triviality}.

\begin{prop}\label{prop:CP to HP}
    There is an $H\mathbb{F}_2$-quotient $\widetilde p: \CP^9_4 \to \mathbb{HP}_2^4$ where $\mathbb H$ is the real quaternion division algebra.
\end{prop}
\begin{proof}
We start from an \emph{unstable} map of spaces. Regard $\mathbb H^5$ as the complex vector space $\mathbb C^{10}$ by restricting right scalars along $\mathbb C=\mathbb R\{1,\mathbf i\} \subset\mathbb H$. The twistor projection is
\[
p:\CP^9\longrightarrow\HP^4,\qquad
[L]\longmapsto L\mathbb H.
\]
Here $[L]$ denotes the projective point represented by the complex line $L\subset\mathbb C^{10}$. Its image $L\mathbb H$ is the quaternionic span of that line, viewed as a point of $\HP^4$.

For its effect on cohomology, let $L$ be the tautological complex line on $\CP^9$ and $H$ the tautological quaternionic line on $\HP^4$. Then $p^*H\cong L\oplus \overline{L}$ where $\overline{L}$ is the conjugate complex line bundle, with $c_1(\overline{L})=-c_1(L)$. By the Whitney sum formula $c(E\oplus F)=c(E)c(F)$, we therefore have
\[
p^*c_2(H)
=c_1(L)c_1(\overline{L})
=-c_1(L)^2.
\]
Consequently $p^\ast$ sends the generator of $H^4(\mathbb{HP}^4)$ to the generator of $H^4(\CP^9)$ up to a unit.

Choosing the standard coordinate inclusion $\mathbb H^2\subset\mathbb H^5$ identifies $\CP^3\subset\CP^9$ and $\HP^1\subset\HP^4$, and $p$ restricts to their twistor projection. It therefore induces a based map of quotient spaces
\[
\widetilde p:\CP^9/\CP^3\longrightarrow\HP^4/\HP^1.
\]
By an abuse of notation, we also denote its $2$-local stabilization by $\widetilde p$, which is our $H\mathbb{F}_2$-quotient map.    
\end{proof}

The bottom-cell inclusion induces an isomorphism $\pi_{10}S^8\cong\pi_{10}\mathbb H\bbP_2^4$. Set $j: S^{10} \to \CP^9_4$ to be the composition
\[S^{10} \xrightarrow{\langle 0,id\rangle} S^8 \vee S^{10} \cong \CP^5_4 \to \CP^9_4 \]

Then $\widetilde p j\simeq0$, while $j$ still induces an isomorphism on $H^{10}$. Set $X=\operatorname{cofib}(j)$. A choice of the indicated nullhomotopy induces a map $X\to\mathbb H\bbP_2^4$, which we also denote by $\widetilde p$; the domain distinguishes it from the map on $\CP_4^9$. Thus there is a commuting diagram
\[
\begin{tikzcd}[column sep=large]
S^{10}\ar[r,"j"] & \CP_4^9\ar[r,"q"]\ar[d,"\widetilde p"] & X\ar[dl,"\exists\,\widetilde p"]\\
& \mathbb H\bbP_2^4 &
\end{tikzcd}
\]
Both $q$ and $\widetilde p$ are $H\F$--quotient maps. The quotient $X$ has cells in dimensions $8,12,14,16,18$.

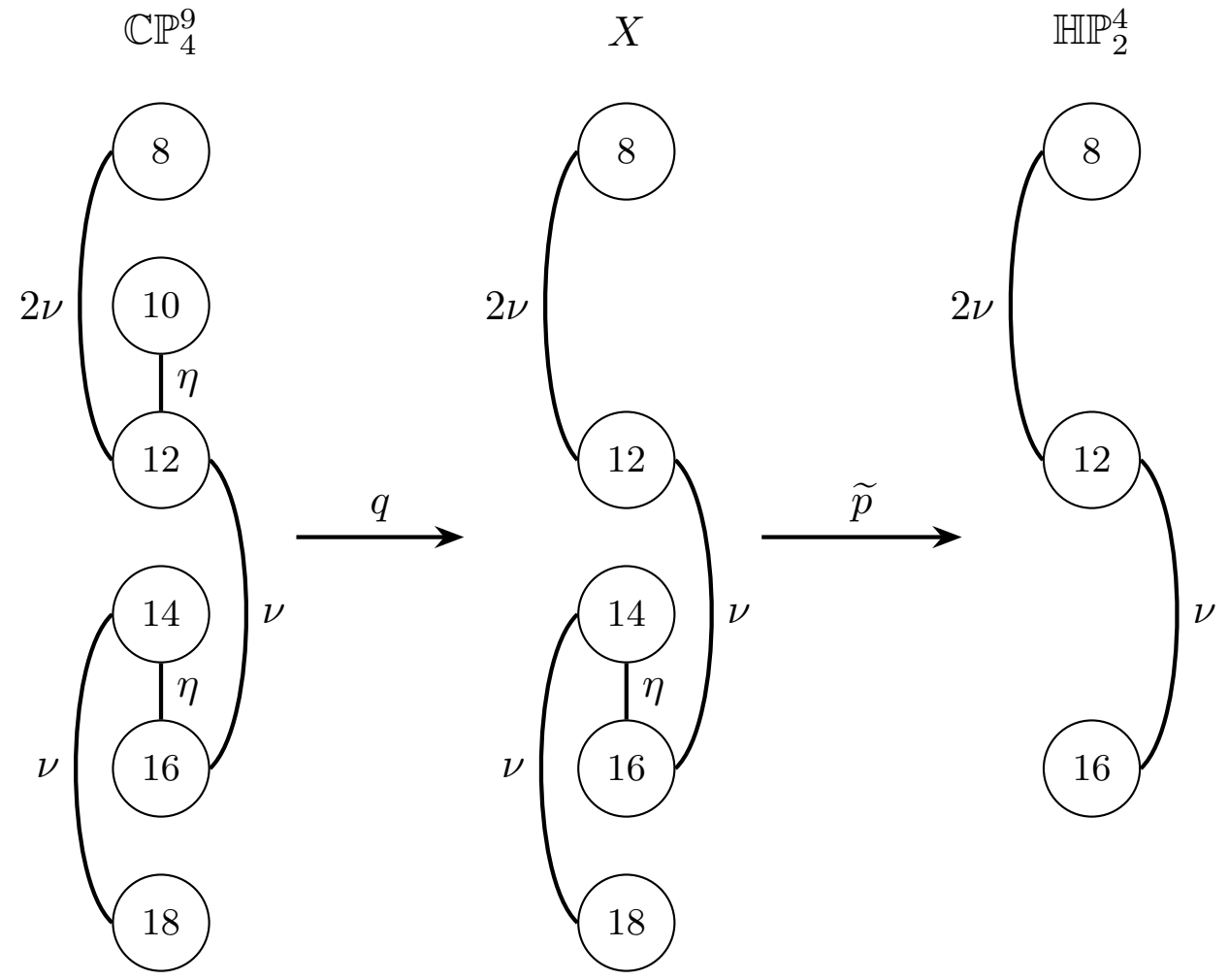
\begin{figure}[htbp]\centering
\begin{tikzpicture}[x=.9cm,y=.95cm,scale=.82,transform shape]
\cells{p}{0}{8,10,12,14,16,18}
\cells{x}{4.3}{8,12,14,16,18}
\cells{y}{8.6}{8,12,16}
\node at (0,1.05) {$\CP_4^9$};
\node at (4.3,1.05) {$X$};
\node at (8.6,1.05) {$\mathbb H\bbP_2^4$};
\etaedge{p10}{p12}\etaedge{p14}{p16}
\twoedge{p8}{p12}{-.85}\nuedge{p12}{p16}{.9}
\nuedge{p14}{p18}{-.9}
\etaedge{x14}{x16}\twoedge{x8}{x12}{-.85}\nuedge{x12}{x16}{.9}
\nuedge{x14}{x18}{-.9}
\twoedge{y8}{y12}{-.85}\nuedge{y12}{y16}{.9}
\draw[map] (1.25,-3.375)--node[above]{$q$}(2.8,-3.375);
\draw[map] (5.55,-3.375)--node[above]{$\widetilde p$}(7.4,-3.375);
\node[align=center,font=\small] at (0,-7.8) {};
\node[align=center,font=\small] at (4.3,-7.8) {};
\node[align=center,font=\small] at (8.6,-7.8) {};
\end{tikzpicture}
\caption{Actual maps underlying the subquotient argument. }\label{fig:cp9-maps}
\end{figure}

\begin{lemma}\label{lem:cp9-fiber}
For $\widetilde p:X\to\mathbb H\bbP_2^4$, there is a cofiber sequence
\[
\operatorname{fib}(\widetilde p)\longrightarrow X
\xrightarrow{\widetilde p}\mathbb H\bbP_2^4\longrightarrow\Sigma \operatorname{fib}(\widetilde p),
\]
with $\operatorname{fib}(\widetilde p)\simeq\Sigma^{14}C\nu$.
\end{lemma}
\begin{proof}
The integral homology sequence gives $H_*(\operatorname{fib}(\widetilde p);\Z)=\Z$ in dimensions $14,18$ and zero elsewhere. The short exact sequence of mod-$2$ cohomology modules is
\[
0\longrightarrow H^*\mathbb H\bbP_2^4\longrightarrow H^*X\longrightarrow H^*\operatorname{fib}(\widetilde p)
\longrightarrow0.
\]
Denote the generators of $\operatorname{fib}(\widetilde p)$'s mod-2 cohomology by $e^{14}$ and $e^{18}$ in respective dimensions. Naturality supplies a nonzero $\Sq^4$-operation, namely $\Sq^4(e^{14}) =e^{18}$. Consequently,
\[
\operatorname{fib}(\widetilde p)\simeq\Sigma^{14}C\nu.
\]
\end{proof}

\begin{figure}[htbp]\centering
\begin{tikzpicture}[x=.9cm,y=.95cm,scale=.82,transform shape]
\cells{f}{0}{14,18}\cells{xx}{4.3}{8,12,14,16,18}\cells{yy}{8.6}{8,12,16}
\node at (0,1.05) {$\Sigma^{14}C\nu$};\node at (4.3,1.05) {$X$};\node at (8.6,1.05) {$\mathbb H\bbP_2^4$};
\nuedge{f14}{f18}{.9}
\etaedge{xx14}{xx16}\twoedge{xx8}{xx12}{-.85}\nuedge{xx12}{xx16}{.9}\nuedge{xx14}{xx18}{-.9}
\twoedge{yy8}{yy12}{-.85}\nuedge{yy12}{yy16}{.9}
\draw[map] (1.4,-5.4)--node[above]{$i$}(3.15,-5.4);
\draw[map] (5.9,-2.7)--node[above]{$\widetilde p$}(7.4,-2.7);
\end{tikzpicture}
\caption{The cofiber sequence where $\mathbb H\bbP_2^4$ is an $H\F$--quotient of $X$.}\label{fig:cp9-fiber}
\end{figure}
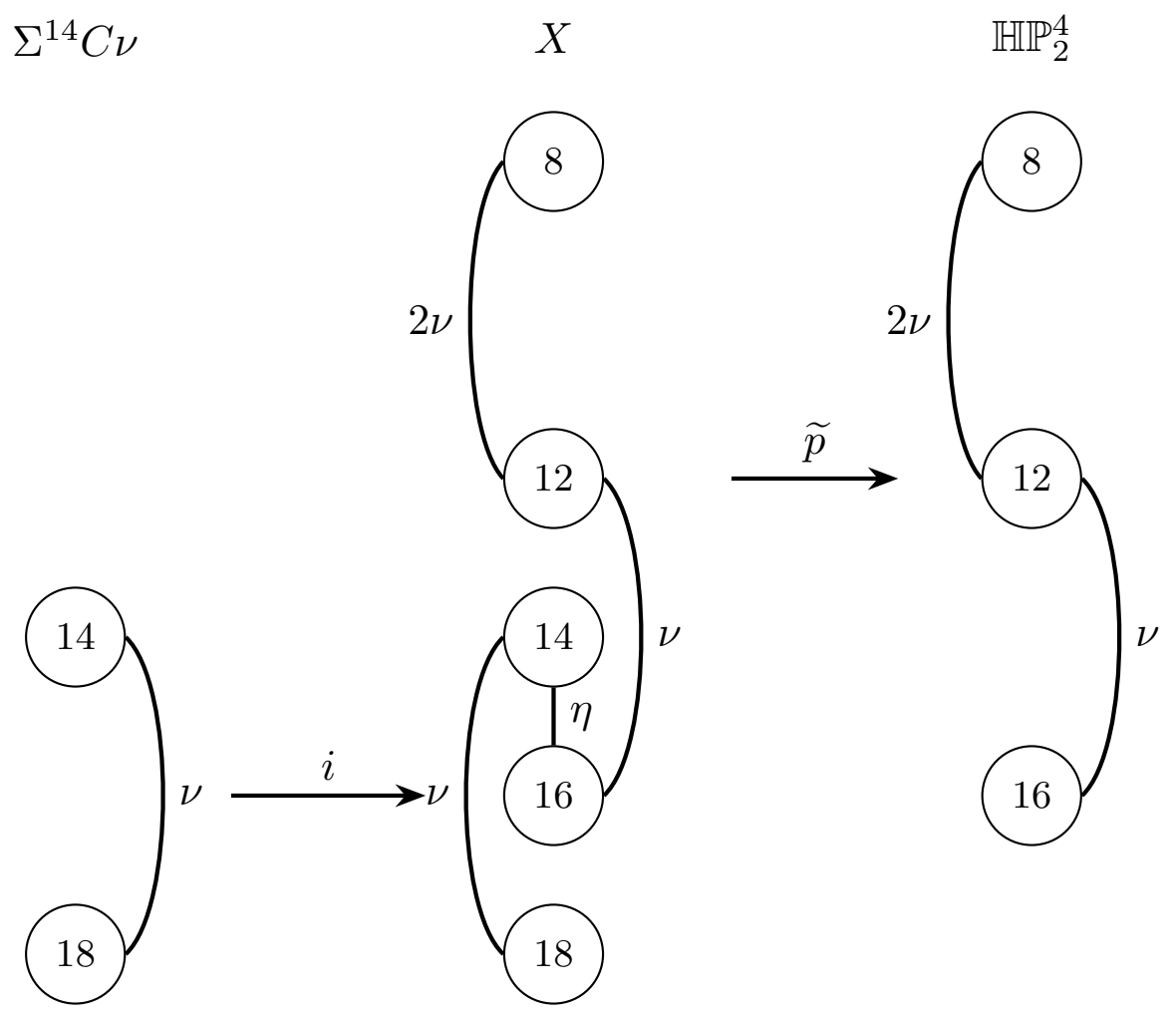

Next we will show that there exist the $d_{10}$-differentials
\[
d_{10}(\epsilon[8])=d_{10}(B[8])=\kappa\nu[18]\ne0.
\]
in the $\tmf$-AHSS of $X$. Naturality then gives the desired differential in $\CP^9_4$. We prove the $d_{10}$ differentials exist in $X$ as follows
\begin{itemize}
    \item We show that the maps 
    \[\tmf^0(\mathbb{HP}^4_2) \to \pi_8\tmf \qquad \text{and} \qquad \tmf(X) \to \pi_8\tmf\]
    induced by inclusion of bottom cells have the same image. The nontriviality of $d_8$-differentials
    \[d_8(\epsilon[8])=d_8(B[8])=\kappa\eta[16]\]
    implies that $\epsilon$ and $B$ are not in the image of $\tmf^0(\mathbb{HP}^4_2) \to \pi_8\tmf$, hence they do not belong to the image of $\tmf(X) \to \pi_8\tmf$. Equivalently, $\epsilon[8]$ and $B[8]$ must support nonzero differentials in the $\tmf$-AHSS of $X$.
    \item The composition
    \[\Sigma^{14}C(\nu) \cong \operatorname{fib}(\widetilde p) \to X \xrightarrow{q} S^{18}\]
    forces the class $\kappa\nu[18]$ to be hit by differential in the $\tmf$-AHSS of $X$.
    \item Finally we rule out the possibility that $\kappa\nu[18]$ is hit by shorter differentials, which concludes our proof.
\end{itemize}

\begin{prop}\label{prop:hp4-d8}
On $\mathbb H\bbP_2^4$, one has
\[
d_8(\epsilon[8])=d_8(B[8])=\kappa\eta[16]\ne0,
\qquad \tmf^1(\mathbb H\bbP_2^4)=0.
\]
The class $B+\epsilon\in[S^8,\tmf]$ extends to $[\mathbb H\bbP_2^4,\tmf]$.
\end{prop}
\begin{proof}
In $\tmf_\ast$ there are the following Toda brackets
\begin{equation*}
\langle \epsilon,2\nu,\nu \rangle=\{\kappa\eta\},\quad
\langle B,2\nu,\nu \rangle=\{\kappa\eta\},\quad
\langle B+\epsilon,2\nu,\nu \rangle=\{0\}.
\end{equation*}
All three brackets have zero indeterminacy since $\pi_7\tmf=0$ and $(\pi_{12}\tmf)\nu=0$. Bauer's bracket $\kappa\eta\in\br{\epsilon,\nu,2\nu}$ \cite[Appendix~A]{BauerTmf2008}, together with scalar shuffling, gives the first equality. For the second, choose $\nu_1\in\br{\bk,\nu,2\nu}$ as in \cite{BauerTmf2008}. The product calculation \cite[Table~8.1]{BGRlocal2021} gives
\begin{equation}\label{eq:cp9-BN}
B\nu_1=\eta\kappa\bk\ne0.
\end{equation}
Toda juggling yields
\begin{equation*}
\bigl(\br{B,2\nu,\nu}\,\bk\bigr)
\cap\bigl(B\,\br{2\nu,\nu,\bk}\bigr)\ne\varnothing.
\end{equation*}
The right-hand side is $\{B\nu_1\}$. Since $\pi_{15}\tmf=\mathbb Z/2\{\kappa\eta\}$, \eqref{eq:cp9-BN} identifies the left bracket. Additivity gives the bracket for $B+\epsilon$.

These classes reach $E_8$, and $\kappa\eta[16]$ is the only nonzero term of total degree $1$. The differentials follow from \cite{wangxu2017triviality}. Furthermore, $\tmf^1(\mathbb H\bbP_2^4)=0$. The zero bracket for $B+\epsilon$ supplies its lift.
\end{proof}

\begin{lemma}\label{lem:cp9-bottomimage}
The composition $S^8 \to X \xrightarrow{\tilde p} \mathbb H\bbP_2^4$ induces a natural inclusion
\begin{equation*}
\operatorname{im}\{\tmf^0(\mathbb{HP}^4_2)\to\pi_8\tmf\} \hookrightarrow \operatorname{im}\{\tmf^0(X)\to\pi_8\tmf,\}
\end{equation*}
which is also an isomorphism.
\end{lemma}
\begin{proof}
Put $K=\Sigma^{14}C\nu$, and let $i:K\to X$ be the inclusion. The cofiber sequence defining $K$ gives
\[
\pi_{18}\tmf \longrightarrow [K,\tmf]
\longrightarrow \pi_{14}\tmf
\xrightarrow{\nu}\pi_{17}\tmf.
\]
The last map is injective, since $\pi_{14}\tmf=\mathbb Z/2\{\kappa\}$ and $\kappa\nu\ne0$. Thus, for every $f:X\to\tmf$, there is an $a\in\pi_{18}\tmf$ such that
\[
f \circ i=a \circ q_K,
\]
where $q_K:K\to S^{18}$ is the top-cell pinch. Choose $q:X\to S^{18}$ so that $q\circ i=q_K$. Then $f-a \circ q$ restricts trivially to $K$ and hence factors through $\widetilde p:X\to\mathbb{HP}_2^4$. Since precomposing $q$ with the bottom-cell inclusion gives a nullhomotopic map, $f-a \circ q$ and $f$ have the same restriction to the bottom cell. This proves the equality of the two images.
\end{proof}

\begin{lemma}\label{lem:cp9-topkernel}
The collapsing map $q:X \to S^{18}$ induces a map with kernel
\begin{equation*}
\ker\bigl(\pi_{17}\tmf\xrightarrow{(q)^*}\tmf^1(X)\bigr)
=\mathbb Z/2\{\kappa\nu\}.
\end{equation*}
In particular $\tmf^1(X)\cong\mathbb Z/2\{\eta B^2[18]\}$.
\end{lemma}
\begin{proof}
The cofiber sequence of Lemma~\ref{lem:cp9-fiber}
\[
\tmf^1(X)\hookrightarrow\tmf^1(\Sigma^{14}C\nu)
=\mathbb Z/2\{\eta B^2[18]\},
\]
The composition
\[\Sigma^{14} C(\nu)\cong\operatorname{fib}(\widetilde p) \to X \xrightarrow{q} S^{18}\]
is the projection to top cell. Thus the injection is an isomorphism, and the composite from $\pi_{17}\tmf$ is the quotient by $\mathbb Z/2\{\kappa\nu\}$.
\end{proof}

The only possible shorter incoming differentials into $\pi_{17}\tmf[18]$ have source in total degree $0$:
\begin{equation*}
\begin{array}{c|c|l}
r&\text{possible source}&\text{reason that it cannot hit }\kappa\nu[18]\\\hline
2&B^2[16]&\Sq^2(x^8)=0\\
4&\kappa[14]&d_2(\kappa[14])=\kappa\eta[16]\ne0\\
6&\{2c_6\}[12]&\text{permanent cycle in the AHSS of }\CP_6^9\\
8&\eta^2B[10]&\text{permanent cycle in the AHSS of }\CP_5^9.
\end{array}
\end{equation*}
The first two rows are immediate. For the last two, it suffices to construct the following finite extensions.

\begin{lemma}\label{lem:cp9-shorter-sources}
The class $\{2c_6\}[12]$ is a nonzero permanent cycle in the $\tmf$--AHSS of $\CP_6^9$.
\end{lemma}
\begin{proof}
$\CP^9_6$ fits into the cofiber sequence
\[
\Sigma^{14}C\nu\longrightarrow\CP_6^9
\longrightarrow\Sigma^{12}C\nu.
\]
A homotopy class detected by $\{2c_6[12]\}$ in $\tmf^0(\Sigma^{12}C \nu)$ pulls back nontrivially to $\tmf^0(\CP^9_6)$.
\end{proof}

\begin{figure}[htbp]\centering
\begin{tikzpicture}[x=.9cm,y=.95cm,scale=.82,transform shape]
\cells{asub}{0}{14,18}\cells{aa}{4.3}{12,14,16,18}\cells{ac}{8.6}{12,16}
\node at (0,-1.65) {$\Sigma^{14}C\nu$};
\node at (4.3,-1.65) {$\CP_6^9$};\node at (8.6,-1.65) {$\Sigma^{12}C\nu$};
\nuedge{asub14}{asub18}{.9}
\etaedge{aa14}{aa16}
\nuedge{aa12}{aa16}{-.9}
\nuedge{aa14}{aa18}{.9}
\nuedge{ac12}{ac16}{.9}
\draw[map] (1.15,-5.4)--node[above,font=\small]{$\text{}$}(3.15,-5.4);
\draw[map] (5.45,-4.7)--node[above,font=\small]{$\text{}$}(7.6,-4.7);
\node[anchor=west,font=\small] at (9.2,-2.7) {$\{2c_6\}$};
\end{tikzpicture}
\caption{The $H\F$--subcomplex on cells $14,18$ and its quotient
on cells $12,16$. }\label{fig:cp9-upper}
\end{figure}
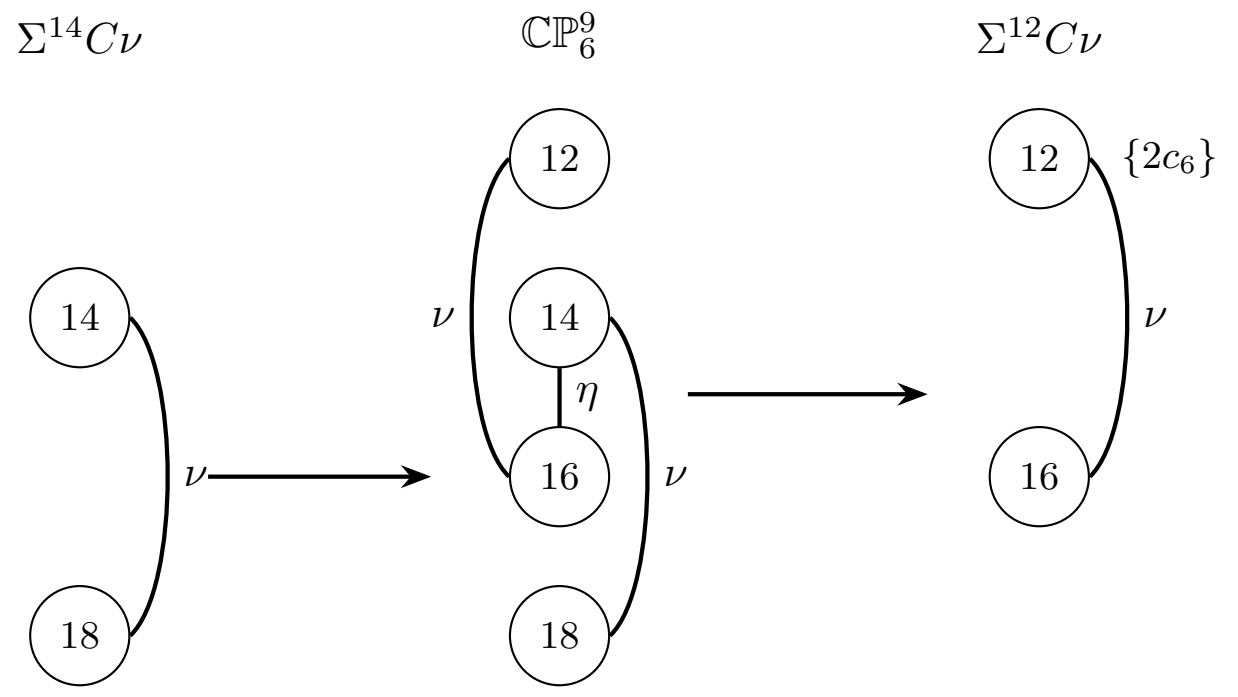

\begin{lemma}\label{lem:cp9-aextension}
The class $\eta^2B[10]$ is a nonzero permanent cycle in the $\tmf$--AHSS of $\CP_5^9$.
\end{lemma}
\begin{proof}

The Toda bracket $\langle \eta^2B,\eta,\nu \rangle = \{0\}$:
\[\langle B\eta^2,\eta,\nu \rangle \cap \langle B,\eta^3,\nu \rangle \neq \varnothing ,\]
where both indeterminacies vanish by $\pi_5\tmf=\pi_7\tmf=0$ and $(\pi_{12}\tmf)\nu=0$. Since $\eta^3=4\nu$, Proposition~\ref{prop:hp4-d8} gives
\[
0\in2\br{B,2\nu,\nu}
\subseteq\br{B,4\nu,\nu}=\br{B,\eta^3,\nu}.
\]
Hence in $\tmf^0(\CP^8_5)$ there is a nonzero homotopy class detected by $\eta^2B[10]$ which we denote by $\{\eta^2B[10]\}$.

Using the splitting
\[
A:=\CP_5^7\simeq\Sigma^{10}C\eta\vee S^{14},
\]
choose an extension $f:\Sigma^{10}C\eta\to\tmf$ of $\eta^2B$, and take the map $(f,0)$ on $A$. The vanishing of $\langle\eta^2B,\eta,\nu\rangle$ allows this map to extend to $f:\CP_5^8\to\tmf$.

Let $\alpha:S^{17}\to\CP_5^8$ be the map attaching the $18$-cell. Its projection to the $16$-cell is zero, since $\Sq^2(x^8)=0$, so $\alpha$ lifts to $A$. The $\Sigma^{10}C\eta$-component of this lift factors through its bottom cell, since $\pi_5S^0=0$. Consequently, after smashing with $\tmf$, the attaching map factors through the $14$-cell:
\[
\alpha\wedge\tmf:
\Sigma^{17}\tmf
\xrightarrow{\nu}
\Sigma^{14}\tmf
\xrightarrow{i\wedge\tmf}
\CP_5^8\wedge\tmf,
\]
since $\pi_7\tmf=0$. As $f\circ i=0$, the $\tmf$-linear extension of $f$ annihilates this composite. By adjunction, $f\alpha=0$. Thus $f$ extends over $\CP_5^9$, with the same nonzero bottom-cell restriction $\eta^2B$.
\end{proof}

\begin{figure}[htbp]\centering
\begin{tikzpicture}[x=.9cm,y=.95cm,scale=.82,transform shape]
\cells{qq}{0}{10,12,14,16,18}
\node at (0,-.3) {$\CP_5^9$};
\etaedge{qq10}{qq12}\etaedge{qq14}{qq16}
\nuedge{qq12}{qq16}{.9}\nuedge{qq14}{qq18}{-.9}
% \node[anchor=west] at (1.7,-.66) {$\eta^2B$};
% \node[anchor=west] at (1.7,-1.5) {$\br{\eta^2B,\eta,\nu}=0$};
% \node[anchor=west] at (1.7,-2.1) {$0$ on cell $14$};
% \node[anchor=west] at (1.7,-3.3) {$\eta^2B\,h(s\sigma)+0\cdot\nu=0$};
\draw[dashed,gray,line width=.8pt] (qq10.west) .. controls (-2,-2.15) and (-2,-5.95) .. node[midway,left,align=center]{possible\\$s\sigma$}(qq18.west);
\end{tikzpicture}
\caption{The class $\eta^2B[10]$ is a permanent cycle in the
$\tmf$--AHSS of $\CP_5^9$, ruling out the possible $d_8$. The dashed attaching map is allowed; its image in $\tmf$ vanishes.}
\end{figure}
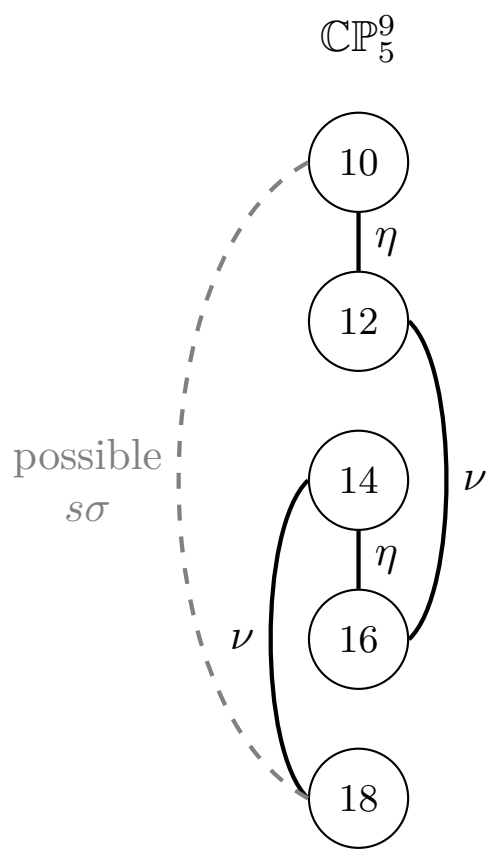

\begin{proof}[Proof of the $d_{10}$-differential in Theorem~\ref{thm: cp4}]
By Lemmas~\ref{lem:cp9-topkernel}, ~\ref{lem:cp9-shorter-sources} and ~\ref{lem:cp9-aextension}, $\kappa\nu[18]$ must be killed by $d_{10}$-differentials in the $\tmf$-AHSS of $X$, whereas by Lemma~\ref{lem:cp9-bottomimage}, each of $\epsilon$ and $B$ must support a nonzero differential. The target of the potential $d_8$ differential $\kappa\eta[16]$ is already killed by
\[d_2(\kappa[14]) = \kappa\eta[16].\]
Consequently
\[
d_{10}(\epsilon[8])=d_{10}(B[8])=\kappa\nu[18]\ne0.
\]
Naturality of $H\mathbb{F}_2$ quotient $\CP^9_4 \to X$ completes the proof.
\end{proof}

\par
\endgroup

Next we would like to choose suitable lifts of $\Delta^{8k}\epsilon[8+192k]$ and $\Delta^{8k+4}\epsilon[104+192k]$ such that they extend to $[\CP^{6+96k},S^0]$ and $[\CP^{54+96k},S^0]$, respectively. Put $m=4+96k$ or $m=52+96k$. In either case, $\CP^{m+1}_m \cong S^{2m} \vee S^{2m+2}$ whereas $\Sigma^{2m}C(2\nu)$ is a quotient of $\CP^{m+2}_m$. Therefore any lift which is a $2\nu$-torsion will extend to $[\CP^{m+2},S^0]$. 

For $2\nu$-torsions in dimension $8+192k$, we periodize $\epsilon+\eta\sigma$ using the finite type $2$ complex $M(3,8)$ \cite{BMQhurewicz2021}. For $2\nu$-torsions in dimension $104+192k$, we use an element in the Modified Adams spectral sequence (MASS) of $M(3,8)$ with $(s,t) = (26, 148)$, which detects a $2$-torsion element in $\pi_{122}M(3,8)$. We postpone their detailed construction to Theorem~\ref{thm:epsilon lift}.

\begin{theorem}\label{thm:CP 52}
The spaces $\CP^{52+96k}$ admit exotic smooth structures for every $k\geq0$.    
\end{theorem}

\begin{proof}
Similar to the proof of Theorem~\ref{thm: cp4}, we can show there exist $d_{10}$-differentials
\[d_{10}(\Delta^{8k+4}\epsilon[104+192k]) = \Delta^{8k}\kappa_4\nu[114+192k]\]
Theorem~\ref{AHSSsmoothing} gives the desired exotic smooth structures. 
\end{proof}

\begin{theorem}\label{thm:CP 6 and 54}
The spaces $\CP^{6+96k}$ and $\CP^{54+96k}$ admit exotic smooth structures for every $k\geq0$.
\end{theorem}

\begin{proof}
    By Theorems~\ref{thm:tmf-detection},~\ref{thm: cp4}, ~\ref{thm:CP 52} and ~\ref{thm:epsilon lift}, $\epsilon_k[8+192k]$ and $\mathcal{E}_k[104+192k]$ are nonzero in $[\CP^{6+96k},S^0]$ and $[\CP^{54+96k},S^0]$, respectively. Theorem~\ref{AHSSsmoothing} gives the desired exotic smooth structures. 
\end{proof}

\begin{proof}[Proof of Theorem~\ref{thm:tmf-primarymain}]
The cases $m\in\{3,7,27,51,55,75\}$ follow from Theorem~\ref{thm:tmf-primary}, the case $m=4$ follows from Theorem~\ref{thm: cp4}, the cases $m=6,54$ from Theorem~\ref{thm:CP 6 and 54}, the cases $m=28,76$ follow from Theorem~\ref{thm:tmf-doubled}, the case $m=52$ follows from Theorem~\ref{thm:CP 52}, and the case $m=56$ follows from Theorem~\ref{thm:tmf-eta-extensions}.
\end{proof}

\section{Computations of \texorpdfstring{$\tmf$}{tmf}-AHSS hidden extensions}
\label{sec: tmf-extension method}

We work $2$-locally in this section. Using the $\tmf$-extension methods, we obtain exotic smooth structures in the following dimensions.

\begin{theorem}\label{thm:delta method}
The spaces $\CP^{m+96k}$ admit exotic smooth structures for
$$
m\in\{10,11,20,22,23,24,35,46,47,58,59,66,67,70,71,72,82,83\},\qquad k\geq0,
$$

\end{theorem}

We begin with a lemma relating Toda brackets to hidden extensions in the AHSS.

\begin{lemma}\label{hiddenextension}
Let $E$ be a commutative ring spectrum with multiplication $\mu$. Let $\alpha,\beta,\gamma,\delta$ be elements in $\pi_\ast(E)$. Suppose $\alpha\beta = \beta\gamma = 0$. The formulas below describe hidden extensions between $E_\infty$-classes, with compatible choices of representatives, attaching maps, and Toda defining systems, and with the usual indeterminacies.
\begin{enumerate}
    \item In the category of $E$-modules, denote the cofiber of
    \[\Sigma^{|\alpha|} E \xrightarrow{\alpha \wedge id} E \wedge E \xrightarrow{\mu} E\]
    by $C(\alpha)$. Then in the $E$-AHSS of $C(\alpha)$ there is a hidden $\gamma$-extension
    \[\gamma \cdot \beta[0] = \langle \gamma,\beta,\alpha\rangle[|\alpha|+1].\]
    \item Suppose $\langle\delta,\gamma,\beta,\alpha\rangle$ is well defined. In the category of $E$-modules, denote the cofiber of
    \[\Sigma^{|\beta|+|\alpha|+1}E \xrightarrow{\alpha[|\beta|+1]} C(\beta)\]
    by $X$, where $\alpha[|\beta|+1]$ is a lift of
    \[\Sigma^{|\beta|+|\alpha|+1}E \xrightarrow{(\Sigma^{|\beta|+1}\alpha) \wedge id} \Sigma^{|\beta|+1}E \wedge E \xrightarrow{\mu} \Sigma^{|\beta|+1}E\]
    along the projection $C(\beta) \to \Sigma^{|\beta|+1}E$. Then in the $E$-AHSS of $X$ there is a hidden $\delta$-extension
    \[\delta \cdot \gamma[0] = \langle \delta,\gamma,\beta,\alpha\rangle[|\alpha|+|\beta|+2].\]
\end{enumerate}
\end{lemma}

\begin{proof}
    The proof is analogous to that of Lemma~5.1 in \cite{wangxu201851stem}.
\end{proof}

\begin{lemma}\label{lem:tmf-periodic-thom}
There exists a $\tmf$-module equivalence
\begin{equation}\label{eq:tmf-periodic-thom}
 \tmf\wedge\CP_{a+96k}^{b+96k}
 \simeq\Sigma^{192k}\tmf\wedge\CP_a^b.
\end{equation}    
\end{lemma}
\begin{proof}
Denote by $L$ the tautological complex line bundle on $\CP^\infty$ and $c_1$ its first Chern class. The real virtual bundle underlying $9L-L^{\otimes3}$ has $w_2=0$ and $p_1=9(c_1)^2-(3c_1)^2=0$. Since $H^4(\CP^\infty;\mathbb Z)$ is torsion-free, its spin characteristic class $p_1/2$ vanishes, so it admits a string structure. On each finite projective skeleton, the Adams conjecture identifies its $2$-local stable spherical fibration with that of $8L$. Transporting the string orientation therefore gives a $\tmf$-orientation of $8L$; see \cite{davismahowald2008nonimmersions}. The Thom description $\CP_a^b\simeq\operatorname{Th}(\CP^{b-a};aL)$ then gives \eqref{eq:tmf-periodic-thom} since $96k$ is divisible by $8$. The Thom class may be normalized on each fiber, so the comparison respects skeletal restrictions and is the identity on cellular coefficients.
\end{proof}

\subsection{Unbranched \texorpdfstring{$\Delta$}{Delta}-extensions}
\label{subsec:tmf-unbranched-delta}

\begin{lemma}\label{lem:tmf-five-two-cell-brackets}
The following Toda cosets hold in $\pi_*\tmf$:
\[
\begin{array}{c|c|c}
\beta&d&\langle2,\beta,\eta\rangle\subseteq\pi_{d+2}\tmf\\\hline
\eta^2\bar{\kappa}&22&8\Delta+2\pi_{24}\tmf\\
\eta\eta_1\bar{\kappa}&46&4\Delta^2+2\pi_{48}\tmf\\
\eta_1^2\bar{\kappa}&70&8\Delta^3+2\pi_{72}\tmf\\
\eta\eta_4\bar{\kappa}&118&8\Delta^5+2\pi_{120}\tmf\\
\epsilon_5\kappa&142&4\Delta^6+2\pi_{144}\tmf.
\end{array}
\]
where $\eta_1=\{\Delta\eta\} \in \pi_{25}\tmf$ and $\eta_4=\{\Delta^4\eta\} \in \pi_{97}\tmf$. Moreover,
\[
\begin{gathered}
8\Delta\in\langle4,2\bar{\kappa},\nu\rangle,\qquad
8\Delta^5\in\langle4,(2\Delta^4)\bar{\kappa},\nu\rangle,\\
8\Delta^7\in\langle2,\nu\nu_6\kappa,\nu\rangle.
\end{gathered}
\]
\end{lemma}
\begin{proof}
For the five rows, use the $d_7$-differentials on $4\Delta,2\Delta^2,4\Delta^3,4\Delta^5,2\Delta^6$, respectively, in \cite[Section~8]{BauerTmf2008}. Their targets are the $h_1$-multiples of the displayed coefficients, up to odd units. The resulting Massey products have the displayed values modulo twice the integral coefficient group. The remaining three brackets follow from
\[
d_5(\Delta)=gh_2,\qquad
d_5(2\Delta^5)=2\Delta^4gh_2,\qquad
d_7(4\Delta^7)=4g\Delta^6h_2,
\]
where $4g\Delta^6$ detects $\nu\nu_6\kappa$. Generalized Moss convergence \cite{moss1970secondary,ullman2013regular,BelmontKong2024Toda} applies with no crossing differentials. For the five rows, the remaining indeterminacy is $2\pi_{d+2}\tmf$: the other term is $(\pi_{d+1}\tmf)\eta=0$.
\end{proof}

\begin{lemma}\label{lem:90 to 96 extension}
Put $w=\eta_1\bar{\kappa}\in\pi_{45}\tmf$. The following fourfold Toda bracket is well-defined and satisfies
\begin{equation*}
\langle2,w^2,\eta,\nu\rangle
=2\Delta^4+2\pi_{96}\tmf.
\end{equation*}
\end{lemma}
\begin{proof}
The bracket is well-defined since $\pi_4\tmf = \pi_{91}\tmf = \pi_{92}\tmf =0$. The indeterminacy is $2\pi_{96}\tmf$. We shuffle by
\[\langle 2,w^2,\eta,\nu\rangle \bar{\kappa}^2 \subseteq\langle 2,w^2, \langle\eta,\nu,\bar{\kappa}^2\rangle\rangle = \langle2,w^2,w\rangle.\]
There is a $d_7$-differential \cite{BauerTmf2008}
\[d_7(\Delta^4\bar{\kappa}^2) = w^3.\]
Applying the Moss convergence theorem gives $2\Delta^4\bar{\kappa}^2 \in \langle2,w^2,w\rangle$, which proves the listed bracket.
\end{proof}

\begin{prop}\label{prop:tmf-integral-local}
\label{prop:tmf-w-square-finite}
For each row of the following table, there is a hidden extension $e\cdot\alpha[2m]=\delta[2N]$ in the $\tmf$-AHSS of $\tmf^0(\CP_m^N)$ where $e$ is $2$ or $4$.
\[
\begin{array}{c|c|c|c|c}
m&\alpha&e&N&\delta\\\hline
10&2\bar{\kappa}&4&12&8\Delta\\
11&\bar{\kappa}\eta^2&2&12&8\Delta\\
35&\eta_1^2\bar{\kappa}&2&36&8\Delta^3\\
45&w^2&2&48&2\Delta^4\\
58&(2\Delta^4)\bar{\kappa}&4&60&8\Delta^5\\
59&\eta\eta_4\bar{\kappa}&2&60&8\Delta^5\\
82&\nu\nu_6\kappa&2&84&8\Delta^7
\end{array}
\]
The same statements hold after adding $96k$ to $m,N$ and multiplying both coefficients by $\Delta^{8k}$.
\end{prop}
\begin{proof}
For $m=11,35,59$, the stunted projective spectrum can be identified with 
\[
\CP_m^{m+1}\simeq \Sigma^{2m}C\eta,
\]
applying Lemma~\ref{hiddenextension}(1) proves the hidden extensions. For $m=10,58,82$, the stunted projective can be identified with
\[
\CP_m^{m+2}\simeq \text{Cof}(S^{2m+3} \xrightarrow{(\nu,\eta)} S^{2m} \vee S^{2m+2}).
\]
Collapse the middle sphere and apply Lemma~\ref{hiddenextension}(1) to $\Sigma^{2m}C(\nu)$, using Lemma~\ref{lem:tmf-five-two-cell-brackets}. Pullback then gives the asserted extensions. For $m=45$, collapse the separate $94$-cell in $\CP_{45}^{47}\simeq\Sigma^{90}C\eta\vee S^{94}$. The quotient is the cofiber of
\[\Sigma^{90}\{ S^5 \xrightarrow{\nu[2]} C\eta \}.\]
Applying Lemma~\ref{hiddenextension}(2) with $(\delta,\gamma,\beta,\alpha)=(2,w^2,\eta,u\nu)$ proves the hidden 2-extension. The equivalences of Lemma~\ref{lem:tmf-periodic-thom} give the statement multiplied by $\Delta^{8k}$.
\end{proof}

\begin{figure}[htbp]\centering
\begin{tikzpicture}[x=.9cm,y=.95cm,scale=.82,transform shape]
\tikzset{cell/.style={circle,draw,fill=white,minimum size=8mm,inner sep=.4pt,font=\small},
 eta/.style={black,line width=.8pt},
 nu/.style={black,line width=.8pt},
 twonu/.style={black,line width=.8pt},
 nuedge/.style={black,line width=.8pt},
 twonuedge/.style={black,line width=.8pt}}
\foreach\d/\y in{90/0,92/-1.35,94/-2.7,96/-4.05}
 {\node[cell](a\d)at(0,\y){$\d$};}
\draw[eta](a90)--node[right]{$\eta$}(a92);
\draw[eta](a94)--node[right]{$\eta$}(a96);
\draw[nu](a92.west)..controls(-.9,-1.75)and(-.9,-3.65)..node[midway,left]{$\nu$}(a96.west);
\foreach\d/\y in{90/0,92/-1.35,96/-4.05}
 {\node[cell](b\d)at(4.3,\y){$\d$};}
\draw[eta](b90)--node[right]{$\eta$}(b92);
\draw[nu](b92.west)..controls(3.4,-1.75)and(3.4,-3.65)..node[midway,left]{$\nu$}(b96.west);
\node[cell](c96)at(8.6,-4.05){$96$};
\draw[->,black,line width=.8pt](1,-1.9)--node[above]{}(2.9,-1.9);
\draw[->,black,line width=.8pt](5.25,-4.05)--node[above]{$q$}(7.85,-4.05);
\node at(0,1.05){$\CP_{45}^{48}$};
\node[align=center]at(4.3,1.05){$\Sigma^{90}(C\eta\cup_{u\widetilde\nu}e^6)$};
\node at(8.6,-3){$S^{96}$};
\end{tikzpicture}
\caption{The actual three-cell quotient used in Proposition~\ref{prop:tmf-w-square-finite}.}
\end{figure}
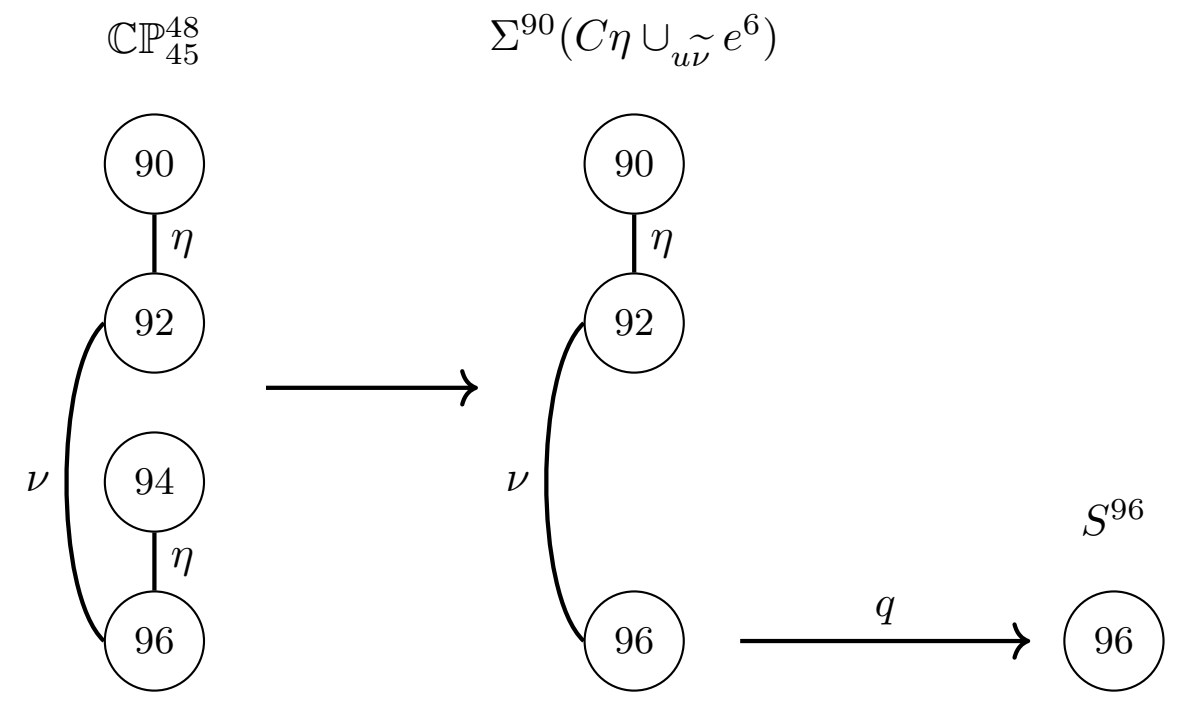
\begin{theorem}\label{thm:tmf-cp-10-45-58-82}
For every $k\geq0$, $\CP^{m+96k}$ admits an exotic smooth structure for
\[
m\in\{10,11,35,45,58,59,82,83\}.
\]
\end{theorem}
\begin{proof}
For $m=11,35,59,82$, applying Lemma~\ref{lem:tmf-integral-no-jump-reformulation} with $s=1$ and combining Lemma~\ref{lem:tmf-five-two-cell-brackets} shows that $\Delta^{8k}\bar{\kappa}\eta^2[22+192k],\Delta^{8k}\eta_1^2\bar \kappa[70+192k], \eta\eta_4\bar{\kappa}[118+192k]$ and $\Delta^{8k}\nu\nu_6\kappa[164+192k]$ are nonzero in $\tmf^0(\CP^{m+96k})$ respectively. For $m=83$ ,$\CP^{83+96k}_{82+96k} \cong S^{164+192k} \vee S^{166+192k}$, hence the lift of $\Delta^{8k}\nu\nu_6\kappa[164+192k]$ in $[\CP^{82+96k},S^0]$ can extend to $[\CP^{83+96k},S^0]$.
 
For $m=10,58$, let $s=2$ in Lemma~\ref{lem:tmf-integral-no-jump-reformulation} and combining Lemma~\ref{lem:tmf-five-two-cell-brackets} provides that $\Delta^{8k}2\bar{\kappa}[20+192k]$ and $2\Delta^{8k+4}\bar{\kappa}[116+192k]$ are nonzero in $\tmf^0(\CP^{m+96k})$ respectively.

For $m=45$, let $a$ be a homotopy class in $\tmf^0(\CP^{48+96k}_{45+96k})$ detected by $\Delta^{8k}w^2[90+192k]$. By Proposition~\ref{prop:tmf-integral-local}, $2a$ is detected by $\Delta^{8k}2\Delta^{4}[96+192k]$. The natural map
\[q^\ast \colon \tmf\text{-AHSS of }\CP^{48+96k}_{45+96k} \to \tmf\text{-AHSS of }(\CP^{48+96k})\]
sends $\Delta^{8k}2\Delta^{4}[96+192k]$ to the element of the same name, hence $0 \neq q^\ast(2a) \in \tmf^0(\CP^{48+96k})$. Consequently $q^\ast(a) \neq 0$, which can only be detected by $\Delta^{8k}w^2[90+192k]$: $\pi_{i}\tmf$ for $i \equiv 91,92,93,94,95,96 (\text{mod } 192)$ are torsion-free, whereas the free components in $\pi_{92+192k}\tmf$ cannot extend to $2a$. 

Theorem~\ref{AHSSsmoothing} thus gives the claimed exotic $\CP^{m+96k}$. 
\end{proof}

\begin{theorem}\label{thm:tmf-cp-46-47}
For every $k\geq0$, $\CP^{46+96k}$ and $\CP^{47+96k}$ admit exotic smooth structures.
\end{theorem}
\begin{proof}
We use the specified lifts $\rho_k\in\pi_{90+192k}S^0$ with
\[
u_*(\rho_k)=\Delta^{8k}w^2,\qquad \eta\rho_k=0,
\]
provided by \cite{BQeta2026}. By Theorem~\ref{thm:tmf-detection}, $\rho_k[90+192k]$ detects a nonzero homotopy class in $[\CP^{45+96k},S^0]$. $\CP^{47+96k}_{45+96k} \cong \Sigma^{90+192k}C\eta \vee S^{94+192k}$, hence this homotopy class extends to $[\CP^{47+96k},S^0]$.
\end{proof}

\subsection{Branched extensions and torsion classes}
\label{subsec:tmf-branched-delta}
We prove Theorem~\ref{thm:cp20-24-finite-survival} in two stages. First we construct the local classes and compute their hidden extensions on stunted projective spectra. We then show that their pullbacks to $\CP^{24+96n}$ preserve the required cellular filtrations.

\begin{prop}\label{prop:cp20-24-local2}
For every $k\ge0$ there is a hidden $2$-extension
\[
2\bigl(\Delta^{8k}\eta\eta_1\bar{\kappa}[46+192k]\bigr)
=((4\Delta^2)\Delta^{8k})[48+192k]
\]
in the $\tmf$--AHSS of $\CP^{24+96k}_{23+96k}$.
\end{prop}
\begin{proof}
First take $k=0$. We compute $\langle 2,\eta\eta_1\bar{\kappa},\eta\rangle$. In the descent spectral sequence, $\eta\eta_1\bar{\kappa}$ is detected by $\Delta gh_1^2$, and $d_7(2\Delta^2)=\Delta gh_1^3$. Thus Moss convergence gives
\[
\langle 2,\eta\eta_1\bar{\kappa},\eta\rangle=4\Delta^2+2\pi_{48}\tmf.
\]
$\CP_{23}^{24}\simeq\Sigma^{46}C\eta$, hence Lemma~\ref{lem:tmf-five-two-cell-brackets} provides the hidden extension. For general $k$ the result follows by multiplying $\Delta^{8k}$ on both the source and the target of the extension.
\end{proof}

Next we show that in the $\tmf$-AHSS of $\CP^{24+96k}_{20+96k}$, there is a hidden $4$-extension
\[4 \cdot \Delta^{8k}\bar{\kappa}^2[40+192k] = \Delta^{8k} 4\Delta^2[48+192k].\]
Similar to Proposition~\ref{prop:CP to HP} there exists an $H\mathbb F_2$-quotient map
\[
\widetilde p:\CP_{20}^{24}\longrightarrow \mathbb{HP}_{10}^{12}.
\]
We first establish the extension in the $\tmf$-AHSS of $\mathbb{HP}_{10}^{12}$, then pullback via $\widetilde p$ to proof the extension in $\CP^{24}_{20}$.

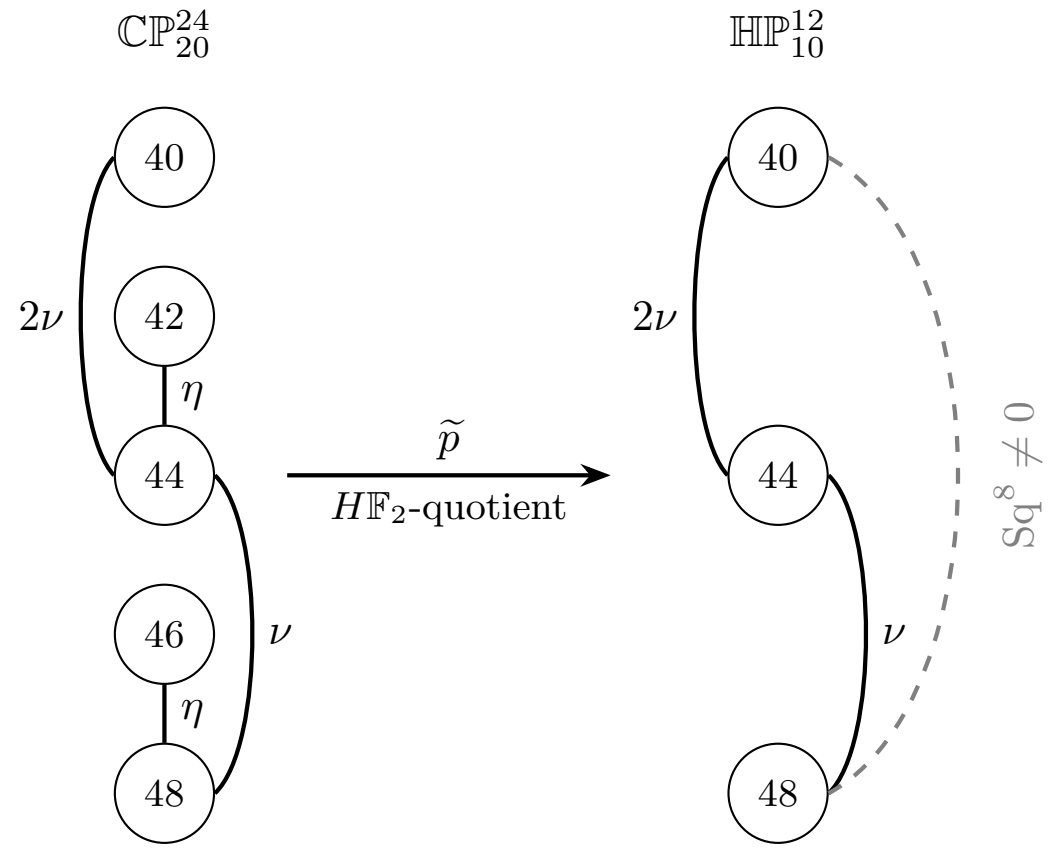
\begin{figure}[htbp]\centering
\begin{tikzpicture}[x=.9cm,y=.95cm,scale=.82,transform shape]
\colorlet{nublue}{black}
\colorlet{twonu}{black}
\colorlet{higher}{gray}
\tikzset{cell/.style={circle,draw,fill=white,minimum size=8mm,inner sep=.4pt,font=\small},
 etaedge/.style={black,line width=.8pt},nuedge/.style={nublue,line width=.8pt},
 twoedge/.style={twonu,line width=.8pt},map/.style={black,-{Stealth[length=2mm]},line width=.8pt}}
\foreach \d/\y in {40/0,42/-1.35,44/-2.7,46/-4.05,48/-5.4}
 {\node[cell] (a\d) at (0,\y) {$\d$};}
\draw[etaedge] (a42)--node[right]{$\eta$}(a44);\draw[etaedge] (a46)--node[right]{$\eta$}(a48);
\draw[twoedge] (a40.west) .. controls (-.85,-.4) and (-.85,-2.3) .. node[midway,left]{$2\nu$}(a44.west);
\draw[nuedge] (a44.east) .. controls (.9,-3.1) and (.9,-5) .. node[midway,right]{$\nu$}(a48.east);
\foreach \d/\y in {40/0,44/-2.7,48/-5.4}
 {\node[cell] (b\d) at (5.5,\y) {$\d$};}
\draw[twoedge] (b40.west) .. controls (4.65,-.4) and (4.65,-2.3) .. node[midway,left]{$2\nu$}(b44.west);
\draw[nuedge] (b44.east) .. controls (6.4,-3.1) and (6.4,-5) .. node[midway,right]{$\nu$}(b48.east);
\draw[higher,dashed,line width=.8pt] (b40.east)..controls(7.5,-.8)and(7.5,-4.6)..(b48.east);
\node[higher,rotate=90] at (7.65,-2.7) {$\mathrm{Sq}^8\ne0$};
\node at (0,1.05) {$\CP_{20}^{24}$};\node at (5.5,1.05) {$\mathbb H\bbP_{10}^{12}$};
\draw[map] (1.1,-2.7)--node[above]{$\widetilde p$}node[below,font=\small]{$H\mathbb F_2$-quotient}(4,-2.7);
\end{tikzpicture}
\caption{The $H\mathbb F_2$ quotient $\widetilde p$ is used. The dashed arc records $\mathrm{Sq}^8$ on the stunted quaternionic projective spectrum, whose effect is negligible in $\tmf$-AHSS computations.}
\end{figure}

\begin{prop}
For every $k\ge0 $ there is a hidden $4$-extension
    \[4 \cdot \Delta^{8k}\bar{\kappa}^2[40+192k] = \Delta^{8k} 4\Delta^2[48+192k]\]
in the $\tmf$--AHSS of $\mathbb{HP}^{12+48k}_{10+48k}$.
\end{prop}

\begin{proof}
Both $\Delta^{8k}\bar{\kappa}^2[40+192k]$ and $\Delta^{8k} 4\Delta^2[48+192k]$ are surviving cycles in the AHSS. We first compute the fourfold Toda bracket
\[
\langle 4,\bar{\kappa}^2,2\nu,\nu\rangle\subseteq\pi_{48}\tmf.
\]

In the descent spectral sequence, $d_5(\Delta dh_2)=4g^2$, $d_5(2\Delta g) = 2g^2h_2$ and $d_5(\Delta^2) = 2\Delta gh_2$. Hence on the DSS $(E_5,d_5)$ the Massey product contains the class
\[2\Delta^2 \in \langle 2,g^2,2h_2,h_2\rangle\]
There are no crossing differentials, so generalized Moss convergence \cite{moss1970secondary,ullman2013regular,BelmontKong2024Toda} supplies
\[
4\Delta^2 \in \langle4,\bar{\kappa}^2,2\nu,\nu \rangle.
\]
In particular, multiplying by $\Delta^{8k}$ on both sides preserves the relation. Lemma~\ref{hiddenextension}(2), with $\alpha=\nu$, $\beta=2\nu$, $\gamma=\bar{\kappa}^2$, $\delta=4$, in turn gives the hidden $4$-extension in the cohomological $\tmf$-AHSS of $\mathbb{HP}_{10}^{12}$.
\end{proof}

\begin{prop}\label{prop:4 extension in CP^24_20}
For every $k\ge0 $ there is a hidden $4$-extension
    \[4 \cdot \Delta^{8k}\bar{\kappa}^2[40+192k] = \Delta^{8k} 4\Delta^2[48+192k]\]
in the $\tmf$--AHSS of $\mathbb{CP}^{24+96k}_{20+96k}$.    
\end{prop}
\begin{proof}
Both $\Delta^{8k}\bar{\kappa}^2[40+192k]$ and $\Delta^{8k} 4\Delta^2[48+192k]$ are surviving cycles in the AHSS of $\CP^{24+48k}_{20+48k}$. The extension then follows by naturality: On the level of $\tmf$-cohomology groups, $\tmf^0(\mathbb{HP}^{12+48k}_{10+48k})$ has an integral subgroup whose generator is detected by $\Delta^8\bar{\kappa}^2[40+192k]$, whose image under $\widetilde{p}^\ast$ is detected by the element of the same name. $4$-times this generator and $4$-times its image are detected by $4\Delta^2[48+192k]$ in both spectral sequences, which amounts to the existence of the hidden extension.
\end{proof}

\begin{prop}\label{prop:2 extensions in CP^24_19}
For every $k\ge0 $ there are successive hidden $2$-extensions
    \[2 \cdot \Delta^{8k}\bar{\kappa}^2[40+192k] = \Delta^{8k} \kappa^3[42+192k] \qquad \text{and} \qquad 2 \cdot \Delta^{8k}\kappa^3[42+192k] = \Delta^{8k} 4\Delta^2[48+192k]\]
in the $\tmf$--AHSS of $\CP^{24+96k}_{19+96k}$.   
\end{prop}
\begin{proof}
Applying Lemma~\ref{lem:tmf-periodic-thom} and multiplication by $\Delta^{8k}$ will reduce all such extensions to the base case where $k=0$. All the aforementioned classes are surviving cycles in the $\tmf$-AHSS of $\CP^{24}_{19}$.

For the first extension, the bottom-cell cofibration of $\CP_{19}^{21}$ gives the exact sequence
\[
\pi_{39}\tmf\xrightarrow{(\eta,\nu)}
\pi_{40}\tmf\oplus\pi_{42}\tmf
\longrightarrow\tmf^0(\CP_{19}^{21})\longrightarrow0.
\]
The first map has image generated by $(2\bar{\kappa}^2,\kappa^3)$, and $2\pi_{42}\tmf=0$. Thus there is a $2$-extension
\[2 \cdot \bar{\kappa}^2[40] = \kappa^3[42].\]
in the $\tmf$-AHSS of $\CP^{21}_{19}$. Consider the consecutive natural maps
\[\CP^{21}_{19} \xrightarrow{i} \CP^{24}_{19} \xrightarrow{q} \CP^{24}_{20}\]
Naturality arguments similar to those in Proposition~\ref{prop:4 extension in CP^24_20} establish the hidden 2-extensions.

Let $q:\CP_{19}^{24}\to\CP_{20}^{24}$ be the quotient map. By Proposition~\ref{prop:4 extension in CP^24_20}, choose $a\in\tmf^0(\CP_{20}^{24})$ detected by $\bar{\kappa}^2[40]$ such that $4a$ is detected by $4\Delta^2[48]$. Put $A=q^\ast a$. The preceding computation, together with $2\pi_{42}\tmf=0$, shows that $A$ and $2A$ are detected by $\bar{\kappa}^2[40]$ and $\kappa^3[42]$, respectively. By naturality,
\[
2(2A)=4A=q^\ast(4a)
\]
is detected by the surviving top-cell class $4\Delta^2[48]$. Taking $\kappa^3[42]$ to denote the representative $2A$ therefore gives the second hidden $2$-extension.

% \begin{figure}[!b]\centering
% \begin{tikzpicture}
% \definecolor{nublue}{RGB}{26,94,165}
% \definecolor{twonu}{RGB}{190,83,18}
% \definecolor{higher}{RGB}{126,63,153}
% \tikzset{cell/.style={circle,draw,fill=white,minimum size=7.5mm,inner sep=.3pt,font=\small},
%  etaedge/.style={line width=.8pt},nuedge/.style={draw=nublue,line width=.9pt},
%  twoedge/.style={draw=twonu,line width=.9pt},map/.style={-{Stealth[length=2mm]},line width=.8pt}}
% \node[cell](a38)at(0,0){$38$};\node[cell](a40)at(0,-1.35){$40$};\node[cell](a42)at(0,-2.7){$42$};
% \draw[etaedge] (a38)--node[right=2pt,font=\scriptsize]{$\eta$}(a40);\draw[nuedge] (a38) .. controls ($(a38)+(-.75,0)$) and ($(a42)+(-.75,0)$) .. (a42);
% \node[nublue,font=\small]at(-.95,-1.35){$u\nu$};
% \node[cell](b40)at(4,-1.35){$40$};\node[cell](b42)at(4,-2.7){$42$};
% \node[cell](c39)at(8,-1.35){$39$};
% \node at(0,.65){$\bbP_{19}^{21}$};\node at(4,.65){$S^{40}\vee S^{42}$};\node at(8,.65){$S^{39}$};
% \draw[map](1,-1.1)--node[above,font=\small]{collapse}(3,-1.1);
% \draw[map](b40.east)--node[above,font=\small]{$\eta$}(c39.west);
% \draw[map](b42.east)--node[below,sloped,font=\small]{$u\nu$}(c39.south west);
% \node[font=\small]at(4,-3.4){$(2\bar{\kappa}^2,\kappa^3)$};\node[font=\small]at(8,-3.4){$\kappa\eta_1$};
% \end{tikzpicture}
% \caption{The actual cofibration identifying the double at cell $42$.
% The coefficient map out of $S^{39}$ accounts for precisely the
% displayed pair.}
% \end{figure}

\end{proof}

\begin{prop}\label{prop:cp20-24-subgroup}
Let $a,b$ be homotopy classes in $\tmf^0(\CP^{24+96k}_{19+96k})$ detected by $\Delta^{8k}\bar{\kappa}^2[40+192k]$ and $\Delta^{8k}\eta\eta_1\bar{\kappa}[46+192k]$, such that
\[4a = 2b = q^\ast(\Delta^{8k}4\Delta^2)\]
where $q \colon \CP^{24+96k} \to S^{48+192k}$ is the quotient map.
\begin{enumerate}
    \item The homotopy class $b-2a$ is detected by $\Delta^{8k}\kappa^3[42+192k]$.
    \item The subgroup generated over $\mathbb Z_{(2)}$ by $a$ and $b$ is isomorphic to $\mathbb Z_{(2)}\oplus\mathbb Z/2$. 
\end{enumerate}

\end{prop}

\begin{proof}
\begin{enumerate}
    \item The image of $b-2a$ under 
    \[i^\ast \colon \tmf^0(\CP^{24}_{19}) \to \tmf^0(\CP_{19}^{21})\]
    is $-i^\ast(2a)$, which is detected by $\Delta^{8k}\bar{\kappa}^2\eta^2[42+192k]$ by Proposition~\ref{prop:2 extensions in CP^24_19}. Naturality then gives that $\Delta^{8k}\bar{\kappa}^2\eta^2[42+192k]$ detects $b-2a$ in the $\tmf$-AHSS of $\CP^{24+96k}_{21+96k}$ . 
    \item On the other hand, according to Propositions~\ref{prop:cp20-24-local2} and ~\ref{prop:4 extension in CP^24_20}, $2b$ and $4a$ are both detected by $\Delta^{8k} 4\Delta^2[48+192k]$. $\Delta^{8k} 4\Delta^2[48+192k]$ has top AH filtration, so the difference has to be zero. By a change of basis, $a$ and $b-2a$ are also generators of this subgroup, where the former has infinite order and the latter has order $2$, hence the subgroup is indeed $\mathbb{Z}_{(2)} \oplus \mathbb{Z}/2$.
\end{enumerate}

\end{proof}

\begin{rmk}
From the viewpoint of the Atiyah--Hirzebruch spectral sequence, the class $\Delta^{8k}\kappa^3[42+192k]$ can detect both a class of order two and of infinite order. We record the AH representatives of the following specific elements in this subgroup:
\begin{center}\small
\begin{tabular}{@{}p{7.5cm}ccc@{}}\toprule
Coordinates & Order & AH representative\\\midrule
$(1,0)$ & $\infty$  & $\Delta^{8k}\bar{\kappa}^2[40+192k]$   \\
$(2,0)$ & $\infty$ & $\Delta^{8k}\kappa^3[42+192k]$   \\
$(2,1)$ & $\infty$ & $\Delta^{8k}\eta\eta_1\bar{\kappa}[46+192k]$  \\
$(4,0)$ & $\infty$ & $\Delta^{8k}4\Delta^2[48+192k]$  \\
$(0,1)$ & $2$ & $\Delta^{8k}\kappa^3[42+192k]$  \\\bottomrule
\end{tabular}
\end{center}
\end{rmk}

% \FloatBarrier

\begin{theorem}\label{thm:cp20-24-finite-survival}
For every $k\geq0$, the classes
\[
\Delta^{8k}\bar{\kappa}^2[40+192k],\qquad
\Delta^{8k}\kappa^3[42+192k],\qquad
\Delta^{8k}\eta\eta_1\bar{\kappa}[46+192k]
\]    
are nonzero surviving cycles in the $\tmf$--AHSS of $\CP^{24+96k}$.
\end{theorem}

\begin{proof}

By Lemma~\ref{lem:tmf-integral-no-jump-reformulation}, $\Delta^{8k}\eta\eta_1\bar{\kappa}[46+192k]$ is nonzero in the $\tmf$-AHSS of $\CP^{24+96k}$. In other words, the image of the homotopy class $b$ in Proposition~\ref{prop:cp20-24-subgroup} does not jump filtration.

Suppose, for a contradiction, that the homotopy class $2a$ in Proposition~\ref{prop:cp20-24-subgroup} jumps above its filtration $42+192k$ under the quotient map
\[q^\ast \colon \tmf^0(\CP^{24+96k}_{19+96k}) \to \tmf^0(\CP^{24+96k}).\] As a result,
$q^\ast(2a)$ cannot jump filtration to $44+192k$. Indeed if it does, then it is detected by some $\alpha[44+192k]$ where $\alpha \in \pi_{44+192k}\tmf$ is integral. It follows that $q^\ast(4a)$ must be detected by $2\alpha[44+192k]$, which is a contradiction since $4a$ is detected by $\Delta^{8k}4\Delta^2[48+192k]$ in the $\tmf$-AHSS of $\CP^{24+96k}_{19+96k}$.

If $q^\ast(2a)$ jumps to filtration $46+192k$, then we have 
\[
\bigl\langle q^\ast(2a)|_{\CP^{23+96k}},[\CP^{23+96k}]_{\tmf}\bigr\rangle
=2\bigl\langle q^\ast(a)|_{\CP^{23+96k}},[\CP^{23+96k}]_{\tmf}\bigr\rangle=0,
\]
where the last equality holds because $\pi_{46+192k}\tmf \cong \mathbb{Z}/2$. On the other hand, $q^\ast(2a)$ is detected by $\Delta^{8k}\eta\eta_1\bar{\kappa}[46+192k]$. Thus the pairing must be nonzero, which is a contradiction.

If $q^\ast(2a)$ jumps filtration to $48+192k$ and it is detected by $\beta[48+192k]$, since $q^\ast(4a)$ is detected by $\Delta^{8k}4\Delta^2[48+192k]$, it follows that $2\beta=\Delta^{8k}4\Delta^2 \in \pi_{48+192k}\tmf$, which is a contradiction. In the $\tmf$-AHSS of $\CP^{24+96k}$, the only class in filtration $42+192k$ that can extend to $\Delta^{8k}4\Delta^2[48+192k]$ is $\Delta^{8k}\kappa^3[42+192k]$, hence it is nonzero and detects $q^\ast(2a)$.

Finally, if $q^\ast(a)$ jumps filtration, it can only jump to filtration $42+192k$. $\pi_{42+192k}\tmf$ is torsion, so $q^\ast(2a)$ must have filtration greater than $42+192k$, which is a contradiction. Hence its filtration is $40+192k$. Furthermore, it must be detected by $\Delta^{8k}\bar{\kappa}^2[40+192k]$: The torsion summand of $\pi_{40+192k}\tmf$ is $\mathbb{Z}/4\{\Delta^{8k}\bar{\kappa}^2\}$ where
\[d_2(\Delta^{8k}\eta_1\kappa[38+192k]) = \Delta^{8k}\eta\eta_1\kappa[40+192k] = \Delta^{8k}2\bar{\kappa}^2[40+192k].\]
The free AH cycles in filtration $40+192k$ cannot extend to filtration $48+192k$, hence $q^\ast(a)$ must be detected by $\Delta^{8k}\bar{\kappa}^2[40+192k]$.

\end{proof}

\begin{theorem}\label{thm:tmf-cp20-24}
For every $k\geq0$, $\CP^{m+96k}$ admits an exotic smooth structure for
\[
m\in\{20,22,23,24\}.
\]
\end{theorem}
\begin{proof}
For every $k\geq0$, the classes
\[
\Delta^{8k}\bar{\kappa}^2[40+192k],\qquad
\Delta^{8k}\kappa^3[42+192k],\qquad
\Delta^{8k}\eta\eta_1\bar{\kappa}[46+192k]
\]    
are nonzero in $\tmf^{0}(\CP^{24+96k})$. Theorem~\ref{thm:tmf-detection} thus provides exotic complex projectives in $20+96k$, $21+96k$ and $23+96k$. 

For $m=22+96n$ and $m=24+96n$, by \cite{BMQhurewicz2021} we can choose a Hurewicz preimage $\kappa_n\in\pi_{14+192n}S^0$ whose $\tmf$-Hurewicz image is $\Delta^{8n}\kappa$, with $\kappa_0=\kappa$. By Theorem~\ref{thm:tmf-detection}, $\kappa^2\cdot\kappa_n[42+192k]$ is a surviving cycle in the $\tmf$-AHSS of $\CP^{21+96k}$.

Now in the 2-primary stable stems,
\[
\langle\kappa^2,\eta,\nu \rangle
=\{\Delta h_1h_3\}\cdot\eta+\nu\cdot\pi_{30}S^0.
\]
where $\{\Delta h_1h_3\}\in\pi_{32}S^0$ is a homotopy class detected by the Adams $E_2$-page element $\Delta h_1h_3$. The bracket follows from applying the Moss convergence theorem to $(E_3,d_3)$ where 
\[
d_3(r)=h_1d_0^2,\qquad h_2r=h_1\Delta h_1h_3
\]
and no crossing differentials are involved. Consequently, shuffling Toda brackets gives
\[
\langle \kappa^2\kappa_n,\eta,\nu\rangle
=\kappa_n\cdot\{\Delta h_1h_3\}\cdot\eta
+\nu\cdot\pi_{44+192n}S^0.
\]
Now we show that in $\CP^{24+96n}_{21+96n}$,
\[d_6(\kappa^2\kappa_n[42+192n]) = 0.\]
Hence the homotopy classes detected by $\kappa^2\cdot\kappa_n[42+192n]$ extend to $[\CP^{24+96n},S^0]$ and Theorem~\ref{AHSSsmoothing} will supply the exotic $\CP^{24+96n}$.

\begin{figure}[htbp]\centering
\begin{tikzpicture}[x=.9cm,y=.95cm,scale=.82,transform shape]
\colorlet{twonu}{black}
\colorlet{higher}{gray}
\colorlet{nublue}{black}
\tikzset{cell/.style={circle,draw,fill=white,minimum size=8mm,inner sep=.4pt,font=\small},
 etaedge/.style={black,line width=.8pt},nuedge/.style={nublue,line width=.8pt}}
\foreach \d/\y in {42/0,44/-1.35,46/-2.7,48/-4.05}
 {\node[cell](b\d)at(0,\y){$\d$};}
\draw[etaedge] (b42)--node[right]{$\eta$}(b44);\draw[etaedge] (b46)--node[right]{$\eta$}(b48);
\draw[nuedge] (b44.west) .. controls (-.9,-1.75) and (-.9,-3.65) .. node[midway,left]{$\nu$}(b48.west);
\node at(0,1.05){$\CP_{21}^{24}$};
\end{tikzpicture}
\caption{The stunted projective spectrum at $n=0$. For general $n$,
add $192n$ to each cell dimension and replace the spectrum by $\CP_{21+96n}^{24+96n}$.}
\end{figure}
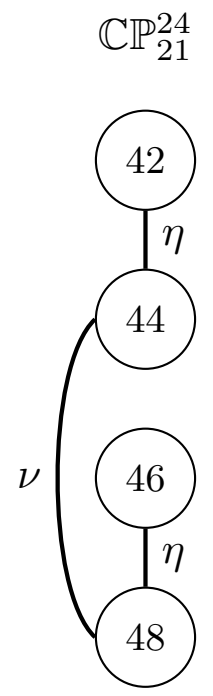

In the AHSS of $\CP^{24+96n}_{21+96n}$, $\kappa^2\cdot\kappa_n[42+192n]$ survives to the $E_6$-page: the $d_2$ vanishes because $\kappa^2\eta=0\in\pi_{29}S^0$, whereas the $d_4$ vanishes because $\CP^{23+96n}_{21+96n}$ splits as $\Sigma^{42+192n}(C\eta\vee S^4)$. $d_6(\kappa^2\cdot\kappa_n[42+192n])$ takes value (modulo earlier boundaries) in the coset $\langle\kappa^2\kappa_n,\eta, \nu \rangle[48+192n]$, where

\[
\langle\kappa^2\kappa_n,\eta, \nu \rangle = \kappa_n\{\Delta h_1h_3\}\eta
+\nu\cdot\pi_{44+192n}S^0.
\]
The first term is killed by the $d_2$ from cell $46+192n$:
\[
d_2\bigl(\kappa_n\{\Delta h_1h_3\}[46+192n]\bigr)
=\kappa_n\{\Delta h_1h_3\}\eta[48+192n].
\]
For every $\alpha \in\pi_{44+192n}S^0$, the class $\alpha[44+192n]$ is a $d_2$-cycle and satisfies
\[
d_4\bigl(\alpha[44+192n]\bigr)=\nu\alpha[48+192n]
\quad\text{modulo $d_2$-boundaries}.
\]
Thus every possible target of the $d_6$ is either a $d_2$-boundary or a $d_4$-boundary, hence $d_6=0$.

\end{proof}

The proof that $\CP^{69+96k}$ to $\CP^{72+96k}$ admits exotic smooth structures is similar to that of $\CP^{20+96k}$ to $\CP^{24+96k}$. More precisely, in both  $\tmf^0(\CP^{24+96k}_{19+96k})$ and $\tmf^0(\CP^{72+96k}_{69+96k})$, there exist subgroups of the form $\mathbb{Z} \oplus \mathbb{Z}/2$, whose generators admit $2$-extensions to $\Delta$-power classes.

\begin{lemma}\label{lem:tmf-kappa138-finite}
In the $\tmf$-AHSS of $\CP^{72+96k}_{69+96k}$, $\Delta^{8k}\kappa^2\kappa_4[138+192k]$, $\Delta^{8k}\epsilon_5\kappa[142+192k]$ and $(4\Delta^6)\Delta^{8k}[144+192k]$ are surviving cycles with hidden $2$-extensions
\[
2\bigl(\Delta^{8k}\kappa^2\kappa_4[138+192k]\bigr)
=2\bigl(\Delta^{8k}\epsilon_5\kappa[142+192k]\bigr)
=(4\Delta^6)\Delta^{8k}[144+192k].
\]
The $\tmf$-cohomology classes detected by these elements generate a subgroup $\mathbb Z_{(2)}\oplus\mathbb Z/2$.
\end{lemma}

\begin{proof}
All the aforementioned classes are AH surviving cycles. According to Lemma~\ref{hiddenextension}, the Toda bracket
\[\langle2,\Delta^{8k}\epsilon_5\kappa,\eta\rangle=\Delta^{8k}4\Delta^6+2\pi_{144+192k}\tmf\]
supplies the hidden $2$-extension $2 \cdot \Delta^{8k}\epsilon_5\kappa[142+192k] =4\Delta^6\Delta^{8k}[144+192k]$ in $\CP^{72+96k}_{71+96k}$, whereas the Toda bracket
\[\langle 2,\Delta^{8k}\kappa_4\kappa^2,\eta,\nu\rangle = \Delta^{8k}\cdot 4\Delta^6\]
supplies the hidden $2$-extension $2 \cdot \Delta^{8k}\kappa^2\kappa_4[138+192k] =\Delta^{8k}4\Delta^6[144+192k]$ in $\CP^{72+96k}_{69+96k}$.

\end{proof}

\begin{prop}\label{prop:cp69-72-subgroup}

Let $a,b$ be homotopy classes in $\tmf^0(\CP^{72+96k}_{69+96k})$ detected by $\Delta^{8k}\kappa_4\kappa^2[138+192k]$ and $\Delta^{8k}\epsilon_5\kappa[142+192k]$, such that
\[2a = 2b = q^\ast(\Delta^{8k}\cdot 4\Delta^6)\]
where $q \colon \CP^{72+96k} \to S^{144+192k}$ is top-cell collapsing.

\begin{enumerate}
    \item The homotopy class $b-a$ is detected by $\Delta^{8k}\kappa_4\kappa^2[138+192k]$.
    \item The subgroup generated over $\mathbb Z_{(2)}$ by $a$ and $b$ is isomorphic to $\mathbb Z_{(2)}\oplus\mathbb Z/2$. 
\end{enumerate}

\end{prop}

\begin{proof}
Building on Lemma~\ref{lem:tmf-kappa138-finite}, the proof is entirely similar to that of Proposition~\ref{prop:cp20-24-subgroup}.
\end{proof}

\begin{theorem}\label{thm:tmf-kappa138-smoothings}
For $k\geq0$ and $m\in\{69,70,71,72\}$, the space $\CP^{m+96k}$ admits an exotic smooth structure.
\end{theorem}
\begin{proof}
Building on Proposition~\ref{prop:cp69-72-subgroup}, the proof is entirely similar to that of Theorem~\ref{thm:tmf-cp20-24}. We replace $\kappa^2\kappa_n[42+192n]$ by $\kappa^2\gamma_k[138+192k]$, where $\gamma_k$ is a lift of $\Delta^{8k}\kappa_4$ by \cite{BMQhurewicz2021}. By almost identical arguments, $\kappa^2\gamma_k[138+192k]$ is a surviving cycle in the AHSS of $\CP^{72+96k}$.
\end{proof}

We conclude this section by showing $\CP^m$ admits an exotic smooth structure for $m=65,66,67,68 (\text{mod } 96)$. Again we work downward inductively, with one distinction from previous arguments: we establish the hidden $2$-extension
\[2 \cdot \Delta^{8k}2\kappa_4\bar{\kappa}[130+192k]=\Delta^{8k}\kappa_4\kappa^2[138+192k]\] 
where $\Delta^{8k}\kappa_4\kappa^2[138+192k]$ detects a torsion element in $\tmf^0(\CP^{69+96k})$. 

One can actually show that there is a hidden $4$-extension 
\[4 \cdot \Delta^{8k}2\kappa_4\bar{\kappa}[130+192k]=\Delta^{8k}\cdot 4\Delta^6[144+192k]\] 
in the $\tmf$-AHSS of $\CP^{72+96k}_{65+96k}$, which would nonetheless require a longer proof. For simplicity, we establish the extensions successively.

\begin{lemma}\label{lem:tmf-w2-shapes}
In the $\tmf$-AHSS of $\CP_{65+96k}^{68+96k}$ there is a hidden $2$-extension
\[2 \cdot \Delta^{8k}w^2[130+192k]=2\Delta^4\Delta^{8k}[136+192k]\] 
\end{lemma}
\begin{proof}
Both $\Delta^{8k}w^2[130+192k]$ and $2\Delta^4\Delta^{8k}[136+192k]$ are surviving cycles. A variation of Lemma~\ref{hiddenextension} shows that should hidden $2$-extension from the bottom cell to the top cell exists, it takes value in
\[
\left\langle2,\Delta^{8k}w^2,(\eta,2\nu),
       \begin{pmatrix}\nu\\\eta\end{pmatrix}\right\rangle, 
\]
which decomposes as $\langle 2,\Delta^{8k}w^2,\eta,\nu\rangle + \langle 2,\Delta^{8k}w^2,2\nu,\eta\rangle$. By Lemma~\ref{lem:90 to 96 extension}, the first bracket is $2\Delta^{8k+4}+2\pi_{96}\tmf$ and the second bracket is $2\pi_{96}\tmf$, by which the hidden extension follows.
\end{proof}

\begin{lemma}\label{lem:tmf-65-to-69 extension}
In the $\tmf$-AHSS of $\CP_{65+96k}^{69+96k}$ there is a hidden $2$-extension
\[2 \cdot \Delta^{8k}2\kappa_4\bar{\kappa}[130+192k]=\Delta^{8k}\kappa_4\kappa^2[138+192k]\]     \end{lemma}

\begin{proof}
In the $\tmf$-AHSS of $\CP_{65+96k}^{68+96k}$, $\Delta^{8k}2\kappa_4\bar{\kappa}[130+192k]$ is a permanent cycle. The relation $\bar{\kappa}^2w^2=2\kappa_4\bar{\kappa}$ thus provides
\begin{equation}\label{barkappa 1}
\bar{\kappa}^2 \cdot \Delta^{8k}w^2[130+192k] = \Delta^{8k}2\kappa_4\bar{\kappa}[130+192k].    
\end{equation}

On the other hand, in the long exact sequence
\[\pi_{135}\tmf\xrightarrow{(\eta,u\nu)}
\pi_{136}\tmf\oplus\pi_{138}\tmf
\xrightarrow{q^*}\tmf^0(\CP_{67}^{69}) \to 0,\]

the image of $\eta_1\kappa_4 \in \pi_{135}\tmf$ is $(\eta\eta_1\kappa_4,\kappa^2\kappa_4)$. Using the multiplicative relation
\[(2\Delta^4)\bar{\kappa}^2=\eta\eta_1\kappa_4\]
we obtain the relation
\begin{equation}\label{barkappa 2}
\bar{\kappa}^2 \cdot \Delta^{8k}\cdot 2\Delta^4[136+192k] = \Delta^{8k}\kappa_4\kappa^2[138+192k]    
\end{equation}
in the $\tmf$-AHSS of $\CP_{67+96k}^{69+96k}$. Combining \eqref{barkappa 1} and \eqref{barkappa 2} and Lemma~\ref{lem:tmf-w2-shapes}, the extension in $\CP^{69+96k}_{65+96k}$ thus follows from naturality. 
   
\end{proof}

% Building on Theorem~\ref{thm:tmf-kappa138-smoothings}, the 2-extension above will supply the survival of $\Delta^{8k}2\kappa_4\bar{\kappa}[130+192k]$ in the $\tmf$-AHSS of $\CP^{65+96k}$. We still need an analogue of Lemma~\ref{lem:tmf-integral-no-jump-reformulation} where we replace the delta power classes on top cell with a torsion class, which has a priori been shown to survive .

% \begin{lemma}\label{lem:tmf-topcell-no-jump-reformulation}
% Let $\alpha[2N]$ be a surviving cycle in the $\tmf$-AHSS of $\CP^N$. Suppose in the $\tmf$-AHSS of $\CP^N_n$, $\beta[2n]$ admits a $2^s$-extension to 
% $\alpha[2N]$. Suppose also that multiplication by $2^{s-1}$ in $\tmf^0(\CP^N)$
% sends every homotopy class of AH filtration greater than $2n$ into filtration $2N$.
% Then $\beta[2n]$ survives in the $\tmf$-AHSS of $\CP^N$; in other words, the homotopy classes detected by $\beta[2n]$ maps nonzeroly under 
% \[q^\ast:\tmf^0(\CP^N_n) \to \tmf^0(\CP^N),\]
% where they are detected by the elements with the same AH names. 
% \end{lemma}

% \begin{proof}

% The proof is similar to that of Lemma~\ref{lem:tmf-integral-no-jump-reformulation}.

% \end{proof}

\begin{theorem}\label{thm:tmf-five-cell-detection}
For every $k\geq0$ and $m\in\{65,66,67,68\}$, the space $\CP^{m+96k}$ admits an exotic smooth structure.

\end{theorem}
\begin{proof}

Under the quotient map $q \colon \CP^{69+96k} \to \CP^{69+96k}_{65+96k}$
\[q^\ast \colon \tmf\text{-AHSS of } \CP ^{69+96k}_{65+96k} \to \tmf\text{-AHSS of } \CP ^{69+96k}\]
sends $\Delta^{8k}\kappa_4\kappa^2[138+192k]$ to the element of the same name, where Theorem~\ref{thm:tmf-kappa138-smoothings} shows the survival of $\Delta^{8k}\kappa_4\kappa^2[138+192k]$ in the latter spectral sequence. It follows that the homotopy class detected by $\Delta^{8k}\cdot 2\kappa_4\bar{\kappa}[130+192k]$ in $\tmf^0(\CP^{69+96k}_{65+96k})$ maps nontrivially to $\tmf^0(\CP^{69+96k})$. By an argument similar to that of Theorem~\ref{thm:tmf-cp20-24} such a homotopy class cannot jump filtration. Thus it has filtration $130+192k$ and is detected by $\Delta^{8k}2\kappa_4\bar{\kappa}[130+192k]$.

By Theorem~\ref{thm:tmf-detection}, $\bar{\kappa}^2\rho_k[130+192k]$ is nonzero in the AHSS of $\CP^{65+96k}$, where $\rho_k \in \pi_{90+192k}S^0$ is the specified lift in Theorem~\ref{thm:tmf-cp-46-47}. $\bar{\kappa}^2\rho_k[130+192k]$ is nonzero in the AHSS of $\CP_{65+96k}^{68+96k}$, Theorem~\ref{AHSSsmoothing} hence provides the exotic smooth structures.

\end{proof}

\begin{proof}[Proof of Theorem~\ref{thm:delta method}]
Combine Theorems~\ref{thm:tmf-cp-10-45-58-82}, \ref{thm:tmf-cp-46-47}, \ref{thm:tmf-cp20-24}, \ref{thm:tmf-kappa138-smoothings}, and \ref{thm:tmf-five-cell-detection}.
\end{proof}

\section{Computations of the algAHSS and $\beta$-families}
\label{sec:anss-detection}

In this section we work $2$-locally and apply the $\beta$-families in the algAHSS method outlined in Section~\ref{section: ANSS methodology}. Subsection~\ref{subsec:anss-cobar} describes the construction of the cobar lifts of the $d_2$-algAHSS cycles of tridegree $(2\ell,1,2)$, by which the survival of $\beta[2n]$ to the algAHSS $E_\infty$ page is guaranteed. Subsection~\ref{subsec:anss-permanence} establishes its permanence in the ANSS and translates between ANSS and AHSS, providing nonzero $b[2n]$'s in the AHSS of $\CP^n$. Then we extend the nonzero $b[2n]$'s further in the AHSS of $\CP^{n+r}$ in Subsection~\ref{subsec:anss-dimensions}. In summary, we obtain the following complex dimensions $n$ for which $\CP^n$ admits an exotic smooth structure:

\begin{theorem}\label{thm:anss-main}
For every integer $k\geq0$, the space $\CP^{m+96k}$ admits an exotic smooth structure for every integer $m$ in
\[
\{8,12,15,19,31,32,34,36,39,40,43,44,60,87,88,90,91,92,94,95,96,98\}.
\]
\end{theorem}

% The geometric proofs are in Subsection~\ref{subsec:anss-dimensions}.

\subsection{$1$-line cobar lifts and survival of top-cell $\beta$-classes}\label{subsec:anss-cobar}

The purpose of this subsection is to prove the following theorem:

\begin{theorem}\label{thm:anss-top-edge}
Consider the algAHSS converging to $\Ext_{BP_*BP}^{\ast,\ast}(BP_\ast \CP^n,BP_\ast)$. Let $\beta \in \Ext_{BP_*BP}^{2,2n+2}(BP_\ast,BP_\ast)$ be an element of the algebraic $\beta$-family. Then the top-cell class $\beta[2n]\in E_1^{2n,2,2}$ survives nontrivially to $E_\infty$, detecting a nonzero class in $\Ext_{BP_*BP}^{2,2}(BP_\ast \CP^n,BP_\ast)$.
\end{theorem}

If $\beta[2n]$ is hit by an algAH differential, then the source must have tridegree $(2\ell,1,2)$ where $\ell < n $. In other words, it takes the form $\alpha[2\ell]$ where $\alpha \in \Ext_{BP_*BP}^{1,2\ell+2}(BP_\ast,BP_\ast)$. The ANSS 1-line of the sphere consists of the $\alpha$-family elements, namely
\begin{equation*}
\Ext^{1,2n}_{BP_*BP}(BP_\ast, BP_\ast)=
\begin{cases}
\mathbb{Z}/2\{\alpha_n\},&n>0\text{ odd},\\
\mathbb{Z}/4\{g_2\},&n=2,\\
\mathbb{Z}/2^{\nu_2(n)+2}\{g_n\},&n\geq4\text{ even}.
\end{cases}
\end{equation*}

There could be nonzero $d_2$-algAH differentials supported on the $\alpha[2\ell]$'s. The $d_2$-differentials are induced by $\alpha_1$-multiplications, namely 
\begin{equation*}
d_2(\gamma[2\ell])=\pm\ell\alpha_1\gamma[2\ell+2],
\qquad \alpha_1 \in \Ext^{1,2}_{BP_\ast BP}(BP_\ast, BP_\ast).
\end{equation*}
For even $\ell$, the differential vanishes because $2\alpha_1=0$. For odd $\ell$, we have $\alpha_1g_2=0$, whereas $\alpha_1g_{\ell+1}$ is nonzero of order $2$ when $\ell\geq3$. These nonzero products lie in different direct summands from the algebraic $\beta$-family elements \cite{ravenel2004complexcobordism}. For each integer $\ell\geq1$, let $K_\ell$ denote the kernel of $d_2$ in tridegree $(2\ell,1,2)$ in the algAHSS for $\CP^\infty$. It follows that
\begin{equation}\label{eq:anss-first-line-kernels}
K_\ell=
\begin{cases}
\mathbb{Z}/2\{\alpha_{\ell+1}\},&\ell\text{ even},\\
\mathbb{Z}/4\{g_2\},&\ell=1,\\
\mathbb{Z}/2^{\nu_2(\ell+1)+1}\{2g_{\ell+1}\},&\ell\geq3\text{ odd}.
\end{cases}
\end{equation}
For $1\leq\ell<n$, restriction to $\CP^n$ identifies $K_\ell$ with the corresponding $d_2$-kernel in the finite algAHSS.

We prove that all elements of the groups $K_\ell$ are permanent algAHSS cycles. The inclusion $BP_\ast(\CP^n) \to BP_\ast(\CP^\infty)$ induces a map between the algAHSSs
\begin{center}
\begin{tikzcd}
{\bigoplus_{i=1}^\infty \Ext_{BP_\ast BP}^{s,t+2i}(BP_\ast,BP_\ast)}  \arrow[rr] \arrow[d, "algAHSS" description, Rightarrow]  &  & {\bigoplus_{i=1}^n \Ext_{BP_\ast BP}^{s,t+2i}(BP_\ast,BP_\ast)} \arrow[d, "algAHSS" description, Rightarrow] \\
{\Ext_{BP_\ast BP}(BP_\ast \CP^\infty, BP_\ast)} \arrow[rr] &  & {\Ext_{BP_\ast BP}(BP_\ast \CP^n, BP_\ast).}    
\end{tikzcd}
\end{center}
By naturality, if, for every $\ell \in \mathbb{N}$, all elements $k \in K_\ell$ survive in the algAHSS converging to $\Ext_{BP_\ast BP}^{s,t}(BP_\ast\CP^\infty, BP_\ast)$, then, for all $\ell < n$, all elements $k \in K_\ell$ must survive in the algAHSS converging to $\Ext_{BP_\ast BP}^{s,t}(BP_\ast\CP^n, BP_\ast)$. We show all elements $k \in K_\ell$ admit lifts to the cobar complex of $BP_\ast(\CP^\infty)$. 

\begin{theorem}\label{thm:CP infty cobar lift}
For every integer $\ell\geq1$, every element of $K_\ell$ defined in \eqref{eq:anss-first-line-kernels} admits a lift to the cobar complex 
\[ \mathcal{C}^{s,t}_{BP_\ast BP}(BP_\ast \CP^\infty, BP_\ast) \cong \Hom_{BP_\ast}^t (BP_\ast \CP^\infty,(BP_\ast BP)^{\otimes_{BP_\ast}^s}) \]
where $(s,t) = (1,2)$. In other words, every such $d_2$-cycle is a permanent cycle. Consequently, for every $n\geq1$ and $1\leq\ell<n$, all elements of $K_\ell$ are also permanent cycles in the algAHSS converging to $\Ext^{1,2}_{BP_\ast BP}(BP_\ast\CP^n,BP_\ast)$.
  
\end{theorem}

\begin{proof}
    The cobar lifts are constructed in Proposition~\ref{prop:ordinary} and Proposition~\ref{prop:divided}.
\end{proof}

\begin{proof}[Proof of Theorem~\ref{thm:anss-top-edge}]
$\beta[2n]$ is a permanent cycle in the algAHSS for filtration reasons. Furthermore, Theorem~\ref{thm:CP infty cobar lift} excludes the possibility that $\beta[2n]$ is hit by an algAH differential. Thus it survives to $E_\infty$ in the algAHSS, detecting some nonzero element in $\Ext^{2,2}_{BP_\ast BP}(BP_\ast \CP^n, BP_\ast)$.
\end{proof}

The remainder of this subsection will be devoted to constructing explicit cobar representatives. At first glance, passing from $\CP^n$ to $\CP^\infty$ may seem to make the computation more difficult. However, the infinite setting allows us to organize the cobar calculations using formal power series associated with the universal $2$-typical formal group law and then recover the finite cases by truncation. Our construction is motivated by Miller's computation of universal Bernoulli numbers and the algebraic $S^1$-transfer \cite{Miller1982Bernoulli}. The short exact sequence of $BP_\ast BP$-comodules
\[
0\longrightarrow \Sigma^{-2}BP_\ast
\longrightarrow BP_\ast\CP^\infty_{-1}
\longrightarrow BP_\ast(\CP^\infty_+)
\longrightarrow 0
\]
induces a long exact sequence of Ext groups. The algebraic transfer is the extension class in $\Ext^{1,2}_{BP_\ast BP}(BP_\ast(\CP^\infty_+),BP_\ast)$ associated with the above short exact sequence. Yoneda composition with this class gives the connecting homomorphisms in the resulting long exact sequence of Ext groups. We denote its restriction to the reduced summand by
\[
e_\tau\in
\Ext^{1,2}_{BP_\ast BP}(BP_\ast\CP^\infty,BP_\ast).
\]
A natural question is how to identify the algAHSS name of $e_\tau$. To answer it, we first describe the filtration on the cobar complex that gives rise to the algAHSS and then apply the rational Bockstein.

We recall some basic constructions and set up the following notation. The cobar complex of $BP_\ast\CP^\infty$ is
\[ 
\begin{aligned}
\mathcal C_{BP_\ast BP}^{s,t}(BP_*\CP^\infty,BP_*)
&=
\Hom_{BP_\ast}^t (BP_\ast \CP^\infty,(BP_\ast BP)^{\otimes_{BP_\ast}^s})\\
&=
\left\{\sum_{j\geq 1}u_jz^j, u_j \in ((BP_\ast BP)^{\otimes_{BP_\ast}^s})_{t+2j}\right\}, |z| =-2.
\end{aligned}
\]

Its cohomology is $\Ext^{s,t}_{BP_\ast BP}(BP_*\CP^\infty,BP_*)$. Every cobar differential has bidegree $(1,0)$:
\[
 d^s:\mathcal C^{s,t}_{BP_\ast BP}(BP_*\CP^\infty,BP_*)
 \longrightarrow\mathcal C^{s+1,t}_{BP_\ast BP}(BP_*\CP^\infty,BP_*),
 \qquad d^{s+1}d^s=0.
\]

In particular, for $s=0$, elements in $\mathcal{C}^{0,t}$ take the form $U(z)=\sum_{j\geq1}u_jz^j$, with $u_j\in BP_{t+2j}$. For $j\geq1$, let $b_j\in BP_{2j}(\CP^\infty)$ be dual to $x^j$, where $x\in BP^2(\CP^\infty)$ is the chosen orientation. We identify zero-cochains of internal degree $t$ with formal power series by
\begin{equation}\label{eq:cobar-power-series-correspondence}
\Phi_U \longleftrightarrow U(z)=\sum_{j\geq1}u_jz^j,
\qquad
\Phi_U(b_j)=u_j\in (BP_\ast)_{t+2j}.
\end{equation}
The explicit formula for the cobar differential $d^0$ is
\[
 d^0U(z)=U_R(S(z))-U_L(z),\qquad
 U_L(z)=\sum_j\eta_L(u_j)z^j,\quad U_R(z)=\sum_j\eta_R(u_j)z^j.
\]
$S(z)$ is the strict transport characterized by
\[
 \lambda_R(S(z))=\lambda_L(z), \qquad S(z)=z-t_1z^2+O(z^3).
\]

Write $(BP_\ast)_{\mathbb Q}= (BP_\ast)\otimes_{\mathbb Z}\mathbb Q$, and regard $(BP_\ast)_{\mathbb Q}$ and $(BP_\ast)_{\mathbb Q}/BP_\ast$ as $BP_\ast BP$-comodules with their induced coactions. Similarly one can define the cobar complex with coefficients in either $(BP_\ast)_{\mathbb Q}$ or $(BP_\ast)_{\mathbb Q}/BP_\ast$, regarded as $BP_\ast BP$-comodules. The short exact sequence of $BP_\ast BP$-comodules
\[
0\longrightarrow BP_\ast
\longrightarrow (BP_\ast)_{\mathbb Q}
\longrightarrow (BP_\ast)_{\mathbb Q}/BP_\ast
\longrightarrow 0
\]
induces the rational Bockstein homomorphism
\[
\delta_{\mathbb Q}:
\Ext_{BP_\ast BP}^{s,t}
\bigl(-,(BP_\ast)_{\mathbb Q}/BP_\ast\bigr)
\longrightarrow
\Ext_{BP_\ast BP}^{s+1,t}(-,BP_\ast).
\]

We now equip the cobar complex with a filtration that gives rise to the algAHSS of Subsection~\ref{section: ANSS methodology}.

\begin{prop}\label{prop: algAHSS by cobar filtration}
The algAHSS can be obtained by filtering the cobar complex by
\[
F^{2i}\mathcal C_{BP_\ast BP}^{s,t}(BP_*\CP^\infty,BP_*)= \{g=\sum_{j\geq i}u_jz^j: u_j \in ((BP_\ast BP)^{\otimes_{BP_\ast}^s})_{t+2j} \}, \qquad F^{2i+1}=F^{2i+2}.
\]
Restriction to $\CP^N$ is truncation modulo $z^{N+1}$, which recovers the finite spectral sequence of Subsection~\ref{section: ANSS methodology}. In particular, let $H(z)=az^m+O(z^{m+1})$ be a rational zero-cochain of internal degree $2$, where $m\geq1$ and $|a|=2+2m$. If $H$ has integral coboundary, then $S(z)=z+O(z^2)$ gives
\[
d^0H=(\eta_R(a)-\eta_L(a))z^m+O(z^{m+1}).
\]
If, moreover, $\delta_{\mathbb Q}(\bar a)\neq0$, then the class $[d^0H]$ has algAH filtration $2m$ and algAHSS name $\delta_{\mathbb Q}(\bar a)[2m]$. Here $\bar a$ denotes the reduction of $a$ modulo $BP_\ast$.
\end{prop}

\begin{proof}
On the associated graded complex, the cobar differential acts only on the coefficients, so
\[
\operatorname{gr}^{2i}_F
\mathcal C^{s,t}_{BP_*BP}(BP_*\CP^\infty,BP_*)
\cong
\mathcal C^{s,t+2i}_{BP_*BP}(BP_*,BP_*).
\]
This is the cobar construction of the algAHSS. Restriction to $\CP^N$ discards the terms of degree greater than $N$.

For the second assertion, the formula $d^0H=H_R(S(z))-H_L(z)$ gives the stated leading term. Its integrality implies that $\bar a$ is invariant, and $\eta_R(a)-\eta_L(a)$ represents $\delta_{\mathbb Q}(\bar a)$. There are no incoming differentials to this nonzero class, since the $0$-line of the ANSS of the sphere is zero for $t-s >0$. Moreover, the integral cocycle $d^0H$ supplies a lift through every page. Its restriction to $\CP^m$ is therefore nonzero, whereas its restriction to $\CP^{m-1}$ vanishes. Hence its algAH filtration is  $2m$, with the asserted algAHSS name.
\end{proof}

We now express $e_\tau$ as the image under $\delta_{\mathbb Q}$ of a class in
\[
\Ext^{0,2}_{BP_\ast BP}
\bigl(BP_\ast\CP^\infty,(BP_\ast)_{\mathbb Q}/BP_\ast\bigr).
\]
Such a class is a $BP_\ast BP$-comodule homomorphism of internal degree $2$ from $BP_\ast\CP^\infty$ to $(BP_\ast)_{\mathbb Q}/BP_\ast$. The formal group logarithm provides an explicit description of this homomorphism.

\begin{defn}\label{def: algebraic transfer Tk classes }
Define the reciprocal logarithm by
\begin{equation}\label{reciprocal logarithm}
P_F(z)=\frac1{\lambda_F(z)}-\frac1z+\lambda_1
=\sum_{j\geq1}c_jz^j,
\end{equation}
where 
\[\lambda_F(z) = \sum^\infty_{i=0} \lambda_i z^{2^i}\]  
is the logarithm of the universal $2$-typical formal group law, expressed in Hazewinkel generators. Consequently, the coefficients $c_j$ lie in $(BP_\ast)_{\mathbb{Q}}$.
\end{defn}

\begin{prop}
The reduction $\overline{P_F}$ of $P_F$ modulo $BP_\ast$ represents a class in
\[
\Ext^{0,2}_{BP_\ast BP}
\bigl(BP_\ast\CP^\infty,(BP_\ast)_{\mathbb Q}/BP_\ast\bigr)
\]
whose image under $\delta_{\mathbb Q}$ is $e_\tau$. Moreover, $e_\tau$ has algAHSS name $g_2[2]$.
\end{prop}

\begin{proof}
Miller's description of the algebraic transfer \cite{Miller1982Bernoulli} identifies
\[
e_\tau
=\delta_{\mathbb Q}(\overline{P_F})
=[d^0P_F].
\]
Expanding $P_F(z)$ gives
\[
P_F(z)=\frac{v_1^2}{4}z+O(z^2),
\qquad
\delta_{\mathbb Q}
\!\left(\overline{\frac{v_1^2}{4}}\right)=g_2.
\]
Thus by Proposition~\ref{prop: algAHSS by cobar filtration}, $e_\tau$ has algAHSS name $g_2[2]$.
\end{proof}

\begin{rmk}
Under the identification \eqref{eq:cobar-power-series-correspondence}, the series $P_F(z)=\sum_{j\geq1}c_jz^j$ corresponds to the $BP_\ast$-linear map
\[
\Phi_F:
BP_\ast\CP^\infty
\longrightarrow \Sigma^{-2}(BP_\ast)_{\mathbb Q},
\qquad
b_j\longmapsto c_j.
\]
We extend $\Phi_F$ to $BP_\ast(\CP^\infty_+)$ by setting $\Phi_F(b_0)=0$, and denote its reduction modulo $BP_\ast$ by $\overline{\Phi_F}$. One can further apply $\Phi_F$ coefficientwise to homology-valued formal power series; in particular,
\[
\Phi_F\left(\sum_{j\geq0}b_jz^j\right)
=
\sum_{j\geq1}c_jz^j
=
P_F(z).
\]
\end{rmk}

\begin{prop}\label{prop:integral}
For every integer $a\geq1$, $d^0\bigl(P_F([a]_F(z))\bigr)$ is an integral cobar one-cocycle. Here $[a]_F$ denotes the multiplication-by-$a$ series of $F$.
\end{prop}

\begin{proof}
First, we compute
\begin{align*}
d^0P_F
&=\left(\frac1{\lambda_R(S(z))}-\frac1{S(z)}+\lambda_{1,R}\right)
 -\left(\frac1{\lambda_L(z)}-\frac1z+\lambda_{1,L}\right)\notag\\
&=\frac1z-\frac1{S(z)}+t_1\in z(BP_\ast BP)[[z]].
\end{align*}
The last inclusion follows from the integrality and strictness of $S$, since $1/S(z)=1/z+t_1+O(z)$.

For any integer $a\geq1$, the $a$-series is natural under strict isomorphisms:
\[
[a]_{F_R}(S(z))=S([a]_{F_L}(z)).
\]
Consequently
\begin{align*}
d^0\bigl(P_F([a]_F(z))\bigr)
&=P_R\bigl([a]_{F_R}(S(z))\bigr)-P_L\bigl([a]_{F_L}(z)\bigr)\\
&=(d^0P_F)([a]_{F_L}(z)).
\end{align*}
The right-hand side is integral because both series being composed are integral and the $a$-series has zero constant term.  
\end{proof}

Now we are ready to construct cobar lifts of the elements of $K_\ell$. 

% \begin{rmk}
% For odd $b>0$, one has
% \begin{equation}\label{eq:oddintegral}
% P_F([b]_F(z))-b^{-1}P_F(z)
% =\frac1{bz}-\frac1{[b]_F(z)}+\lambda_1(1-b^{-1})\in zR[[z]].
% \end{equation}
% To check integrality, write
% $[b]_F(z)=bz+(b-b^2)\lambda_1z^2+O(z^3)$. Its linear coefficient is a unit, so all positive coefficients of the reciprocal difference are integral. The displayed correction cancels the constant.  
% \end{rmk}

\begin{lemma}\label{lem:denominator}
For the coefficients $c_j$ in \eqref{reciprocal logarithm},
\begin{equation}\label{eq:denominator}
2^{j+1}c_j\in BP_\ast,
\qquad
2^jc_j\equiv(-1)^{j+1}\frac{v_1^{j+1}}2\pmod{BP_\ast} .
\end{equation}
\end{lemma}
\begin{proof}
Expand
\[
\frac1{\lambda_F(z)}
=\frac1z\left(1+\sum_{i\geq1}\lambda_i z^{2^i-1}\right)^{-1}.
\]
A monomial $\lambda_{i_1}\cdots\lambda_{i_q}$ contributing to $c_j$ satisfies
\[
\sum_{h=1}^q(2^{i_h}-1)=j+1.
\]
By Hazewinkel's recursion we have $2^i\lambda_i\in BP_\ast$; see \cite{ravenel2004complexcobordism}. If $i_h\geq2$ for some $h$, then $\sum_h i_h<\sum_h(2^{i_h}-1)=j+1$, so $2^j\prod_h\lambda_{i_h}\in BP_\ast$. Thus only the term with all $i_h=1$ contributes to $2^jc_j$ modulo $BP_\ast$. Its contribution is $(-1)^{j+1}2^j\lambda_1^{j+1}=(-1)^{j+1}v_1^{j+1}/2$, proving both assertions.
\end{proof}

\begin{prop}\label{prop:ordinary}
For every even $n\geq2$, the class $2^nd^0P_F$ has algAHSS filtration $2n$, with leading symbol $\alpha_{n+1}[2n]$. 
\end{prop}

\begin{proof}
By \eqref{eq:denominator}, the truncation $Q_{n}(z)=\sum_{j=1}^{n-1}2^nc_jz^j$ is integral. Subtracting $d^0Q_{n}$ therefore gives a representative
\[
2^nd^0P_F-d^0Q_{n}
=d^0\bigl(2^nP_F-Q_{n}\bigr).
\]
By Lemma~\ref{lem:denominator}, the leading term has degree $n$, with coefficient $2^nc_n\equiv-v_1^{n+1}/2\pmod{BP_*}$. Consequently, by Proposition~\ref{prop: algAHSS by cobar filtration}, the algAH filtration of this cocycle is precisely $2n$, and its algAH name is $\alpha_{n+1}[2n]$, since $\delta_{\mathbb Q}\!\left(\overline{v_1^{n+1}/2}\right) = \alpha_{n+1}$.
\end{proof}

Constructing cobar representatives of $2g_{n+1}[2n]$ and its nonzero multiples by powers of $2$ in the same algAH filtration requires a more involved argument. By Proposition~\ref{prop: algAHSS by cobar filtration}, it suffices to find a formal power series $G_n(z)$ with coefficients in $(BP_\ast)_{\mathbb Q}$ and internal degree $2$ such that
\begin{itemize}
    \item $d^0G_n$ is integral;
    \item $G_n(z)=a_nz^n+O(z^{n+1})$, where $\delta_{\mathbb Q}(\overline{a_n})=2g_{n+1}$.
\end{itemize}
One natural candidate for $a_n$ follows from Miller's universal Bernoulli coefficients in \cite{Miller1982Bernoulli}. Let $e_F$ be the inverse of $\lambda_F$. The universal Bernoulli coefficients are defined by
\begin{equation*}
\frac{u}{e_F(u)}
=
\sum_{r\geq0}\frac{B_r(F)}{r!}u^r.
\end{equation*}
Consequently,
\begin{equation*}
P_F(e_F(u))
=
\frac1u-\frac1{e_F(u)}+\lambda_1
=
-\sum_{r\geq1}
\frac{B_{r+1}(F)}{(r+1)!}u^r.
\end{equation*}
Furthermore, Miller shows that for every even $w\geq2$,
\begin{equation*}
\overline{\frac{B_w(F)}{w}}
\in
\Ext_{BP_\ast BP}^{0,2w}
\bigl(BP_*,(BP_*)_{\mathbb Q}/BP_*\bigr)
\quad\text{has order }2^{\nu_2(w)+1}.
\end{equation*}
This motivates us to seek a suitable formal power series with the required leading coefficient, as in the following proposition.

\begin{prop}\label{prop:divided}
For odd $n\geq3$, set
\begin{equation}\label{eq:Gn}
G_n(z)
=
-\sum_{a=1}^{n}(-1)^{n-a}\binom{n}{a}
P_F([a]_F(z)).
\end{equation}
The class $d^0G_n$ is integral. It has algAH filtration $2n$, and its algAH name is $2g_{n+1}[2n]$.
\end{prop}

\begin{proof}
Indeed, $d^0G_n$ is integral by Proposition~\ref{prop:integral}. Equip $BP_\ast(\CP^\infty_+)$ with the Pontryagin product and write
\[
\beta(z)=\sum_{j\geq0}b_jz^j.
\]
Here $b_0=1$ is the unit for the Pontryagin product. The series satisfies the product formula
\[
\beta(z)\beta(w)=\beta(F(z,w)).
\]
In particular,
\[
\beta(z)^a=\beta([a]_F(z))
\qquad(a\geq0).
\]
Applying $\Phi_F$ coefficientwise and using $\Phi_F(\beta(z))=P_F(z)$, we obtain from the binomial theorem
\[
\begin{aligned}
-\Phi_F\bigl((\beta(z)-1)^n\bigr)
&=
-\sum_{a=0}^{n}(-1)^{n-a}\binom{n}{a}
  \Phi_F\bigl(\beta(z)^a\bigr)\\
&=
-\sum_{a=0}^{n}(-1)^{n-a}\binom{n}{a}
  P_F([a]_F(z))\\
&=G_n(z).
\end{aligned}
\]

To determine the leading coefficient, observe that
\[
(\beta(z)-1)^n=b_1^nz^n+O(z^{n+1}),
\]
and hence
\[
G_n(z)=-\Phi_F(b_1^n)z^n+O(z^{n+1}).
\]
Thus it remains to compute $\Phi_F(b_1^n)$. Consider the series
\[
\beta(e_F(u))=\sum_{r\geq0}\frac{p_r}{r!}u^r.
\]
Miller identifies $p_r$ with $b_1^r$ \cite{Miller1982Bernoulli}. Extending $\Phi_F$ rationally therefore gives
\[
\sum_{r\geq1}\frac{\Phi_F(b_1^r)}{r!}u^r
=
\Phi_F\bigl(\beta(e_F(u))\bigr)
=
P_F(e_F(u)) 
= -\sum_{r\geq1}
\frac{B_{r+1}(F)}{(r+1)!}u^r.
\]
It follows that
\begin{equation*}
G_n(z)
=
\frac{B_{n+1}(F)}{n+1}z^n+O(z^{n+1}).
\end{equation*}

Consequently, Proposition~\ref{prop: algAHSS by cobar filtration} shows that $d^0G_n$ has algAH filtration $2n$. Furthermore, its algAH name is
\[
\delta_\mathbb Q \left(
\overline{\frac{B_{n+1}(F)}{n+1}}
\right)[2n]
=
2g_{n+1}[2n].
\]
\end{proof}

% \begin{center}
% \renewcommand{\arraystretch}{1.2}
% \begin{tabular}{@{}crrrl@{}}
% \toprule
% $n$ & $A_{n,0}$ & $A_{n,1}$ & $A_{n,2}$ & Leading symbol\\
% \midrule
% $3$ & $-10/3$     & $3$    & $0$  & $2g_4[6]$\\
% $5$ & $-128/15$   & $10$   & $5$  & $2g_6[10]$\\
% $7$ & $-2416/105$ & $70/3$ & $35$ & $2g_8[14]$\\
% \bottomrule
% \end{tabular}
% \end{center}
% For example,
% \[
% X_3=3\mathcal T_1-\frac{10}{3}\mathcal T_0,
% \qquad
% X_5=5\mathcal T_2+10\mathcal T_1-\frac{128}{15}\mathcal T_0.
% \]

\subsection{From ANSS detection to AHSS detection}
\label{subsec:anss-permanence}

Next we trace through the Mahowald square and show that the homotopy class detected by $b[2n]$ of Theorem~\ref{thm:anss-top-edge} is detected by $\beta[2n]$ in the AHSS.

\begin{prop}\label{prop:anss-top-cell-detection} 
Use the notation in Theorem~\ref{thm:anss-top-edge}. Let $b \in\pi_{2n}S^0_{(2)}$ be a homotopy class detected by a permanent $\beta$-element $\beta \in\Ext_{BP_\ast BP}^{2,2n+2}(BP_*,BP_*)$. The class $b[2n]$ survives in the AHSS of $\CP^n$.
\end{prop}

\begin{proof}

\[
\begin{tikzcd}
{\Ext^{s,t+2n}_{BP_\ast BP}(BP_\ast,BP_\ast)} \arrow[rrr, "q^\ast"] \arrow[dd, "ANSS" description, Rightarrow] &  &  & {\Ext^{s,t}_{BP_\ast BP}(BP_\ast\CP^n,BP_\ast)} \arrow[dd, "ANSS" description, Rightarrow] \\
                                                                      &  &  &                                               \\
{\pi_{t+2n-s}S^0} \arrow[rrr, "q^\ast"']                                           &  &  & {[\CP^n,\Sigma^{s-t}S^0]}                                           
\end{tikzcd}
\]

By naturality the class $q^\ast(\beta)$ is detected by $b[2n]$ in the ANSS of $\CP^n$, which is detected by $\beta[2n]$ in the AHSS of $\CP^n$.

\end{proof}

\subsection{Extending $\beta$-family classes from the top cell}\label{subsec:anss-dimensions}

We first give a brief account of some $2$-primary $\beta$-families in the stable stems, which we will use for the detection of exotic $\CP^n$. Their detailed construction will appear in an upcoming paper by Li-Ma-Xu \cite{LiMaXu2026v2periodic}.

\begin{center}
\small
\setlength{\tabcolsep}{5pt}
\renewcommand{\arraystretch}{1.15}
\begin{tabular}{@{}llc@{}}
\toprule
ANSS name & Sphere stem &
\shortstack{Periodizing\\spectrum}\\
\midrule
$\beta_{(4+32k)/3}$
 & $16+192k$ & $M(1,4)$ \\
$\beta_{(8+32k)/8}+\epsilon_{30}\beta_{(6+32k)/2}$
 & $30+192k$ & $M(3,8)$ \\
$\beta_{(8+32k)/4}$
 & $38+192k$ & $M(3,8)$ \\
$\beta_{(12+32k)/4}+\epsilon_{62}\beta_{11+32k}$
 & $62+192k$ & $M(3,8)$ \\
$\beta_{(12+32k)/3}$
 & $64+192k$ & $M(3,8)$ \\
$\beta_{(16+32k)/8}+\epsilon_{78}\beta_{(14+32k)/2}$
 & $78+192k$ & $M(3,8)$ \\
$\beta_{(16+32k)/4}+\epsilon_{86}\beta_{15+32k}$
 & $86+192k$ & $M(3,8)$ \\
$\beta_{(16+32k)/3}$
 & $88+192k$ & $M(3,8)$ \\
$\beta_{(20+32k)/3}$
 & $112+192k$ & $M(3,8)$ \\
\midrule
$\beta_{(32+32k)/8}$ & $174+192k$ & $M(3,8)$ \\
$\beta_{(32+32k)/7}$ & $176+192k$ & $M(3,8)$ \\
% $\beta_{(32+32k)/6}$ & $178+192k$ & $M(3,8)$ \\
% $\beta_{(32+32k)/4}$ & $182+192k$ & $M(3,8)$ \\
% $\beta_{(32+32k)/3}$ & $184+192k$ & $M(3,8)$ \\
$\beta_{(32+32k)/2}$ & $186+192k$ & $M(3,8)$ \\
\bottomrule
\end{tabular}
\end{center}

All rows have $k\geq0$ and Adams--Novikov filtration two. For the first nine rows, the ANSS names are those of Li--Ma--Xu \cite{LiMaXu2026v2periodic}. Each $\epsilon_d\in\mathbb F_2$ is fixed throughout its row, and the displayed sum is a permanent cycle.

The first column records their ANSS names. The second column records their stems. The third column records the type-$2$ spectrum used for the periodization. These families are all $192$-periodic, constructed using a type-$2$ complex $M$ where either $M = M(1,4)$ of \cite{BHHM2008selfmap} or $M = M(3,8)$ of \cite{BMQhurewicz2021}. These four-cell complexes are known to admit $v_2$-periodic self-maps
\[v^{32}_2:\Sigma^{192}M \to M.\]

We first briefly recall the constructions of the families in the first nine rows. Let $b_0\in\pi_{2n}S^0$ denote the initial sphere class ($k=0$). Denote its lift along $q:M \to S^{10}$ (resp. $q: M \to S^{18}$) by
\[\widetilde{b_0}: S^{10+2n} \to M \qquad (\text{resp. } \widetilde{b_0}: S^{18+2n} \to M).\]
The periodization of $b_0$ is given by
\[b_k= q\circ (v_2^{32})^k\circ\Sigma^{192k}\widetilde{b_0} \in\pi_{2n+192k}S^0,\qquad k\ge0.\]
This class is nonzero by the geometric boundary theorem \cite{ravenel2004complexcobordism}.

In favorable cases, the relation $\alpha \cdot b_0 = 0$ is inherited by its periodizations, which we prove in the following lemma. 

\begin{lemma}\label{lem:anss-periodic-annihilation}
Let $M$ be either $M(1,4)$ or $M(3,8)$, and write $v^{32}_2:\Sigma^{192}M\to M$ for the respective $192$-periodic self-map. Use the notation $b_0,\widetilde{b_0},b_r$ as in the construction above.
\begin{itemize}
    \item If $\alpha\in\pi_uS^0$ satisfies $\alpha \cdot \widetilde{b_0}=0$ in $\pi_{*}M$, then $b_r\alpha=0$ for every $r\ge0$.
    \item For $b_0$ periodized via $M(1,4)$ (resp. $M(3,8)$), $2b_k =0$ (resp. $8b_k =0$) $\forall k \geq 0$.
\end{itemize}
\end{lemma}

\begin{proof}
The first part follows by construction, whereas for the second part, we factor $b_k$ as
\[S^{192k+2n+10} \xrightarrow{\Sigma^{192k}\widetilde b_0} \Sigma^{192k}M \xrightarrow{v_2^{32k}} M \to \Sigma^9C(2) \to S^{10},\]
which supplies a lift of $b_k$ in $\pi_{192k+2n+1}C(2)$; hence $2b_k =0$ by exactness.
\end{proof}

Lemma~\ref{lem:anss-periodic-annihilation} reduces the $2\nu$-annihilation required for Theorem~\ref{thm:CP 6 and 54} to a calculation on the initial lifts in $M(3,8)$. For the family lifting $\Delta^{8n}\epsilon$, we choose $\epsilon+\eta\sigma\in\pi_8S^0$ to be a representative of $\epsilon\in\pi_8\tmf$. The top-cell class $(\epsilon+\eta\sigma)[18]$ survives the homological AHSS for $M(3,8)$ (whereas $\epsilon[18]$ does not) and detects a lift in $\pi_{26}M(3,8)$.

For the family lifting $\Delta^{8n+4}\epsilon$, we use the modified Adams spectral sequence (MASS) for $M(3,8)$ \cite{behrenshillhopkinsmahowald2020detecting}. In stem $122$ and filtration $26$, we identify a permanent cycle $x$ detecting a class $a\in\pi_{122}M(3,8)$ whose $\tmf$-Hurewicz image is detected by $\Delta^4\epsilon[18]$ in the $\tmf$-homology AHSS for $M(3,8)$. To prove that $a$ has order two, we show that $h_0x$ is a $d_2$-boundary in the MASS for $M(3,8)$, and no class in stem $122$ and filtration at least $28$ survives to $E_3$.

\begin{theorem}\label{thm:epsilon lift}
\begin{enumerate}
    \item There is a periodic family $\{\epsilon_n\}$ in the stable stems $8+192n$, whose $\tmf$-Hurewicz images are $\Delta^{8n}\epsilon$. These classes satisfy 
\[\epsilon_n \cdot 2\nu =0.\]
    \item There is a periodic family $\{\mathcal{E}_n\}$ in the stable stems $104+192n$, whose $\tmf$-Hurewicz images are $\Delta^{8n+4}\epsilon$. These classes satisfy 
\[\mathcal{E}_n \cdot 2\nu =0.\]
\end{enumerate}

\end{theorem}
\begin{proof}

\begin{enumerate}
    \item Define the class $\epsilon_0 = \epsilon+\eta\sigma \in \pi_8S^0$. It has AN filtration 2 and its $\tmf$-Hurewicz image is $\epsilon \in \pi_8\tmf$. In the (homological) AHSS of $M(3,8)$, the class $\epsilon_0[18] \in \pi_{26}M(3,8)$ is a nonzero permanent cycle since 
    \[0 \in \langle \epsilon+\eta\sigma, 8, 4\rho_{15}, 8 \rangle, \]
     where $\rho_{15} \in \pi_{15}S^0$ is detected by $h_0^3h_4$ in the Adams spectral sequence. We denote a homotopy class detected by $\epsilon_0[18]$ by $\widetilde{\epsilon_0}$, which satisfies
    \[2\nu \cdot \widetilde\epsilon = 0 \in \pi_{29}M(3,8).\]
    Indeed in stem $29$, the Atiyah-Hirzebruch $E_1$-page consists of $\kappa^2[1]$ and $\{Ph_2\}[18]$, where $\{Ph_2\} \in \pi_{11}S^0_{(2)}$ is the image-of-$J$ generator detected by $Ph_2$. $\kappa^2[1]$ is a surviving cycle:
    \begin{itemize}
        \item $d_1(\kappa^2[1]) = 8\kappa^2[0] = 0$.
        \item $\pi_{12}S^0 = \pi_{13}S^0 = 0$ so $\kappa^2[1]$ cannot be hit by a differential for degree reasons.
    
    \end{itemize}
    As for $\{Ph_2\}[18]$, a combination of the May convergence theorem and the Moss convergence theorem provides $\{P^3h_2\} \in \langle \{Ph_2\},8,4\rho_{15}\rangle$; consequently,
    
    \[d_{17}(\{Ph_2\}[18]) = \{P^3h_2\}[1].\]
    This differential extends nontrivially over $2\{Ph_2\}$ and $4\{Ph_2\}$ by the Leibniz rule. Hence $\pi_{29}M(3,8) = \mathbb{Z}/2$, so the desired $2\nu$-annihilation holds for $\widetilde{\epsilon_0}$. Lemma~\ref{lem:anss-periodic-annihilation} then establishes the $2\nu$-annihilation for the entire family.
    
\item We first recall the modified Adams spectral sequences used below. The object $H(8,v_1^8)$ is obtained by taking successive cofibers of $h_0^3$ and $v_1^8$ in the derived category of $A_\ast$-comodules \cite{behrenshillhopkinsmahowald2020detecting}. The modified Adams spectral sequences for $M(3,8)$ and $\tmf\wedge M(3,8)$ take the form
\[
\begin{aligned}
E_2^{s,t}(M(3,8))
  &=\Ext^{s,t}_{A_\ast}\bigl(\mathbb F_2,H(8,v_1^8)\bigr)
    \Longrightarrow \pi_{t-s}M(3,8),\\
E_2^{s,t}(\tmf\wedge M(3,8))
  &=\Ext^{s,t}_{A(2)_\ast}\bigl(\mathbb F_2,H(8,v_1^8)\bigr)
    \Longrightarrow \tmf_{t-s}M(3,8).
\end{aligned}
\]
The unit $u\colon S^0\to\tmf$ induces a map of spectral sequences that agrees with restriction on the $E_2$-page.

To study restriction, consider the algebraic $\tmf$-resolution
\[
{}^{\mathrm{alg}}E_1^{n,r,t}
=
\Ext_{A(2)_\ast}^{r,t}
\left(
  \mathbb F_2,\,
  \operatorname{coker}\bigl(\mathbb F_2\to(A/\!/A(2))_\ast\bigr)
       ^{\otimes n}\otimes H(8,v_1^8)
\right)
\Longrightarrow
\Ext_{A_\ast}^{r+n,t}\bigl(\mathbb F_2,H(8,v_1^8)\bigr).
\]
Its edge homomorphism is restriction. Using the charts in \cite{behrenshillhopkinsmahowald2020detecting}, we take a lift
\[
x_0\in\Ext_{A_\ast}^{10,36}\bigl(\mathbb F_2,H(8,v_1^8)\bigr)
\qquad\text{of}\qquad v_1^4e_0[1].
\]
Set
\[
x=v_2^{16}x_0
\in\Ext_{A_\ast}^{26,148}\bigl(\mathbb F_2,H(8,v_1^8)\bigr),
\qquad
u_\ast x=v_2^{16}v_1^4e_0[1].
\]
The rest of the proof falls into the following four parts.

\begin{enumerate}

\item \emph{Permanence of $x$ in the modified Adams spectral sequence for $M(3,8)$.} In the modified Adams spectral sequence for $\tmf\wedge M(3,8)$, the class $v_2^{16}v_1^4e_0[1]$ detects $\Delta^4\epsilon[18]\in\tmf_{122}M(3,8)$. Thus $x$ cannot be hit by a differential, and any nonzero differential supported by $x$ would have to map to zero under $u_\ast$. The proof of Theorem~11.1 of \cite{behrenshillhopkinsmahowald2020detecting} rules out all such targets. We may therefore choose a class $a\in\pi_{122}M(3,8)$ detected by $x$ whose projection to the top cell has $\tmf$-Hurewicz image $\Delta^4\epsilon$.

\item \emph{The differential $d_2(v_2^{16}j[1])=h_0x$.} The class $h_0x$ is nonzero on the $E_2$-page because its image $h_0v_2^{16}v_1^4e_0[1]$ is nonzero. The class $h_0x_0$ has algebraic Atiyah--Hirzebruch spectral sequence name $h_0Pe_0[1]$. In the classical Adams spectral sequence, $d_2(j)=h_0Pe_0$ \cite{isaksen2019stable}. Since $h_0^3j=0$, the class $j$ lifts to a class $j[1]$ in the modified Adams spectral sequence for the mod-$8$ Moore spectrum $C(8)$. Naturality under the projection $C(8)\to S^1$ gives
\[
d_2(j[1])=h_0Pe_0[1].
\]
Mapping to the modified Adams spectral sequence for $M(3,8)$ and applying the Leibniz rule, we obtain
\[
d_2(v_2^{16}j[1])=h_0x\ne0,
\]
since $d_2(v_2^{16})=0$. Consequently, the homotopy class $2a$ lies in filtration at least $28$.

\item \emph{Injectivity of restriction in filtration at least $28$.} In stem $122$ and total filtration $s\geq28$, the Brown--Gitler decomposition and the charts in \cite{behrenshillhopkinsmahowald2020detecting} give
\[
{}^{\mathrm{alg}}E_1^{n,s-n,122+s}
=
\Ext_{A(2)_\ast}^{s-n,122+s}
\left(
  \mathbb F_2,\,
  \operatorname{coker}\bigl(\mathbb F_2\to(A/\!/A(2))_\ast\bigr)
       ^{\otimes n}\otimes H(8,v_1^8)
\right)
=0
\qquad(n\geq1).
\]
Thus only the length-zero column contributes, and restriction is injective in bidegree $(s,122+s)$ for every $s\geq28$.

\item \emph{The three remaining classes are $d_2$-boundaries.} In stem $122$, the $A(2)_\ast$-chart vanishes above filtration $28$, and the group in filtration $28$ has basis
\[
h_0^2u_\ast x,\qquad v_1^6h_{2,1}^{22}[0],\qquad
v_1^3h_{2,1}^{23}[1].
\]
Here we write the latter two classes by their names in the $g$-inverted modified Adams spectral sequence, where $g=h_{2,1}^4$ \cite{BMQhurewicz2021}. In the modified Adams spectral sequence for $M(3,8)$, the Leibniz rule gives
\[
d_2(h_0v_2^{16}j[1])=h_0^2x.
\]
For the other two classes, we use the following differentials in the $g$-inverted modified Adams spectral sequence for $\tmf\wedge M(3,8)$:
\[
\begin{aligned}
d_2(v_2^8v_1^5h_{2,1}^{13}[0]) &= v_1^6h_{2,1}^{22}[0],\\
d_2(v_2^8v_1^2h_{2,1}^{14}[1]) &= v_1^3h_{2,1}^{23}[1].
\end{aligned}
\]
The wedge classes constructed in \cite{BMQhurewicz2021} give lifts of both sources to the $E_2$-page of the unlocalized modified Adams spectral sequence for $\tmf\wedge M(3,8)$. In bidegree $(28,150)$, the two wedge classes remain linearly independent after $g$-localization, while $h_0^2u_\ast x$ maps to zero. Hence by naturality, the unlocalized differentials for $\tmf\wedge M(3,8)$ have the form
\[
\begin{aligned}
d_2(v_2^8v_1^5h_{2,1}^{5}g^2[0])
  &= v_1^6h_{2,1}^{14}g^2[0]+\lambda_0h_0^2 u_*x,\\
d_2(v_2^8v_1^2h_{2,1}^{6}g^2[1])
  &= v_1^3h_{2,1}^{15}g^2[1]+\lambda_1h_0^2u_*x,
\end{aligned}
\]
where $\lambda_0,\lambda_1\in\mathbb F_2$. The wedge construction further lifts both the sources and the targets to the modified Adams spectral sequence for $M(3,8)$. Since restriction is injective in bidegree $(28,150)$, naturality gives the same formulas for these lifts. It follows that
\[
E_3^{s,s+122}(M(3,8))=0
\qquad\text{for every }s\geq28.
\]

\end{enumerate}

Consequently, there exists neither an explicit $h_0$-extension nor a hidden extension on $x\in\Ext^{26,148}_{A_\ast}\bigl(\mathbb F_2,H(8,v_1^8)\bigr)$. Hence $2a=0$ in homotopy. According to Lemma~\ref{lem:anss-periodic-annihilation}, the periodization
\[
q \circ (v^{32}_2)^k \circ \Sigma^{192k}a
\]
also has order $2$. Hence the family is annihilated by $2\nu$.

\end{enumerate} 
\end{proof}

The next lemma supplies $\sigma, 2\sigma$ and $8\sigma$-attaching maps in certain subquotients $\CP^{n+3}_n$, which lays the foundation for further Atiyah-Hirzebruch analysis.

\begin{lemma}\label{lem:anss-cell-attachments}

\begin{enumerate}
\item\label{item:anss-cell-bottom-split} If $8\mid m$, then $\CP_m^{m+3} \cong S^{2m} \vee \CP^{m+3}_{m+1}$. Denote the quotient map by
\[
 q:\Sigma^{-2m}\CP_m^{m+3}\longrightarrow S^0.
\]
If $m\equiv8\pmod{16}$ and $g:S^7\to\Sigma^{-2m}\CP_m^{m+3}$ is the successive attaching map, we have
\[
 q \circ g=2\sigma.
\]
\item\label{item:anss-cell-32} If $32\mid m$, then, with the notation of \textup{(\ref{item:anss-cell-bottom-split})},
\[
q\circ g=
\begin{cases}
8\sigma,&\nu_2(m)=5,\\
0,&\nu_2(m)\geq6.
\end{cases}
\]
    \item\label{item:anss-cell-mod8-five} If $m\equiv5\pmod8$, then
\[Sq^8(z^m) = z^{m+4}, \qquad Sq^8(z^{m+1}) = z^{m+5}\]
in $H^\ast(\CP^{m+5}_m;\mathbb F_2)$. 
\item If $m\equiv5\pmod8$, then $\CP^{m+5}_m$ has a $4$-cell $H\mathbb{F}_2$-subquotient which we denote by $X$. The $\mathbb{F}_2$-cohomology of $X$ is given by
\[H^n(X) = \mathbb{F}_2 \qquad \text{for } n=2m,2m+2,2m+8,2m+10.\]
Denote the cohomology generators by $e^i$ where $i=2m,2m+2,2m+8,2m+10$. 
\[Sq^2(e^{2m}) = e^{2m+2}, Sq^2(e^{2m+8}) = e^{2m+10}, Sq^8(e^{2m}) = e^{2m+8}, Sq^8(e^{2m+2}) = e^{2m+10}.\]

\end{enumerate}

\end{lemma}
\begin{proof}
    The first three parts follow from \cite{MosherProjective1968}. For the last part, $S^{2m+4}$ is an $H\mathbb{F}_2$-subcomplex of $\CP^{m+5}_m$. Take the cofiber of the inclusion
    \[S^{2m+4} \to \CP^{m+5}_m \to Y\]
    and denote it by $Y$. $S^{2m+6}$ is an $H\mathbb{F}_2$-quotient of $\CP^{m+5}_m$. Take the fiber of the quotient
    \[X \to Y \to S^{2m+6}\]
    and denote it by $X$. The Steenrod squares act on $H^\ast(X)$ as desired.
\end{proof}

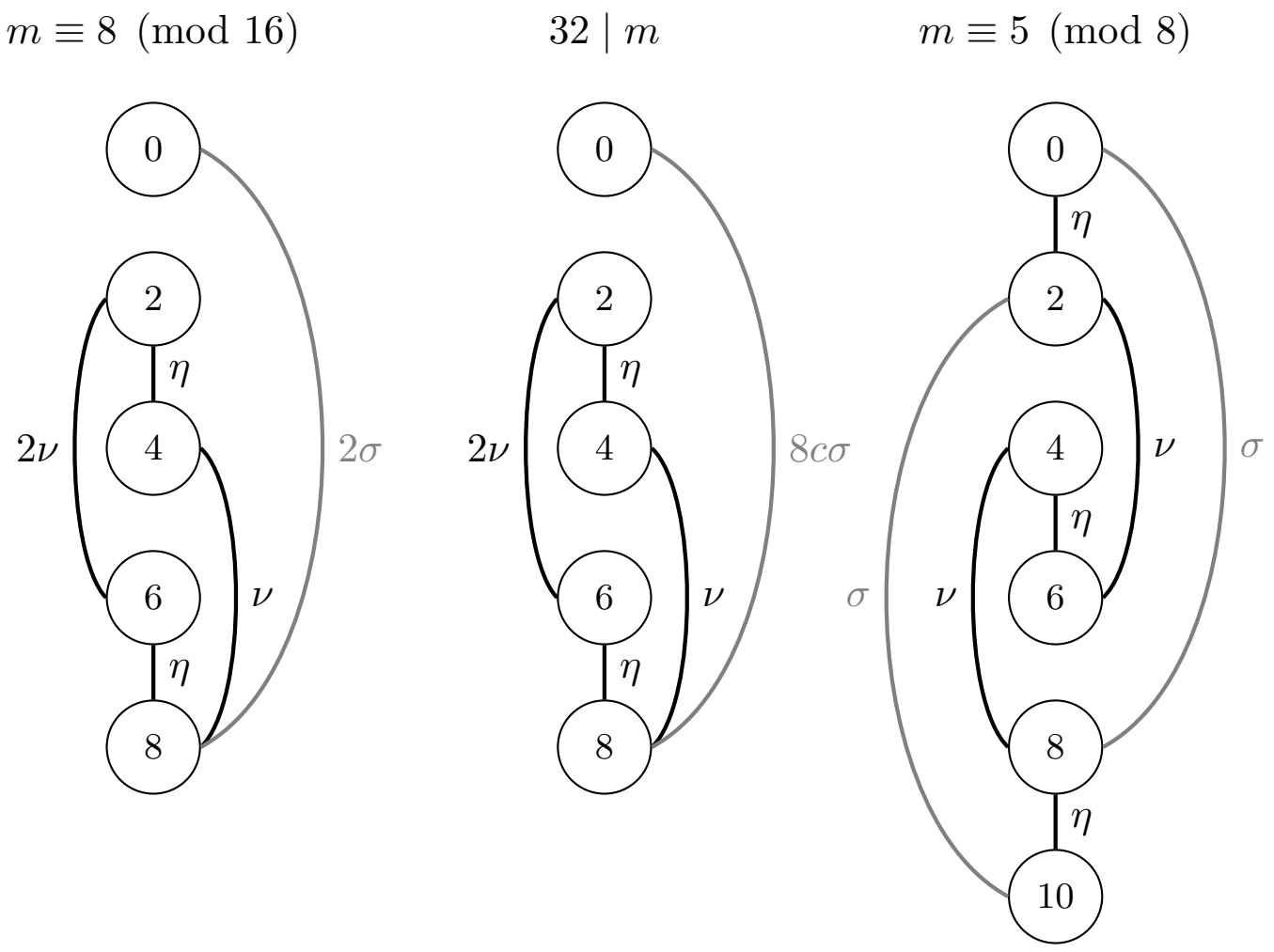
\begin{figure}[htbp]
\centering
\begingroup
\tikzset{
  cell/.style={circle,draw,fill=white,minimum size=8mm, inner sep=.4pt,font=\small}, eta/.style={black,line width=.8pt}, nuedge/.style={black,line width=.8pt}, twonuedge/.style={black,line width=.8pt}, sigmaedge/.style={gray,line width=.8pt}
}
\begin{tikzpicture}[x=.9cm,y=.95cm,scale=.82,transform shape]

% Left panel
\begin{scope}
  \foreach \d/\y in {0/0,2/-1.35,4/-2.7,6/-4.05,8/-5.4} {
    \node[cell] (a\d) at (0,\y) {$\d$};
  }
  \draw[eta] (a2) -- node[right] {$\eta$} (a4);
  \draw[eta] (a6) -- node[right] {$\eta$} (a8);
  \draw[twonuedge] (a2.west)
    .. controls (-.85,-1.75) and (-.85,-3.65) ..
    node[midway,left] {$2\nu$} (a6.west);
  \draw[nuedge] (a4.east)
    .. controls (.9,-3.1) and (.9,-5) ..
    node[midway,right] {$\nu$} (a8.east);
  \draw[sigmaedge] (a0.east)
    .. controls (2,-.8) and (2,-4.6) ..
    node[midway,right] {$2\sigma$} (a8.east);
  \node at (0,1.05) {$m\equiv8\pmod{16}$};
\end{scope}

% Middle panel
\begin{scope}[shift={(4.3,0)}]
  \foreach \d/\y in {0/0,2/-1.35,4/-2.7,6/-4.05,8/-5.4} {
    \node[cell] (b\d) at (0,\y) {$\d$};
  }
  \draw[eta] (b2) -- node[right] {$\eta$} (b4);
  \draw[eta] (b6) -- node[right] {$\eta$} (b8);
  \draw[twonuedge] (b2.west)
    .. controls (-.85,-1.75) and (-.85,-3.65) ..
    node[midway,left] {$2\nu$} (b6.west);
  \draw[nuedge] (b4.east)
    .. controls (.9,-3.1) and (.9,-5) ..
    node[midway,right] {$\nu$} (b8.east);
  \draw[sigmaedge] (b0.east)
    .. controls (2,-.8) and (2,-4.6) ..
    node[midway,right] {$8c\sigma$} (b8.east);
  \node at (0,1.05) {$32\mid m$};
\end{scope}

% Right panel
\begin{scope}[shift={(8.6,0)}]
  \foreach \d/\y in {0/0,2/-1.35,4/-2.7,6/-4.05,8/-5.4,10/-6.75} {
    \node[cell] (c\d) at (0,\y) {$\d$};
  }
  \draw[eta] (c0) -- node[right] {$\eta$} (c2);
  \draw[eta] (c4) -- node[right] {$\eta$} (c6);
  \draw[eta] (c8) -- node[right] {$\eta$} (c10);
  \draw[nuedge] (c2.east)
    .. controls (.9,-1.75) and (.9,-3.65) ..
    node[midway,right] {$\nu$} (c6.east);
  \draw[nuedge] (c4.west)
    .. controls (-.9,-3.1) and (-.9,-5) ..
    node[midway,left] {$\nu$} (c8.west);
  \draw[sigmaedge] (c0.east)
    .. controls (2,-.8) and (2,-4.6) ..
    node[midway,right] {$\sigma$} (c8.east);
  \draw[sigmaedge] (c2.west)
    .. controls (-2,-2.15) and (-2,-5.95) ..
    node[midway,left] {$\sigma$} (c10.west);
  \node at (0,1.05) {$m\equiv5\pmod8$};
\end{scope}

\end{tikzpicture}
\par\endgroup
\caption{Stable $2$-local cell diagrams desuspended by $2m$. The short edges and arcs record the $\eta,\nu$ and $2\nu$-attaching maps, whereas the long grey arcs record $\sigma,2\sigma$ and $8c\sigma$-attaching maps (where $c$ is either $0$ or $1$).}
\label{fig:anss-five-cell-sigma-attachments}
\end{figure}

Now, according to Proposition~\ref{prop:anss-top-cell-detection}, the classes $b_k[2n+192k]\in[\CP^{n+96k},S^0]$ are nonzero, where $b_k$ is listed in the preceding table. We first apply the $b$-families in the first nine rows, extending the classes $b_k[2n+192k]$ to $\CP^{n+96k+r}$.

\begin{theorem}\label{thm:anss-dimensions-8-12}
For every integer $k\ge0$, the spaces $\CP^{m+96k}$ admit exotic smooth structures for $8\le m\le12$.
\end{theorem}

\begin{proof}

Let $b_k$ be the class detected by $\beta_{(4+32k)/3}$. Proposition~\ref{prop:anss-top-cell-detection} shows that $0 \neq b_k[16+192k] \in [\CP^{8+96k},S^0]$. Lemma~\ref{lem:anss-periodic-annihilation} gives $2b_k=0$.

By Lemma~\ref{lem:anss-cell-attachments}, the first potentially nonzero differential on $b_k[16+192k]$ in the AHSS of $\CP^{\infty}_{8+96k}$ is a $d_8$ differential, whereas
\[d_8(b_k[16+192k]) = 2\sigma \cdot b_k[24+192k] = 0.\]
Hence $b_k[16+192k]$ extends to $[\CP^{12+96k},S^0]$ and Theorem~\ref{AHSSsmoothing} supplies the exotic complex projective spaces.

\end{proof}

\begin{theorem}\label{thm:anss-dimensions-15}
For every integer $k\geq0$, the space $\CP^{15+96k}$ admits an exotic smooth structure.
\end{theorem}
\begin{proof}

Let $b_k$ be the class detected by $\beta_{(8+32k)/8}+\epsilon_{30}\beta_{(6+32k)/2}$. Proposition~\ref{prop:anss-top-cell-detection} shows that $0 \neq b_k[30+192k] \in [\CP^{15+96k},S^0]$. For $k\geq1$, $15+96k$ is greater than $15$ and is neither $31$ nor $63$ modulo $96$. Theorem~\ref{AHSSsmoothing} proves the claim. For $k=0$, Theorem~\ref{thm:anss-dimensions-8-12} shows that $\eta_4[16]$ survives nontrivially in the AHSS of $\CP^{12}$. Any nontrivial differential supported by this class in the AHSS of $\CP^{15}$ would therefore have the form
\[
d_{2r}(\eta_4[16])=\alpha[16+2r],
\qquad
\alpha\in\pi_{16+2r-1}S^0,
\qquad
5\le r\le7.
\]
Here $\eta_4 \in \pi_{16}S^0$ is detected by $h_1h_4$. All such candidates $\alpha$ lie in the image of $J$. Consequently, $\eta_4[16]$ survives nontrivially in the AHSS of $\CP^{15}$ and detects a nonzero class in $[\CP^{15},S^0]$. Combining Proposition~\ref{LES} and Theorem~\ref{map of cohomology} and arguing as in the proof of Theorem~\ref{AHSSsmoothing}, we conclude that this class determines an exotic smooth structure on $\CP^{15}$.
\end{proof}

\begin{theorem}\label{thm:anss-dimensions-19}
\label{thm:anss-dimensions-19-23}
For every integer $k\geq0$, the space $\CP^{19+96k}$ admits an exotic smooth structure.
\end{theorem}
\begin{proof}

Let $b_k$ be the class detected by $\beta_{(8+32k)/4}$. Proposition~\ref{prop:anss-top-cell-detection} shows that $0 \neq b_k[38+192k] \in [\CP^{19+96k},S^0]$. Theorem~\ref{AHSSsmoothing} proves the claim.
\end{proof}

\begin{theorem}\label{thm:anss-dimensions-31}
For every integer $k \geq 0$, the space $\CP^{31+96k}$ admits an exotic smooth structure.
\end{theorem}

\begin{proof}
For $k \geq 1$, let $b_k\in\pi_{62+192k}S^0$ be detected by $\beta_{(12+32k)/4}+\epsilon_{62}\beta_{11+32k}$. Proposition~\ref{prop:anss-top-cell-detection} gives
\[
0\ne b_k[62+192k]\in[\CP^{31+96k},S^0].
\]
The nonexistence theorem of Hill--Hopkins--Ravenel \cite{hillhopkinsravenel2016nonexistence} implies that the Kervaire-invariant homomorphism
\[
i_{62+192k}:\pi_{62+192k}S^0\cong\pi_{62+192k}SG
\longrightarrow\pi_{62+192k}(G/Top)\cong\mathbb Z/2
\]
is zero. Naturality with respect to the top-cell collapse, together with Remark~\ref{translation preserves filtration}, gives
\[
0\ne {t}_*(b_k[62+192k])\in
\ker\bigl(i_*:[\CP^n,SG]\longrightarrow[\CP^n,G/Top]\bigr).
\]
The claim follows from Proposition~\ref{LES} and Theorem~\ref{concordanceCP}. For $k=0$, choose a class $a\in\pi_{58}S^0$ detected by $h_1Q_2$ such that $\eta a=0$. Such a choice exists by \cite{isaksenwangxu2023stable}. The product $a\sigma\in\pi_{65}S^0$ is nonzero and is detected by $h_1h_3Q_2$.

Consider the $2$-local AHSS for $\CP^\infty_{29}$. Since $\eta a=0$ and $(\pi_{61}S^0)_{(2)}=0$, the class $a[58]$ survives to the $E_6$-page. If $d_6(a[58])\neq0$, then Theorem~\ref{finiteinfinitecriterion} already implies that $a[58]$ detects a nonzero class in $[\CP^{29},S^0]$. Otherwise, $a[58]$ survives to the $E_8$-page. On this page, $(a\sigma)[66]$ is nonzero, and by Lemma~\ref{lem:anss-cell-attachments}
\[
d_8(a[58])=(a\sigma)[66]\neq0.
\]
Thus, in either case, Theorem~\ref{finiteinfinitecriterion} shows that $a[58]$ detects a nonzero class in $[\CP^{29},S^0]$, which extends over $[\CP^{30},S^0]$. Applying Theorem~\ref{AHSSsmoothing}, we conclude that $\CP^{31}$ admits an exotic smooth structure.

\end{proof}

\begin{theorem}\label{thm:anss-dimensions-32-36}
\label{thm:anss-dimensions-32-35}
For every integer $k\geq0$, the spaces $\CP^{m+96k}$ admit exotic smooth structures for $32\leq m\leq36$.
\end{theorem}

\begin{proof}
Let $b_k$ be the class detected by $\beta_{(12+32k)/3}$. Proposition~\ref{prop:anss-top-cell-detection} gives
\[
0\ne b_k[64+192k]
\in[\CP^{32+96k},S^0].
\]
Lemma~\ref{lem:anss-periodic-annihilation} gives $8b_k=0$.

Since $32\mid32+96k$, according to Lemma~\ref{lem:anss-cell-attachments}\textup{(\ref{item:anss-cell-32})}, the first potentially nonzero differential on $b_k[64+192k]$ is
\[
d_8\bigl(b_k[64+192k]\bigr)
=(8c\sigma)b_k[72+192k]=0.
\]
The last equality holds by Lemma~\ref{lem:anss-periodic-annihilation}. Hence $b_k[64+192k]$ extends to $[\CP^{36+96k},S^0]$. Theorem~\ref{AHSSsmoothing} thus supplies the exotic complex projective spaces.
\end{proof}

\begin{theorem}\label{thm:anss-dimensions-39}
\label{thm:anss-dimensions-39-40}
For every integer $k\geq0$, the spaces $\CP^{39+96k}$ and $\CP^{40+96k}$ admit exotic smooth structures.
\end{theorem}
\begin{proof}
Let $b_k$ be the sphere class detected by $\beta_{(16+32k)/8}+\epsilon_{78}\beta_{(14+32k)/2}$. Proposition~\ref{prop:anss-top-cell-detection} gives $0\ne b_k[78+192k]\in[\CP^{39+96k},S^0]$. By \cite{LiMaXu2026v2periodic}, the classes $b_k$ are $2$-divisible, so that
\[d_2(b_k[78+192k]) = \eta b_k[80+192k] = 0.\]
Theorem~\ref{AHSSsmoothing} proves both claims.
\end{proof}

\begin{theorem}\label{thm:anss-dimensions-43}
For every integer $k\geq0$, the space $\CP^{43+96k}$ admits an exotic smooth structure.
\end{theorem}
\begin{proof}
Let $b_k$ be the sphere class detected by $\beta_{(16+32k)/4}+\epsilon_{86}\beta_{15+32k}$. Proposition~\ref{prop:anss-top-cell-detection} gives
\[
 0\ne b_k[86+192k]\in[\CP^{43+96k},S^0].
\]
Theorem~\ref{AHSSsmoothing} proves the claim.
\end{proof}

\begin{theorem}\label{thm:anss-dimensions-44-47}
For every integer $k\geq0$, the spaces $\CP^{m+96k}$ admit exotic smooth structures for $44\leq m\leq47$.
\end{theorem}
\begin{proof}
Put $b_k=\beta_{(16+32k)/3}$. Proposition~\ref{prop:anss-top-cell-detection} gives $0\ne b_k[88+192k]\in[\CP^{44+96k},S^0]$. By \cite{LiMaXu2026v2periodic}, $b_k$ can be chosen so that $b_k(2\nu)=0$ for every $k$.

Put $n=44+96k$. Since $n\equiv4\pmod8$, we have 
\[
 \Sigma^{-88-192k}\CP_{44+96k}^{47+96k}\simeq
 \operatorname{Cof}\bigl(S^3\xrightarrow{(2\nu,\eta)}
                         S^0\vee S^2\bigr)\vee S^6
\]
by \cite{MosherProjective1968}. Consequently, 
\[
 d_4(b_k[88+192k])= b_k(2\nu)[92+192k]=0 .
\]
The splitting of the top cell in $\CP^{47+96k}_{44+96k}$ shows $b_k[88+192k]$ further extends to $[\CP^{47+96k},S^0]$, and Theorem~\ref{AHSSsmoothing} gives all four exotic smooth structures.
\end{proof}

\begin{theorem}\label{thm:anss-dimensions-56-60}
For every integer $k\geq0$, the spaces $\CP^{m+96k}$ admit exotic smooth structures for $56\leq m\leq60$.
\end{theorem}
\begin{proof}
Put $b_k=\beta_{(20+32k)/3}$. Proposition~\ref{prop:anss-top-cell-detection} gives $0\ne b_k[112+192k]\in[\CP^{56+96k},S^0]$.

Let $\widetilde{b_0}$ be the homotopy class in $\pi_{130}M(3,8)$ detected by $b_0[18]$. Once we prove $2\widetilde{b_0} = 0$, by Lemma~\ref{lem:anss-periodic-annihilation}, $2b_k=0$ for every $k\geq0$. In $\pi_{112}S^0$, $2b_0 =0$. The relations
\[\langle 2,b_0,8 \rangle \supseteq  \langle 2,b_0,2 \rangle 4\ni 4\eta b_0 = 0\]
exclude the hidden extension from dimension $18$ to $17$. $b_0$ is detected by the class $x_{112,18}$ in $\Ext_A^{18,130}(\mathbb F_2,\mathbb F_2)$. A Massey product computation shows that
\[0 = Ph_1c_0 \cdot x_{112,18} \in \langle h_0,x_{112,18},h_0^3,h_0^5h_4 \rangle \]
in $\Ext^{26,155}_A(\mathbb F_2,\mathbb F_2)$. All the classes in stem $129$ that have filtration higher than $29$ are $K(1)$-local. Hence the hidden extension from filtration $18$ to $1$ is also excluded. By a similar argument one can show that there is no hidden extension from $b_0[18]$ to classes with AHSS filtration $0$. We thus conclude that 2$\widetilde{b_0} = 0$.

Put $n=56+96k$. Since $n\equiv8\pmod{16}$, Lemma~\ref{lem:anss-cell-attachments} shows that
\[
 d_8(b_k[112+192k])=(2\sigma)b_k[120+192k]=0 .
\]
Hence $b_k[112+192k]$ extends to $[\CP^{60+96k},S^0]$. Theorem~\ref{AHSSsmoothing} gives all five exotic smooth structures.
\end{proof}

% BLOCK 2: Replace the paragraph beginning
% "For the following three results, put M=M(3,8)" up to Theorem 7.17.
% The y_6 and b_{r,2} setup belongs to the separately extracted 93--98 proof.

We now recall the construction of the families in the last three rows of the preceding chart, namely $\beta_{32r/j}$ for $r\geq1$ and $j=8,7,2$. Unlike those of the preceding families, the initial classes in $\pi_{16-2j}M(3,8)$ are detected below the top AH filtration.

Let $\{Ph_2\}\in\pi_{11}S^0_{(2)}$ be an image-of-$J$ generator detected by $Ph_2$. The following table records the AHSS names and $BP$-Hurewicz images of the nonzero classes $e_j\in\pi_{16-2j}M(3,8)$:
\[
\begin{array}{c|c|c}
j & \text{AHSS name of }e_j & h_{BP}(e_j) \\ \hline
8 & 4[0]                    & 4          \\
7 & \eta[1]                 & 4v_1       \\
2 & \bigl(2\{Ph_2\}\bigr)[1] & 4v_1^6
\end{array}
\]

For $r\geq1$, define
\[
b_{r,j}
=\Sigma^{-18}\bigl(q\circ v_2^{32r}\circ\Sigma^{192r}e_j\bigr)
\in\pi_{192r-2j-2}S^0.
\]
For each of these families, we have
\[
h_{BP}\bigl(v_2^{32r}\circ\Sigma^{192r}e_j\bigr)
=4v_1^{8-j}v_2^{32r}.
\]
The geometric boundary theorem then shows that Ravenel's classes $\beta_{32r/j}$ are permanent cycles detecting $b_{r,j}$ in ANSS filtration two \cite{ravenel2004complexcobordism}.

In particular, the classes $\beta_{32r/2}$ are $2$-divisible: by definition on the ANSS $E_2$-page we have
\[2\beta_{32r/2,2} = \beta_{32r/2}.\]
The classes $\beta_{32r/2,2}$ are the divided $\beta$-elements, which detect the periodization of $z \in \pi_{12}M$ in the ANSS of the sphere. We denote this family by $b_{r/2,2}$.

\begin{theorem}\label{thm:anss-dimensions-87-88}
For every integer $k\geq0$, the spaces $\CP^{m+96k}$ admit exotic smooth structures for $87\leq m\leq88$.
\end{theorem}
\begin{proof}
Put $r=k+1$. Proposition~\ref{prop:anss-top-cell-detection} gives
\[
0\ne b_{r,8}[192r-18]\in[\CP^{96r-9},S^0].
\]
Since $4\eta=0$, $\eta \cdot e_8 = 0 \in \pi_1M(3,8)$. Lemma~\ref{lem:anss-periodic-annihilation} therefore gives $b_{r,8}\eta=0$.

Since $96r-9$ is odd, the attaching map between the first two cells is $\eta$, as detected by $Sq^2(z^m)=m z^{m+1}$. Hence 
\[
d_2\bigl(b_{r,8}[192r-18]\bigr)
=\eta b_{r,8}[192r-16]=0.
\]
Hence the class extends nontrivially to $[\CP^{96r-8},S^0]$, and Theorem~\ref{AHSSsmoothing} proves both claims.
\end{proof}

% BLOCK 4: Replace the current Theorem 7.18 (93--98) and its proof
% by the missing 88--92 result. The 93--98 argument is kept separately.
% Do NOT retain the label thm:anss-dimensions-93-98 on this theorem.
\begin{theorem}\label{thm:anss-dimensions-88-92}
For every integer $k\geq0$, the spaces $\CP^{m+96k}$ admit exotic smooth structures for $88\leq m\leq92$.
\end{theorem}
\begin{proof}
Put $r=k+1$. Proposition~\ref{prop:anss-top-cell-detection} gives
\[
0\ne b_{r,7}[192r-16]\in[\CP^{96r-8},S^0].
\]

Since $96r-8\equiv8\pmod{16}$, according to Lemma~\ref{lem:anss-cell-attachments}, the first potentially nonzero differential on $b_{r,7}[192r-16]$ is
\[
d_8\bigl(b_{r,7}[192r-16]\bigr)
=(2\sigma)b_{r,7}[192r-8]=0.
\]
The last equality follows from the fact that $2b_{r,7}=0$, which is implied by Lemma~\ref{lem:anss-periodic-annihilation}. Hence the class extends to $[\CP^{96r-4},S^0]$ and restricts nontrivially to every intermediate projective space. Theorem~\ref{AHSSsmoothing} proves all five claims.
\end{proof}

\begin{theorem}\label{thm:anss-dimensions-93-98}
For every integer $k\geq0$, the spaces $\CP^{m+96k}$ admit exotic smooth structures for $93\leq m \leq 98$.
\end{theorem}

\begin{proof}
Put $r=k+1$. Proposition~\ref{prop:anss-top-cell-detection} gives
\[
0\ne b_{r,2}[192r-6] = 2b_{r/2,2}[192r-6] \in[\CP^{96r-3},S^0].
\]
\begin{enumerate}
    \item The $d_2$-differential on $b_{r,2}[192r-6]$ is zero since it is $2$-divisible.
    \item The $d_4$-differential on $b_{r,2}[192r-6]$ is zero since the top cell splits off in $\CP_{96r-3}^{96r-1}$.
    \item The $d_6$-differential on $b_{r,2}[192r-6]$ is zero: the value of the differential lies in
    \[\langle 2b_{r/2,2}, \eta, \nu\rangle[192r]. \]
    Shuffling provides $\langle 2b_{r/2,2}, \eta, \nu\rangle = \{0\} + \nu \pi_{192r-4}$.
    \item The $d_8$-differential on $b_{r,2}[192r-6]$ is zero since $0 = \sigma\cdot \{Ph_2\}[1] \in \pi_{19}M$.
    \item The $d_{10}$-differential on $b_{r,2}[192r-6]$ is zero. Besides $\sigma\cdot \{Ph_2\}[1] = 0$, low-dimensional computations also give $\alpha\cdot\{Ph_2\}[1] = 0$ where $\alpha$ can be any element of $\pi_9S^0$. Lemma~\ref{lem:anss-periodic-annihilation} therefore gives $b_{r/2,2}\sigma=0$ and $b_{r/2,2}\pi_9S^0=0$. Consequently,
\[\left\langle2b_{r/2,2},(\eta,\sigma),
       \begin{pmatrix}\sigma\\\eta\end{pmatrix}\right\rangle \qquad \text{and} \qquad \left\langle b_{r/2,2},(2\eta,2\sigma),
       \begin{pmatrix}\sigma\\\eta\end{pmatrix}\right\rangle \]
have the same indeterminacy
\[
 \pi_{192r-4}S^0\cdot\sigma+\pi_{192r+2}S^0\cdot\eta.
\]
By shuffling,
\[\left\langle2b_{r/2,2},(\eta,\sigma),
       \begin{pmatrix}\sigma\\\eta\end{pmatrix}\right\rangle = \left\langle b_{r/2,2},(2\eta,2\sigma),
       \begin{pmatrix}\sigma\\\eta\end{pmatrix}\right\rangle \supseteq \langle  b_{r/2,2},2\sigma,\eta\rangle \supseteq \langle  b_{r/2,2}\sigma,2,\eta\rangle \ni 0.\]
It follows that 
\[d_{10}(b_{r,2}[192r-6]) = 0\]
in the four cell complex $X$ of Lemma~\ref{lem:anss-cell-attachments}. It follows by naturality that $d_{10}(b_{r,2}[192r-6]) = 0$ also holds in $\CP^{96r+2}_{96r-3}$.

\end{enumerate}
 
\end{proof}

\begin{proof}[Proof of Theorem~\ref{thm:anss-main}]
For $k\geq0$, combine Theorems~\ref{thm:anss-dimensions-8-12}, \ref{thm:anss-dimensions-15} \ref{thm:anss-dimensions-19}, \ref{thm:anss-dimensions-31} \ref{thm:anss-dimensions-32-36}, \ref{thm:anss-dimensions-39-40}, \ref{thm:anss-dimensions-43}, \ref{thm:anss-dimensions-44-47}, \ref{thm:anss-dimensions-56-60}, \ref{thm:anss-dimensions-87-88}, \ref{thm:anss-dimensions-88-92} and \ref{thm:anss-dimensions-93-98}.
\end{proof}

\section{Three-primary detection}\label{sec:three-primary-detection}
In this section, we work $3$-locally.

\begin{theorem}\label{thm:anss3-main}
For every integer $k\geq0$, the space $\CP^{n+72k}$ admits an exotic smooth structure whenever
$$
n\in [5,6]\cup[13,14]\cup[19,24]\cup[37,38]\cup[41,50]\cup[53,72].
$$
\end{theorem}

We apply the Adams--Novikov method of Section~\ref{sec:anss-detection} to the classes listed in Table~\ref{tab:anss3-detectors}. Proposition~\ref{prop:anss3-top-cell} establishes nonzero top-cell pullbacks, and Corollary~\ref{cor:anss3-extension} extends to the intervals.

\begin{table}[htbp]
\centering
\small
\setlength{\tabcolsep}{2.5pt}
\renewcommand{\arraystretch}{1.2}
\caption{Three-primary detectors in sphere stem $2m$. The column $s$ specifies the interval $[m,m+s]$
for which exotic smooth structures exist.
%supplied by Corollary~\ref{cor:anss3-extension}.
}
\label{tab:anss3-detectors}

\begin{tabular*}{\linewidth}{@{\extracolsep{\fill}}
  lcc>{\raggedright\arraybackslash}p{0.19\linewidth}
  @{\hspace{1em}}
  lcc>{\raggedright\arraybackslash}p{0.19\linewidth}@{}}
\toprule
Detector & $m\bmod72$ & $s$ & Reference
& Detector & $m\bmod72$ & $s$ & Reference \\
\midrule

$\beta_{9t+1}$ & $5$ & $1$
& \multirow{2}{=}{\raggedright
    \cite{behrens2004v2selfmap}}
& $\beta_{9t+6}$ & $45$ & $5$
& \cite{belmont2023beta} \\
\cmidrule[0.3pt](lr){5-8}
$\beta_{9t+2}$ & $13$ & $1$ &
& $\Xi_{t+1}$ & $53$ & $1$
& Proposition~\ref{prop:anss3-corrected-family} \\
\cmidrule[0.3pt](lr){1-4}
\cmidrule[0.3pt](lr){5-8}

$\beta_{(9t+3)/2}$ & $19$ & $1$
& \multirow{2}{=}{\raggedright
    \cite{belmont2023beta}}
& $\beta_{(9t+9)/8}$ & $55$ & $1$
& \multirow{5}{=}{\raggedright
    \cite{belmont2023beta}} \\
$\beta_{9t+3}$ & $21$ & $3$ &
& $\beta_{(9t+9)/7}$ & $57$ & $3$ & \\
\cmidrule[0.3pt](lr){1-4}

$\beta_{9t+5}$ & $37$ & $1$
& \cite{behrens2004v2selfmap}
& $\beta_{(9t+9)/5}$ & $61$ & $1$ & \\
\cmidrule[0.3pt](lr){1-4}

$\beta_{(9t+6)/3}$ & $41$ & $1$
& \multirow{2}{=}{\raggedright
    \cite{belmont2023beta}}
& $\beta_{(9t+9)/4}$ & $63$ & $5$ & \\
$\beta_{(9t+6)/2}$ & $43$ & $1$ &
& $\beta_{9t+9}$ & $69$ & $3$ & \\
\bottomrule
\end{tabular*}
\end{table}

The degrees are
\begin{equation*}
|\beta_{I/j}|=16I-4j-2=2m,\qquad m=8I-2j-1.
\end{equation*}

The following proposition holds at every odd prime and will also be used in Section~\ref{sec:density-one}.

\begin{prop}\label{prop:anss3-top-cell}
Let $p$ be an odd prime.

1)~~For every $m\geq 1$, the top-cell collapse $q:\CP^m\to S^{2m}$ induces an injection
$$
q^*:\Ext_{BP_*BP}^{2,2m+2}(BP_*,BP_*)
\longrightarrow \Ext_{BP_*BP}^{2,2}(BP_*\CP^m,BP_*).
$$

2)~~Suppose that the ANSS class $b\in \Ext^{2,2m+2}$ detects $\beta\in\pi_{2m}S^0$. Then the class $q^*b$ detects the nonzero pullback $q^*\beta$ in the ANSS computing $[\CP^m,S^0]$.
\end{prop}
\begin{proof}
\textup{(1)} As in the proof of Theorem~\ref{thm:anss-top-edge}, a differential hitting $b[2m]$ in the algAHSS must start on $y[2\ell]$ for some $1\leq\ell<m$. Put $w=\ell+1$. The coefficient $y$ lies on the sphere ANSS $1$-line, for which Ravenel's calculation \cite{ravenel2004complexcobordism} gives
$$
\Ext_{BP_*BP}^{1,2w}(BP_*,BP_*)\cong
\begin{cases}
\mathbb Z/p^{1+\nu_p(w)}, & (p-1)\mid w,\\
0, & (p-1)\nmid w.
\end{cases}
$$
It suffices to lift a generator on each cell for which $(p-1)\mid w$.

To construct the lift, apply \eqref{eq:Gn} to the odd-primary series $P_F(z)=1/\lambda_F(z)-1/z$. Since $S(z)=z+O(z^p)$, the integrality argument in the proof of Proposition~\ref{prop:integral} gives an integral cobar cocycle $d^0G_\ell$. The same leading-term calculation as in the proof of Proposition~\ref{prop:divided} shows that $d^0G_\ell$ lifts the algAHSS class $a_w[2\ell]$, where
$$
a_w=\delta_{\mathbb Q}\bigl(\overline{B_w(F)/w}\bigr)
\in \Ext_{BP_*BP}^{1,2w}(BP_*,BP_*).
$$
%
%Put $u=B_w(F)/w$. By \eqref{eq:Gleading},
%$G_\ell(z)=uz^\ell+O(z^{\ell+1})$. Since $S(z)$ has linear term $z$,
%$$
%d^0G_\ell(z)\equiv
%\bigl(\eta_R(u)-\eta_L(u)\bigr)z^\ell \pmod{z^{\ell+1}}.
%$$
%By definition of the coefficient Bockstein, the sphere cobar
%cocycle $\eta_R(u)-\eta_L(u)$ has cohomology class
%$a_w=\delta_{\mathbb Q}(\bar u)$. Thus $d^0G_\ell$ lifts the
%algAHSS class $a_w[2\ell]$.

For $(p-1)\mid w$, $a_w$  has order $p^{1+\nu_p(w)}$ by Miller's denominator calculation \cite{Miller1982Bernoulli} and specialization to a $p$-typical coordinate on the multiplicative formal group \cite{ravenel2004complexcobordism}. Thus $a_w$ generates this group.

It follows that no differential hits $b[2m]$. Since $b[2m]$ lies on the top cell, it also supports no differential. It follows that $q^*b\ne0$ for every $b\ne0$, proving injectivity.

\smallskip
\noindent\textup{(2)} By \textup{(1)}, $q^*b$ is nonzero. It supports no differential by naturality, and cannot be hit for degree reasons. Thus $q^*b$ survives. By naturality it detects
$q^*\beta\neq 0$.

\end{proof}

We now return to $p=3$ and construct sphere classes detected by $\Xi_r$.

\begin{prop}\label{prop:anss3-corrected-family}
For a fixed choice of sign $c\in\{1,-1\}$ and every $r\geq1$, there is a class $\xi_r\in\pi_{144r-38}S^0$ detected by the nonzero permanent cycle
\begin{equation*}
\Xi_r=\beta_{9r/9}+cr\beta_{9r-2},\qquad r\geq1.
\end{equation*}
\end{prop}
\begin{proof}
Put $V=S/(3,v_1^9)$. By \cite{belmont2023beta}, there is a fixed sign $c\in\{1,-1\}$ for which
$$
w=v_2^9+cv_1^8v_2^7
$$
is a permanent cycle in the $3$-primary Adams--Novikov spectral sequence for $V$. Since $V$ is a homotopy ring \cite{oka1979ring} multiplication by a homotopy class detected by $w$ gives a self-map
$$
Q:\Sigma^{144}V\longrightarrow V
$$
whose action on $BP_*V$ is multiplication by $w$ \cite{belmont2023beta}.

Let $i:S^0\to V$ and $q:V\to S^{38}$ be the bottom inclusion and top quotient. Define $\xi_r\in\pi_{144r-38}S^0$ by desuspending the composite
$$
S^{144r}
\xrightarrow{\Sigma^{144r}i}
\Sigma^{144r}V
\xrightarrow{Q^r}
V
\xrightarrow{q}
S^{38}.
$$
The $BP$-Hurewicz image of $Q^r\circ\Sigma^{144r}i$ is 
$$
w^r=v_2^{9r}+crv_1^8v_2^{9r-2}
\quad\text{in }BP_*V=BP_*/(3,v_1^9).
$$

 Consider the cofiber sequences
$$
\begin{gathered}
\Sigma^{36}S/3\xrightarrow{v_1^9}S/3
\longrightarrow V\xrightarrow{\partial_9}\Sigma^{37}S/3,\\
S^0\xrightarrow{3}S^0\longrightarrow S/3
\xrightarrow{\partial}S^1.
\end{gathered}
$$
We have $q=(\Sigma^{37}\partial)\circ\partial_9$. We use the same symbols for the corresponding connecting homomorphisms in $\Ext$.

By the ANSS 2-line computation \cite{ravenel2004complexcobordism} and naturality of the connecting homomorphisms \cite{belmont2023beta},
$$
\begin{aligned}
\partial\partial_9(w^r)
&=\partial\partial_9(v_2^{9r})
  +cr\,\partial\partial_9(v_1^8v_2^{9r-2})\\
&=\beta_{9r/9}+cr\beta_{9r-2}\\
&=\Xi_r\neq 0.
\end{aligned}
$$

By the geometric boundary theorem \cite{ravenel2004complexcobordism}, these connecting homomorphisms commute with ANSS differentials and agree with the geometric boundary maps on associated graded homotopy. Since $w^r$ is a permanent cycle, $\Xi_r$ supports no differentials. There are no incoming differentials for degree reasons. Thus $\Xi_r$ survives and detects $\xi_r$.

\end{proof}

\begin{lemma}\label{lem:anss3-retraction}
For $m\geq1$, the bottom sphere of $(\Sigma^{-2m}\CP_m^{m+s})$ splits off for $s=1$ with arbitrary $m$, for $s=3$ if $3\mid m$, and for $s=5$ if $9\mid m$. In each case there is a retraction
\[
\rho:(\Sigma^{-2m}\CP_m^{m+s})\longrightarrow S^0.
\]
\end{lemma}
\begin{proof}
Mosher's $S$-coreducibility theorem \cite{MosherProjective1968} states that $T(mH_k)=\CP_m^{m+k-1}$ is $S$-coreducible if and only if $M_k\mid m$, where $H_k$ is the canonical complex line bundle over $\CP^{k-1}$, $T(-)$ is the Thom space construction and $M_k$ is the $k$th Atiyah--Todd number. Here $S$-coreducibility means that the bottom sphere splits off.

Passing to the $3$-local setting, the bottom sphere splits off whenever $3^{\nu_3(M_k)}\mid m$. Taking $k=s+1$, the Atiyah--Todd formula \cite{AdamsWalker1965stiefel} gives
$$
\nu_3(M_2)=0,\qquad
\nu_3(M_4)=1,\qquad
\nu_3(M_6)=2.
$$

These values give the three stated conditions.
\end{proof}

\begin{cor}\label{cor:anss3-extension}
Let $m\geq3$, and let $s$ satisfy the hypotheses of Lemma~\ref{lem:anss3-retraction}. If $\beta\in\pi_{2m}S^0$ satisfies $q^*\beta\ne0$, then $\CP^n$ admits an exotic smooth structure for every $m\leq n\leq m+s$.
\end{cor}

\begin{proof}
Let $N=m+s$ and let $\rho$ be the retraction in Lemma~\ref{lem:anss3-retraction}. Consider the maps
$$\CP^m\xrightarrow{i} \CP^N\xrightarrow{Q} \CP^N_m \xrightarrow{\Sigma^{2m}\rho} S^{2m} \xrightarrow{\beta} S^0. $$
Denote the composite $\beta\circ\Sigma^{2m}\rho\circ Q$ by $a.$ The composite $a|_{\CP^m}:= a\circ i=q^*\beta\ne0$. Thus $a|_{\CP^n}\ne0$ for every $m\leq n\leq N$. Theorem~\ref{oddAHSSsmoothing} gives the claimed exotic smooth structures.
\end{proof}

\begin{proof}[Proof of Theorem~\ref{thm:anss3-main}]
For each class listed in Table~\ref{tab:anss3-detectors}, choose a homotopy class $\alpha\in\pi_{2m}S^0$ detected by it. Proposition~\ref{prop:anss3-top-cell} gives $q^*\alpha\ne0$.

The listed values of $m$ and $s$ in Table~\ref{tab:anss3-detectors} satisfy the hypotheses of Lemma~\ref{lem:anss3-retraction}. Therefore Corollary~\ref{cor:anss3-extension} gives exotic smooth structures for dimensions in the interval $[m,m+s]$.
\end{proof}

\section{Density one at $p\geq 5$}
\label{sec:density-one}

In this section we consider the primes $p\geq5$. We use divided beta families to detect exotic smooth structures of $\CP^n$. We will show that complex dimensions with exotic smooth structures form a set of density one. Again, we implicitly localize at $p$.

\begin{theorem}\label{thm:density-one}
Let $\mathcal D_{\mathrm{ex}}$ be the set of positive integers $n\geq3$ for which $\CP^n$ admits an exotic smooth structure.
\begin{enumerate}[label=(\roman*),ref=\roman*,leftmargin=*]
\item\label{part:density-construction} Consider a prime $p\geq5$. Suppose
$$
s\geq1,\qquad p\nmid s,\qquad 3\leq j\leq p,
\qquad(s,j)\neq(1,p).
$$
Let $m=sp(p^2-1)-(p-1)j-1$. Every complex dimension in $[m,m+p-2]$ belongs to $\mathcal D_{\mathrm{ex}}$.
\item\label{part:density-density} Let $\mathcal D_\beta$ be the union of the intervals in part~\ref{part:density-construction}. Then
$$
d(\mathcal D_\beta)=d(\mathcal D_{\mathrm{ex}})=1,
$$
where $d$ denotes natural density:
$$
d(A)=\lim_{X\to\infty}\frac{\#(A\cap[1,X])}{X},
\qquad X\in\bbN.
$$
\end{enumerate}
\end{theorem}

\subsection{Construction at a fixed prime}
\label{subsec:density-geometric}

By Theorem~\ref{oddAHSSsmoothing}, it suffices to construct nonzero stable cohomotopy classes on the projective spaces in question. We use Proposition~\ref{prop:anss3-top-cell} to obtain nonzero pullbacks of divided beta-family elements, then extend these pullbacks over short intervals of dimensions.

\begin{lemma}\label{lem:density-sphere-family}
Let $p\geq5$ be prime, and suppose
$$
s\geq1,\qquad p\nmid s,\qquad 3\leq j\leq p,
\qquad (s,j)\neq(1,p).
$$
Put $m=sp(p^2-1)-(p-1)j-1$. There is a class $x_{s,j}\in\pi_{2m}S^0$ detected by the nonzero ANSS class $\beta_{sp/j}$.
\end{lemma}

\begin{proof}
Under these hypotheses, $\beta_{sp/j}$ is a nonzero permanent cycle by \cite{ravenel2004complexcobordism}. Choose a homotopy class $x_{s,j}$ that it detects.
\end{proof}

% \begin{lemma}\label{lem:density-sphere-family}
% Let $p\geq5$ be prime, and suppose
% $$
% s\geq1,\quad p\nmid s,\quad3\leq j\leq p,
% \quad(s,j)\neq(1,p),
% $$
% and put $m=sp(p^2-1)-(p-1)j-1$.
% There is a sphere class $x_{s,j}\in\pi_{2m}S^0$ detected by $\beta_{sp/j}$, such that
% $\alpha_1x_{s,j}\neq0$. \end{lemma}

% \begin{proof}
% By \cite{ravenel2004complexcobordism}, 
% the class
% $\beta_{sp/j}$ is a nonzero permanent cycle.
% Choose a homotopy class $x_{s,j}$ detected by $\beta_{sp/j}$.

% Under the assumption on $s,j$, Miller--Ravenel--Wilson's product criterion gives a nonzero
% algebraic class $\alpha_1\beta_{sp/j}$
% \cite{Miller1977periodic}.

% The ANSS is concentrated in internal degrees that are multiples of $2(p-1)$. 
% Therefore nonzero differentials have length $\geq 2p-1 \geq 9$.
% Thus by degree reasons, there are no incoming differentials hitting $\alpha_1\beta_{sp/j}$ which is in ANSS filtration 3. 
% Outgoing differentials on the product
% vanish by the Leibniz rule. Thus $\alpha_1\beta_{sp/j}$ detects the nonzero product $\alpha_1x_{s,j}$.

% \end{proof}

The following lemma combines the nonvanishing result of Proposition~\ref{prop:anss3-top-cell} with a splitting of short stunted projective spectra.

\begin{lemma}\label{lem:density-projective-detection}
Let $p\geq5$ be prime and $m\geq1$. Suppose that $x\in\pi_{2m}S^0$ is detected by a nonzero class on the ANSS $2$-line. Then $q^*x\in[\CP^m,S^0]$ is nonzero and extends to a nonzero class in $[\CP^{m+e},S^0]$ for every $0\leq e\leq p-2$.
\end{lemma}

\begin{proof}
Proposition~\ref{prop:anss3-top-cell} gives $q^*x\neq0$.

For $0\leq e\leq p-2$, there is a splitting \cite{Beben2010Stiefel}
$$
\CP_m^{m+e}\simeq\bigvee_{i=0}^{e}S^{2m+2i}.
$$
Let $\rho:\CP_m^{m+e}\to S^{2m}$ be a retraction onto the bottom cell. The composite
$$
\CP^{m+e}\longrightarrow\CP_m^{m+e}
\xrightarrow{\rho}S^{2m}\xrightarrow{x}S^0
$$
restricts to $q^*x$ on $\CP^m$. This gives the required nonzero extension.
\end{proof}

% \begin{lemma}\label{lem:density-projective-detection}
% Consider prime $p\geq5$ and $m\geq1$ with $p\nmid m$.
% Suppose $x\in\pi_{2m}S^0$ satisfies $\alpha_1x\neq0$.
% Then $q_m^*x\in  [\CP^{m}, S^0]$ is nonzero and extends to a class in
% $[\CP^{m+e}, S^0]$ for every $0\leq e\leq p-2$.
% \end{lemma}

% \begin{proof}

% Let $h=p-1$. Choose $a\equiv m \pmod {p-1}.$
% The natural Mimura--Nishida--Toda splitting gives
% $$
% \CP^n\simeq\bigvee_{a=1}^{h}X_a(n),
% $$
% where $X_a(n)$ contains the cells of dimension $2r$ with
% $r\equiv a\pmod h$, compatibly with skeletal inclusions
% \cite{Beben2010Stiefel}.

% Consider 
% $$f_m: S^{2(m+h)-1}\xrightarrow{\gamma_m} X_a(m)\xrightarrow{q} S^{2m},$$
% where $\gamma_m$ is top cell attaching map in $X_a(m+h)$ and $q$ is collapsing to the top cell.
% Therefore
% $
% f_m\in\pi_{2p-3}S^0
% \cong\mathbb Z/p\{\alpha_1\}.
% $

% The cofiber of $f_m$ is the two-cell quotient of $X_a(m+h)$ from the top.
% Let $u$ and $v$ denote the cohomology generators in degrees $2m$ and
% $2(m+h)$. Naturality and the Cartan formula give
% $$
% \mathcal P^1u=mv,
% \qquad
% \mathcal P^1(z^m)=mz^{m+p-1}.
% $$
% Since $p\nmid m$, this operation is nonzero. Therefore
% $f_m=\lambda\alpha_1$ for some
% $\lambda\in\mathbb Z^\times$, and
% $$
% (xq)\gamma_m
% =xf_m
% =\lambda\alpha_1x
% \neq0.
% $$
% Thus $xq\neq0$. By compatibility of the splitting with the
% top-cell quotient, extending $xq$ by zero on the other summands
% gives $q_m^*x\neq0$.

% For $0\leq e\leq p-2$, inclusion induces an equivalence
% $
% i_e:X_a(m)\xrightarrow{\simeq}X_a(m+e),
% $ and the same arguments follow.

% \end{proof}

\begin{proof}[Proof of Theorem~\ref{thm:density-one}\,
(\ref{part:density-construction})] The result follows by combining Theorem~\ref{oddAHSSsmoothing}, Lemma~\ref{lem:density-sphere-family} and Lemma~\ref{lem:density-projective-detection}.
\end{proof}

\subsection{Density one counting}
\label{subsec:density-counting}

In this subsection, we use $[a,b]$ to denote the set of all integers $a\leq n\leq b.$ This section uses only elementary arithmetic arguments, not homotopy theory.

Fix a prime $p\geq5$ and an integer $s$. Then we have
$$\cup_{j\in [3,p]}[m,m+p-2]= [sp(p^2-1)-p(p-1)-1, sp(p^2-1)-2p].$$
Here we temporarily include the excluded pair $(s,j)=(1,p)$. Allowing also negative $s$, define
$$
B_p=\bigcup_{\substack{s\in\bbZ\\p\nmid s}}
[sp(p^2-1)-p(p-1)-1,\ sp(p^2-1)-2p].
$$
Translation by $p^2(p^2-1)$ replaces $s$ by $s+p$, so $p^2(p^2-1)$ is a period of $B_p$.

For an integer $n$, we write
$$
n+p(p-1)+1=(p^2-1)I_p(n)+r_p(n),\qquad0\leq r_p(n)<p^2-1.
$$

\begin{lemma}\label{lem:density-membership}
For every integer $n$,
$$
n\in B_p\quad\iff \quad
r_p(n)<(p-1)(p-2),\quad p\mid I_p(n),\quad p^2\nmid I_p(n).
$$
\end{lemma}

\begin{proof}
The result can be proved by direct computation.
\end{proof}

By Lemma~\ref{lem:density-membership}, $B_p$ has density 
$$d(B_p) = \frac{p-1}{p^2}\cdot \frac{(p-1)(p-2)}{p^2-1}= \frac{(p-1)(p-2)}{p^2(p+1)}$$

To count the union over several primes, we must also account for overlaps between these periodic sets.

Enumerate the primes $p_1=5<p_2=7<p_3=11<\cdots$ and write
$$
\begin{aligned}
\mathcal B_k&=\bigcup_{i=1}^k B_{p_i},&
L_k&=\operatorname{lcm}_{1\leq i\leq k}\bigl(p_i^2(p_i^2-1)\bigr),\\
\delta_k&=\frac{\#(\mathcal B_k\cap[0,L_k))}{L_k}.
\end{aligned}
$$
We define $L_0=1$, $\mathcal B_0=\varnothing$. The finite union $\mathcal B_k$ is periodic with period $L_k$, so $\delta_k=d(\mathcal B_k)$.

\begin{prop}\label{prop:density-finite-unions}
We have
$$
\delta_{k}-\delta_{k-1}\geq
\frac{p_k-1}{p_k^2}(\frac{p_k-2}{p_k+1} - \delta_{k-1}).
$$
\end{prop}

\begin{proof}
At stage $k$, put $p=p_k$ and $ M_k=\operatorname{lcm}(L_{k-1},p^2-1). $ It follows that
$$
L_k=M_kp^2.
$$
Arrange the integers in one $L_k$ period into
$$
\Omega_a=\{a+M_kt:0\leq t<p^2\},\qquad0\leq a<M_k.
$$
Since both $L_{k-1}$ and $p^2-1$ divide $M_k$, membership in the old union $\mathcal B_{k-1}$ and the value of $r_p$ are constant on each row. And we have
$$
I_p(a+M_kt)=I_p(a)+\frac{M_k}{p^2-1}\,t.
$$
The integer $M_k/(p^2-1)$ is prime to $p$, so $\Omega_a$ contains $p-1$ members of $B_p$ if $r_p(a)<(p-1)(p-2)$, and none otherwise.

There are $(1-\delta_{k-1})M_k$ rows that do not belong to $\mathcal B_{k-1}$. Since $p^2-1$ divides $M_k$, $\frac{3}{p+1}M_k$ rows have $r_p(a)\geq(p-1)(p-2)$. Therefore at least $(1-\delta_{k-1}-\frac{3}{p+1})M_k$ rows are not in $\mathcal B_{k-1}$ and satisfy $r_p(a)<(p-1)(p-2)$. Each such row contributes $p-1$ new dimensions to $\mathcal B_{k}$, so
$$
L_k(\delta_{k}-\delta_{k-1})\geq(p-1)(1-\delta_{k-1}-\frac{3}{p+1})M_k.
$$
The result follows.

\end{proof}

\begin{proof}[Proof of Theorem~\ref{thm:density-one}\,(\ref{part:density-density})]

By Proposition~\ref{prop:density-finite-unions}, $\delta_k\to 1.$

We need to exclude the negative part and the case when $(s,j)=(1,p)$.

When $s\leq 0$, the dimensions in $[m,m+p-2]$ are also negative. When $(s,j)=(1,p)$, the interval $N_p:= [m,m+p-2]$ is finite.

Therefore 
$$d(\mathcal D_\beta)\geq d((\mathcal B_k\cap \mathbb Z_+)\setminus (\cup_{i=1}^k N_{p_i}))=\delta_k,$$ and thus $d(\mathcal D_\beta)=d(\mathcal D_{\mathrm{ex}})=1.$

\end{proof}

\newpage

\begingroup
\raggedright
\bibliography{refs}
\bibliographystyle{alpha}
\endgroup

\end{document}